\documentclass[11pt]{article}
\usepackage [margin= 0.85in]{geometry}

\usepackage{amsmath}
\usepackage{wasysym}
\usepackage{latexsym}
\usepackage{amsfonts}
\usepackage{mathrsfs}
\usepackage{amssymb}
\usepackage{ifsym}
\usepackage{dsfont}
\usepackage[all]{xy}
\usepackage{authblk}
\usepackage{lipsum}

\newcommand\blfootnote[1]{%
  \begingroup
  \renewcommand\thefootnote{}\footnote{#1}%
  \addtocounter{footnote}{-1}%
  \endgroup
}

\usepackage{amsfonts}
\usepackage{amssymb}

\begin{document}

\title{Set theory with a proper class of indiscernibles}
\author{Ali Enayat

\textsc{\bf Dedicated to the memory of Ken Kunen}\\}
\maketitle

\begin{abstract}
We investigate an extension of $\mathrm{ZFC}$ set theory, denoted $\mathrm{%
ZFI}_{\mathrm{<}}\mathrm{,}$ which is equipped with a well-ordering $<$ of
the universe \textrm{V} of set theory, and a proper class $\mathrm{I}$ of
indiscernibles over the structure $(\mathrm{V},\in ,<)$. \medskip

Our main results are Theorems A, B, and C below. Note that the equivalence
of condition $(ii)$ and $(iii)$ in Theorem A was established in an earlier
(2004) published work of the author. In what follows \textrm{GBC} is the G%
\"{o}del-Bernays theory of classes with global choice. In Theorem C the
symbol $\rightarrow $ is the usual Erd\H{o}s-arrow notation for partition
calculus.\medskip

\noindent \textbf{Theorem A.}~\textit{The following are equivalent for a
sentence }$\varphi $\textit{\ in the language }$\{=,\in \}$ \textit{of set
theory}:\medskip

\noindent $(i)$ $\mathrm{ZFI}_{\mathrm{<}}\vdash \varphi .$\medskip

\noindent $(ii)$\textbf{\ }$\mathrm{ZFC}+\Lambda \vdash \varphi ,$ \textit{%
where} $\Lambda =\{\lambda _{n}:n\in \omega \}$, \textit{and} $\lambda _{n}$
\textit{is the sentence asserting the existence of an }$n$\textit{-Mahlo
cardinal} $\kappa $ \textit{such that} $\mathrm{V}(\kappa )$ \textit{is a} $%
\Sigma _{n}$\textit{-elementary submodel of the universe }$\mathrm{V}$%
\textit{.}\medskip

\noindent $(iii)$ \textrm{GBC} $+$ \textit{\textquotedblleft }\textrm{Ord is
weakly compact}\textit{\textquotedblright } $\vdash \varphi $.\medskip

\noindent \textbf{Theorem B.}~\textit{Every }$\omega $-\textit{model of} $%
\mathrm{ZFI}_{\mathrm{<}}$\textit{\ satisfies }$\mathrm{V}\neq \mathrm{L}$%
\textit{.}\medskip

\noindent \textbf{Theorem C.}~\textit{The sentence expressing }$\forall
m,n\in \omega \left( \mathrm{Ord}\rightarrow \left( \mathrm{Ord}\right)
_{m}^{n}\right) $ is \textit{not provable in the theory} $T=$ \textrm{GBC} $%
+ $ \textit{\textquotedblleft }\textrm{Ord is weakly compact}\textit{%
\textquotedblright , assuming }$T$\textit{\ is consistent.}

\medskip

\noindent The paper also includes results about the interpretability relationship
between the theories $\mathrm{ZFC}+\Lambda $, $\mathrm{ZFI}_{\mathrm{<}}$,
and GBC + ``Ord is weakly compact''.

\medskip
\end{abstract}

\blfootnote{\textit{Key Words}. Zermelo-Fraenkel set theory, G\"{o}del-Bernays class theory, indiscernibles, Mahlo cardinal, weakly compact cardinal, satisfaction class.}

\blfootnote {\textit{2010 Mathematical Subject Classification}. Primary: 03E55, 03F25, 03C62;
Secondary: 03E02, 03H15.}

\begin{center}
{\small TABLE\ OF\ CONTENTS}

{\small 1.
INTRODUCTION.......................................................................................................................2%
}

{\small 2.
PRELIMINARIES.......................................................................................................................3%
}

{\small \quad 2.1. Models of set
theory..............................................................................................................3%
}

{\small \quad 2.2. Satisfaction
classes................................................................................................................5%
}

{\small \quad 2.3.
Indiscernibles........................................................................................................................6%
}

{\small \quad 2.4. The theory GBC + \textquotedblleft Ord is weakly
compact\textquotedblright
......................................................................8}

{\small 3. THE BASIC FEATURES OF ZFI AND ZFI}$_{<}${\small %
..........................................................................14}

{\small 4. WHAT ZFI}$_{<}${\small \ KNOWS ABOUT SET
THEORY.......................................................................20%
}

{\small 5. INTERPRETABILITY ANALYSIS OF ZFI}$_{<}${\small %
...........................................................................27%
}

{\small 6. SOME VARIANTS OF ZFI}$_{<}${\small %
....................................................................................................30%
}

{\small 7. OPEN
QUESTIONS..................................................................................................................33%
}
\end{center}

\pagebreak

\begin{center}
\textbf{1. INTRODUCTION\bigskip }
\end{center}

\noindent The principal focus of this paper is on an extension $\mathrm{ZFI}%
_{\mathrm{<}}$ of Zermelo-Fraenkel set theory $\mathrm{ZF}$ that is equipped
with a global well-ordering $<$ and a proper class $I$ of ordinals such that
$\left( I,\in \right) $ is a collection of order indiscernibles over the
structure $(\mathrm{V},\in ,<).$ Moreover, the axioms of $\mathrm{ZFI}_{%
\mathrm{<}}$ stipulate that the expanded universe $(\mathrm{V},\in ,<,I)$
satisfies the axioms of $\mathrm{ZF}$ in the extended language incorporating
$<$ and $I$. Thus $\mathrm{ZFI}_{\mathrm{<}}$ is a system of set theory that
can be described as strongly `anti-Leibnizian': The Leibniz dictum on the
identity of indiscernibles bars the existence of a single pair of distinct
indiscernibles in the universe $(\mathrm{V},\in )$ of sets, but models of $%
\mathrm{ZFI}_{\mathrm{<}}$ are endowed, intuitively speaking, with an
unnameable number of such objects that are grouped into a proper class $I$\
that can be used in set-theoretical reasoning.\footnote{%
The impact of Leibnizian motifs in set theory and its model theory is
explored in \cite{Ali-LM} and \cite{Ali-Leibnizian}.}\medskip

The precise definition of $\mathrm{ZFI}_{\mathrm{<}}$ is given in Section 3.
The definition makes it clear that (a) if $\kappa $ is a weakly compact
cardinal, then $\left( \mathrm{V}(\kappa ),\in \right) $ has an expansion
that satisfies any prescribed finitely axiomatized subtheory of $\mathrm{ZFI}%
_{\mathrm{<}}$, and (b) if $\kappa $ is a Ramsey cardinal, then $\left(
\mathrm{V}(\kappa ),\in \right) $ has an expansion to a model of $\mathrm{ZFI%
}_{\mathrm{<}}.$ One of our main results is Theorem 4.1 (a refinement of
Theorem A of the abstract) that shows that the purely set-theoretical
consequences of $\mathrm{ZFI}_{\mathrm{<}}$ coincide with the theorems of
the theory obtained by augmenting ZFC with the Levy scheme\footnote{%
The Levy scheme $\Lambda $ was denoted $\Phi $ in earlier work of the
author, and in particular in \cite{Ali NFUA}. The new notation is occasioned
by the author's appreciation of the role played by Azriel Levy in the
investigations of the Mahlo hierarchy and reflection phenomena, masterfully
overviewed in Kanamori's portraiture \cite{Kanamori on Levy}. In Subsection
2.4 we review the basic features of the Levy\ Scheme.} $\Lambda $, a scheme
that ensures that the class of ordinals behaves like an $\omega $-Mahlo
cardinal (the precise definition of $\Lambda $ is given in Definition
2.4.10). Theorem 4.1 complements the main results in \cite{Ali powerlike}
and \cite{Ali NFUA} that exhibit the surprising ways in which $\mathrm{ZFC}%
+\Lambda $ manifests itself as a canonical theory, especially in the context
where the model theory of $\mathrm{ZF}$ is compared with the model theory of
$\mathrm{PA}$ (Peano Arithmetic). In contrast, parts (e) and (f) of Theorem 3.8 (which refine
Theorem B of the abstract) show that an $\omega $-model of $\mathrm{ZFI}%
_{\mathrm{<}}$ (i.e., a model of $\mathrm{ZFI}_{\mathrm{<}}$ whose $\omega $
is well-founded in the real world) satisfies large cardinal hypotheses
significantly stronger than the existence of $\omega $-Mahlo cardinals. Our
third main result is Theorem 4.9 (Theorem C of the abstract), which should
be contrasted with the fact that \textrm{GBC} $+$ \textquotedblleft \textrm{%
Ord is weakly compact}\textquotedblright\ can prove sentences of the form $%
\forall \kappa \left( \mathrm{Ord}\rightarrow \left( \mathrm{Ord}\right)
_{\kappa }^{n}\right) ,$ where $n$ ranges over nonzero natural numbers in
the real world. We also include some interpretability-theoretic results
concerning $\mathrm{ZFI}_{\mathrm{<}}$ and variants of $\mathrm{ZFI}_{\mathrm{<}}$. \medskip

There is a notable series of papers investigating combinatorial features of $%
n$-Mahlo cardinals, beginning with the groundbreaking work of Schmerl \cite%
{Schmerl}, which eventually culminated in the Hajnal-Kanamori-Shelah paper
\cite{Hajnal-Kanamori-Shelah}. The relationship between $n$-Mahlo cardinals
and various types of sets of indiscernibles has also been extensively
studied by many researchers including McAloon, Ressayre, Friedman, Finkel
and Todor\v{c}evi\'{c} (see, e.g., \cite{Finkel-Ressayre} and \cite%
{Finkel-Stevo}). However, the proofs of our results dominantly employ
techniques from the model theory of set theory together with classical
combinatorial ideas, thus they do not rely on the machinery developed in the
above body of work. Of course it would be interesting to work out the
relationship between our results and the aforementioned literature.\medskip

The organization of the paper is as follows. Section 2 contains a mix of
preliminary material employed in the paper; the reader is advised to pay
special attention to Subsections 2.3 and 2.4. Section 3 introduces $\mathrm{%
ZFI}$ and $\mathrm{ZFI}_{\mathrm{<}}$ and mostly focuses on their model
theory. Section 4 is devoted to the calibration of the purely
set-theoretical consequences of $\mathrm{ZFI}_{\mathrm{<}}$, and Section 5
studies $\mathrm{ZFI}_{\mathrm{<}}$ from an interpretability-theoretic point
of view. In Section 6 we discuss four systems that are closely related to $%
\mathrm{ZFI}_{<}$. Finally, we close the paper by presenting a few open questions in
Section 7.\medskip

\textbf{History and} \textbf{Acknowledgments.}~This paper might appear as a
natural sequel to \cite{Ali NFUA}, but the work reported
here arose in a highly indirect way as a result of an engagement with certain
potent ideas proposed by Jan Mycielski \cite{Mycielski JSL} concerning
Leibnizian motifs in set theory, an engagement that culminated in the
trilogy of papers \cite{Ali-LM}, \cite{Ali-Leibnizian}, and \cite{Enayat-DO}%
. Informed by Bohr's aphorism \textquotedblleft The opposite of a correct
statement is a false statement. But the opposite of a profound truth may
well be another profound truth\textquotedblright , and as if to maintain a
cognitive balance, upon the completion of the aforementioned trilogy my
attention and curiosity took an opposite turn towards the highly
`anti-Leibnizian' systems of set theory studied here. The protoforms of the
results of this paper were first presented at the New York Logic Conference
(2005), IPM Logic Conference (2007, Tehran, Iran), the Kunen Fest
Meeting (2009, Madison, Wisconsin), and most recently at the Oxford Set
Theory Seminar (2020). \medskip

I am grateful to Kentaro Fujimoto, Philip Welch, Kentaro Sato, and an anonymous reader for
reading an earlier draft of this paper and offering their detailed comments
and suggestions for improvements. Thanks also to Neil Barton, Andreas Blass,
Cezary Cie\'{s}li\'{n}ski, Vika Gitman, Joel David Hamkins, Roman Kossak,
Mateusz \L e\l yk, Jim Schmerl, and Bartosz Wcis\l o for their keen interest
in this work.

\medskip
The research presented in this paper was supported by the National Science Centre, Poland (NCN), grant number 2019/34/A/HS1/00399.
\textbf{\bigskip }

\begin{center}
\textbf{2.~PRELIMINARIES\bigskip }
\end{center}

In this section we collect the basic definitions, notations, conventions,
and results that will be used in the remaining sections.
\medskip

\begin{center}
\textbf{2.1.~Models of set theory}\medskip
\end{center}

\noindent \textbf{2.1.1.}~\textbf{Definitions and basic facts.~(}Models,
languages, and theories) Let $\mathcal{L}_{\mathrm{Set}}=\{=,\in \}$ be the usual language of set theory. In what follows we make the blanket assumption that $\mathcal{M}$, $\mathcal{N}$, etc. are $\mathcal{L}$-structures, where $\mathcal{L} \supseteq \mathcal{L}_{\mathrm{Set}}$. By \textit{a model of set theory}, we mean an $\mathcal{L}$-structure that satisfies enough of \textrm{ZF} set theory to support a decent theory of ordinals, and of
the von Neumann levels $\mathrm{V}(\alpha )$ of the universe $\mathrm{V}$ of $\textrm{ZF}$.

\medskip

\noindent \textbf{(a) }We follow the convention of using $M$, $M^{\ast },$ $%
M_{0}$, etc.~to denote (respectively) the universes of discourse of
structures $\mathcal{M}$, $\mathcal{M}^{\ast },$ $\mathcal{M}_{0},$ etc.
Given a structure $\mathcal{M}$, we write $\mathcal{L}(\mathcal{M})$ for the
\textit{language of} $\mathcal{M}$. Given some relation symbol $R\in
\mathcal{L}(\mathcal{M})$, we often write $R_{M}$ for the $\mathcal{M}$%
-interpretation of $R$. In particular, we denote the membership relation of $%
\mathcal{M}$ by $\in _{M}$; thus an $\mathcal{L}_{\mathrm{Set}}$-structure $\mathcal{M}$ is of the form $(M,\in_{M})$. Sometimes when there is no risk of
confusion, we conflate formal symbols with their denotations.\medskip

\noindent \textbf{(b) }For $c\in M$, $\mathrm{Ext}_{\mathcal{M}}(c)$ is the $%
\mathcal{M}$-extension of $c$, i.e.,

\begin{center}
$\mathrm{Ext}_{\mathcal{M}}(c):=\{m\in M:m\in _{M}c\}.$
\end{center}

\noindent We say that a subset $X$ of $M$ \textit{is coded in} $\mathcal{M}$
if there is some $c\in M$ such that $\mathrm{Ext}_{\mathcal{M}}(c)=X.$ $X$\
is \textit{piecewise coded} in $\mathcal{M}$ if $X\cap \mathrm{Ext}_{%
\mathcal{M}}(m)$ is coded for each $m\in M.$ For $A\subseteq M$, $\mathrm{Cod%
}_{A}(\mathcal{M)}$ is the collection of sets of the form $A\cap \mathrm{Ext}%
_{\mathcal{M}}(c)$, where $c\in M.$\medskip

\noindent \textbf{(c) }$\mathrm{Ord}^{\mathcal{M}}$ is the class of
\textquotedblleft ordinals\textquotedblright\ of $\mathcal{M}$, i.e., $%
\mathrm{Ord}^{\mathcal{M}}:=\left\{ m\in M:\mathcal{M}\models \mathrm{Ord}%
(m)\right\} ,$ where $\mathrm{Ord}(x)$ expresses \textquotedblleft $x$ is
transitive and is well-ordered by $\in $\textquotedblright . More generally,
for a formula $\varphi (\overline{x})$, where $\overline{x}=\left(
x_{1},\cdot \cdot \cdot ,x_{k}\right) $, we write $\varphi ^{\mathcal{M}}$
for $\left\{ \overline{m}\in M^{k}:\mathcal{M\models \varphi }\left(
m_{1},\cdot \cdot \cdot ,m_{k}\right) \right\} .$ We write $\mathbb{\omega }%
^{\mathcal{M}}$ for the set of finite ordinals (i.e., natural numbers) of $%
\mathcal{M}$, and $\mathbb{\omega }$ for the set of finite ordinals in the
real world, whose members we refer to as\textit{\ metatheoretic natural
numbers. }$\mathcal{M}$ is said to be\textit{\ }$\omega $-\textit{standard}
if $\left( \mathbb{\omega },\mathbb{\in }\right) ^{\mathcal{M}}\cong \left(
\mathbb{\omega },\mathbb{\in }\right) .$ For $\alpha \in \mathrm{Ord}^{%
\mathcal{M}}$ we often use $\mathcal{M}(\alpha )$ to denote the substructure of $\mathcal{M}$ whose universe is $\left( \mathrm{V}(\alpha )\right) ^{\mathcal{M}}.$\medskip

\noindent \textbf{(d)} $\mathcal{N}$ is said to \textit{end extend} $%
\mathcal{M}$ (equivalently: $\mathcal{M}$ is an \textit{initial} submodel of
$\mathcal{N}$), written $\mathcal{M}\subseteq _{\mathrm{end}}\mathcal{N},$
if $\mathcal{M}$ is a submodel of $\mathcal{N}$ and for every $a\in M,%
\mathrm{Ext}_{\mathcal{M}}(a)=\mathrm{Ext}_{\mathcal{N}}(a)$. We often write
\textquotedblleft e.e.e.\textquotedblright\ instead of \textquotedblleft
elementary end extension\textquotedblright . It is easy to see that if $%
\mathcal{N}$ is an e.e.e.~of a model $\mathcal{M}$ of $\mathrm{ZF}$, then $%
\mathcal{N}$ is a \textit{rank extension} of $\mathcal{M}$, i.e., whenever $%
a\in M$ and $b\in N\backslash M$, then $\mathcal{N}\models \rho (a)\in \rho
(b)$, where $\rho $ is the usual ordinal-valued rank function defined by $%
\rho (x)=\sup \{\rho (y)+1:y\in x\}.$\medskip

\noindent \textbf{(e)} We treat $\mathrm{ZF}$ as being axiomatized as usual,
except that instead of including the scheme of replacement among the axioms
of \textrm{ZF}, we include the schemes of separation and collection, as in
\cite[Appendix A]{Chang-Keisler}. Thus, in our set-up the axioms of Zermelo
set theory $\mathrm{Z}$ are obtained by removing the scheme of collection
from the axioms of $\mathrm{ZF}$. More generally, we construe $\mathrm{ZF}(%
\mathcal{L})$ to be the natural extension of $\mathrm{ZF}$ in which the
schemes of separation and collection are extended to $\mathcal{L}$-formulae,
and we will denote $\mathrm{Z}(\mathcal{L})$ (Zermelo set theory over $%
\mathcal{L}$) as the result of extending $\mathrm{Z}$ with the $\mathcal{L}$%
-separation scheme $\mathrm{Sep}(\mathcal{L})$, which consists of the
universal closures of $\mathcal{L}$-formulae of the form:

\begin{center}
$\forall v\exists w\forall x(x\in w\longleftrightarrow x\in v\wedge \varphi
(x,\overline{y})).$
\end{center}

\noindent Thus $\mathrm{ZF}(\mathcal{L})$ is the result of augmenting $%
\mathrm{Z}(\mathcal{L})$ with the $\mathcal{L}$-collection scheme $\mathrm{%
Coll}(\mathcal{L})$, which consists of the universal closures of $\mathcal{L}
$-formulae of the form:

\begin{center}
$\left( \forall x\in v\text{ }\exists y\text{\ }\varphi (x,y,\overline{z}%
)\right) \rightarrow \left( \exists w\text{ }\forall x\in v\text{ }\exists
y\in w\text{ }\varphi (x,y,\overline{z})\right) .$
\end{center}

\noindent When $\mathcal{L=L}_{\mathrm{%
Set}}\cup \{X\}$, we will write $\mathrm{Sep}(X)$, $\mathrm{Coll}(X)$, etc.
instead of $\mathrm{Sep}(\mathcal{L})$, $\mathrm{Coll}(\mathcal{L})$, etc.
(respectively).

\medskip

\noindent \textbf{(f)} Suppose $n\in \omega $. $\Sigma _{n}(\mathcal{L})$ is
the natural extension to $\mathcal{L}$-formulae of the usual Levy hierarchy.
Thus $\Sigma _{0}(\mathcal{L})$ is the smallest family of $\mathcal{L}$%
-formulae that contains all atomic $\mathcal{L}$-formulae and is closed
under Boolean operations and bounded quantification. We write $\mathcal{M}%
\prec _{\mathrm{\Sigma }_{n}(\mathcal{L})}\mathcal{N}$ to indicate that $%
\mathcal{M}$ is a proper $\mathrm{\Sigma }_{n}(\mathcal{L})$-elementary
submodel of $\mathcal{N}$, i.e., $\mathcal{M}$ is a proper submodel of $%
\mathcal{N}$, and for each $k$-ary $\varphi (\overline{x})\in \Sigma _{n}(%
\mathcal{L})$ and each $k$-tuple $\overline{m}$ from $\mathcal{M}$, $%
\mathcal{M}\models \varphi (\overline{m})$ iff $\mathcal{N}\models \varphi (%
\overline{m})$.\medskip

\noindent \textbf{(g)} Given a language $\mathcal{L}$ and a predicate symbol
$X$, we often write $\mathcal{L}(X)$ instead of $\mathcal{L}\cup \{X\}.$
Similarly, we write $\Sigma _{n}(X)$ instead of $\Sigma _{n}(\mathcal{L}_{\mathrm{Set}}(X))$%
, and $\mathrm{ZF}(X)$ instead of $\mathrm{ZF}(\mathcal{L_{\mathrm{Set}}(}X\mathcal{)}).$
Given $\mathcal{M}\models \mathrm{ZF}(\mathcal{L})$, we say that a subset $%
X_{M}$ of $M$ is $\mathcal{M}$-\textit{amenable} if $(\mathcal{M}%
,X_{M})\models \mathrm{ZF}(\mathcal{L}(X)).$\footnote{%
We will often conflate $X$ and $X_{M}$ to lighten the notation. Also note
that some authors use the expression `$X$ is a class of $\mathcal{M}$'
instead of `$X$ is amenable over $\mathcal{M}$'.} \ It is well-known that if
$\mathcal{M}\models \mathrm{ZF}$ and $X_{M}\subseteq M$, then $(\mathcal{M}%
,X_{M})\models \mathrm{ZF}(X)$ iff $X_{M}$\ is piecewise coded in $\mathcal{%
M\ }$and $(\mathcal{M},X_{M})\models \mathrm{Coll}(X)$.\footnote{%
This fact is essentially due to Keisler, its proof is implicit in the proof
of Theorem C of \cite{Keisler (tree)}.}\medskip

\noindent \textbf{(h)} Suppose $X\subseteq M^{k},$ for $1\leq k\in \mathbb{%
\omega }$, $X$ is $\mathcal{M}$-\textit{definable} if $X=\varphi ^{\mathcal{M%
}}$ for some $\mathcal{L}(\mathcal{M})$-formula. $X$ is \textit{%
parametrically} $\mathcal{M}$-\textit{definable} if $X=\varphi ^{\mathcal{M}%
^{+}}$ for some $\mathcal{L}(\mathcal{M}^{+})$-formula, where $\mathcal{M}%
^{+}$ is the expansion $\left( \mathcal{M},m\right) _{m\in M}$ of $\mathcal{M%
}$. A \textit{parametrically} $\mathcal{M}$-\textit{definable} function is a
function $f:M^{k}\rightarrow M$ (where $1\leq k\in \mathbb{\omega })$ such
that the graph of $f$ is parametrically $\mathcal{M}$-definable. If $%
\mathcal{M}^{\ast }$ is an elementary extension of $\mathcal{M}$, then any
such $f$ extends naturally to a parametrically $\mathcal{M}^{\ast }$%
-definable function according to the same definition; we may also denote
this extension as $f$.\medskip

\noindent \textbf{(i) }$\mathcal{M}$ has \textit{definable Skolem functions}
if for every $\mathcal{L}(\mathcal{M})$-formula $\varphi (x,y_{1},\ldots
,y_{k})$, whose free variable(s) include a distinguished free variable $x$
and whose other free variables (if any) are $y_{1},\ldots ,y_{k}$, there is
an $\mathcal{M}$-definable function $f$ such that (abusing notation
slightly):

\begin{center}
$\mathcal{M}\models \forall y_{1}\ldots \forall y_{k}\left[ \exists x\
\varphi (x,y_{1},\cdot \cdot \cdot ,y_{k})\rightarrow \varphi \left(f(y_{1},\cdot
\cdot \cdot ,y_{k}),y_{1},\cdot \cdot \cdot ,y_{k} \right) \right] .$
\end{center}

\noindent \textbf{(j)} If $\mathcal{M}$ has definable Skolem functions, then
given any $X\subseteq M$, there is a least elementary substructure $\mathcal{%
M}_{X}$ of $\mathcal{M}$ that contains $X$, whose universe is the set of all
applications of $\mathcal{M}$-definable functions to tuples from $X$. We
will refer to $\mathcal{M}_{X}$ as the \textit{submodel of}\emph{\ }$%
\mathcal{M}$\emph{\ }\textit{generated by} $X$.\medskip

\noindent \textbf{(k)} Given a distinguished binary relation symbol $<$, the
\textit{global well-ordering axiom}, denoted $\mathrm{GW}$ is the
conjunction of the sentences \textquotedblleft $<$ is a linear
order\textquotedblright\ and \textquotedblleft every nonempty set has a $<$%
-least element\textquotedblright . It is well-known that within $\mathrm{ZF}(<)+\mathrm{GW}$ there is global well-ordering $<^{*}$ that is set-like and thus is of order-type $\mathrm{Ord}$ (by defining $ x<^{*} y$ iff [($\rho (x)=\rho (y)$ and $x< y$) or $ \rho (x)\in \rho (y)]$, where $\rho $ is the usual ordinal-valued rank function).
\medskip

\noindent \textbf{(l)} Given a distinguished unary function
symbol $f$, \textit{the global choice axiom}, denoted $\mathrm{GC}$, is the
axiom $\forall x\left( x\neq \varnothing \rightarrow f(x)\in x\right) .$%
\medskip

The following two theorems are well-known. A proof of Theorem 2.1.2 can be
found in \cite[Section V.4]{Levy}; for Theorem 2.1.3 see \cite{Felgner}.%
\textit{\medskip }

\noindent \textbf{2.1.2.~Theorem.}~\textit{The theories} $\mathrm{ZF(}%
\mathcal{L}(<))+\mathrm{GW}$\textrm{\ }\textit{and}\textrm{\ }$\mathrm{ZF(}%
\mathcal{L}(f))+\mathrm{GC}$\textrm{\ }\textit{are definitionally equivalent
for every language }$\mathcal{L\supseteq L}_{\mathrm{Set}}$\textrm{\textit{.}%
}\footnote{%
Two theories $T_{1}$ and $T_{2}$ are said to be definitionally equivalent if
they have a common definitional extension. Definitional equivalence is also
commonly referred to as \textit{synonymy}, see \cite{Albert-Tehran}.}\medskip

\noindent \textbf{2.1.3.~Theorem.}~\textit{Suppose }$\mathcal{M}\models
\mathrm{ZFC(}\mathcal{L}\mathrm{)}$ \textit{for some countable language} $%
\mathcal{L}$, \textit{and} $\mathrm{Ord}^{\mathcal{M}}$ \textit{has
countable cofinality.} \textit{Then} $\mathcal{M}$ \textit{has an expansion}
$\left( \mathcal{M},<_{M}\right) \models \mathrm{ZF(}\mathcal{L}\mathrm{)+GW}%
.$\textit{\medskip }

The following proposition provides us with a large class of models of set
theory that have definable Skolem functions.\textit{\medskip }

\noindent \textbf{2.1.4.~Proposition.}~\textit{For any language }$\mathcal{L}
$ \textit{that includes} $<$, \textit{every model of} $\mathrm{ZF(}\mathcal{L%
})+\mathrm{GW}$\textrm{\ }\textit{has definable Skolem functions}.\medskip

\noindent \textbf{Proof.}~Given $\varphi =\varphi (x,y_{1},\ldots ,y_{k})$,
we can define a Skolem function $f$ for $\varphi $ by first choosing $\alpha$ to be the first ordinal such that $\exists x\in \mathrm{V}(\alpha )\ \varphi
(x,y_{1},\cdot \cdot \cdot ,y_{k})$, if $\exists x\ \varphi (x,y_{1},\cdot
\cdot \cdot ,y_{k}),$ and then defining $f(y_{1},\cdot \cdot \cdot ,y_{k})$
to be the $<$-first element of:

\begin{center}
$\left\{ x:x\in \mathrm{V}(\alpha )\wedge \varphi (x,y_{1},\cdot \cdot \cdot
,y_{k})\right\}.$
\end{center}

\noindent We define $f(y_{1},\cdot \cdot \cdot ,y_{k})=0$ if $\lnot
\exists x\ \varphi (x,y_{1},\cdot \cdot \cdot ,y_{k}).$ \hfill $\square$

\medskip

\noindent For models of $\mathrm{ZF}$, the $\mathcal{L}_{\mathrm{Set}}$%
-sentence $\exists p\left( \mathrm{V}=\mathrm{HOD}(p)\right) $ expresses:
\textquotedblleft there is some $p$ such that every set is first order
definable in some structure of the form $\left(\mathrm{V}(\alpha ),\in ,p,\beta\right)_{\beta< \alpha}$ with $%
p\in \mathrm{V}(\alpha )$\textquotedblright . The following theorem is
well-known; the equivalence of (a) and (b) will be revisited in Remark
4.3.\medskip

\noindent \textbf{2.1.5.~Theorem.~}\textit{The following statements are
equivalent for} $\mathcal{M}\models \mathrm{ZF}$\textrm{:}\textbf{\medskip }

\noindent \textbf{(a) }$\mathcal{M}\models \exists p\left( \mathrm{V}=%
\mathrm{HOD}(p)\right) .$\textbf{\medskip }

\noindent \textbf{(b) }\textit{For some} $p\in M$ \textit{and some
set-theoretic formula} $\varphi (x,y,z)$, $\mathcal{M}$ \textit{satisfies }%
\textquotedblleft $\varphi (x,y,p)$ \textit{well-orders the universe}%
\textquotedblright .\textbf{\medskip }

\noindent (\textbf{c)} \textit{For some} $p\in M$ \textit{and some
set-theoretic formula} $\psi (x,y,z)$, $\mathcal{M}$ \textit{satisfies }%
\textquotedblleft $\psi (x,y,p)$ \textit{is the graph of a global choice
function}\textquotedblright .\medskip

\begin{center}
\textbf{2.2. Indiscernibles}\medskip
\end{center}

This subsection includes the basic notation and facts about indiscernibles
that will be used in later sections.

\begin{itemize}
\item Given a linear order $(X,<)$, and nonzero $n\in \omega $, we use $%
[X]^{n}$ to denote the set of all \textit{increasing} sequences $x_{1}<\cdot
\cdot \cdot <x_{n}$ from $X$.
\end{itemize}

\noindent \textbf{2.2.1.} \textbf{Definition. }Given a structure $\mathcal{M}
$ and some linear order $(I,<)$ where $I\subseteq M$, we say that $(I,<)$ is%
\textit{\ a set of order indiscernibles in} $\mathcal{M}$ if for any $%
\mathcal{L}(\mathcal{M})$-formula $\varphi (x_{1},\cdot \cdot \cdot ,x_{n})$%
, and any two $n$-tuples $\overline{i}$ and $\overline{j}$ from $[I]^{n}$,
we have:

\begin{center}
$\mathcal{M}\models \varphi (i_{1},\cdot \cdot \cdot ,i_{n})\leftrightarrow
\varphi (j_{1},\cdot \cdot \cdot ,j_{n}).$
\end{center}

The following classical result is due to Ehrenfeucht and Mostowski; see,
e.g.,\ Theorem 3.3.11 of \cite{Chang-Keisler}. In what follows we use the
notation $\mathcal{M}_{I}$ introduced in part (j) of Definition 2.1.1 to
denote the elementary submodel of $\mathcal{M}$ generated by $I$.\medskip

\noindent \textbf{2.2.2.} \textbf{Theorem}.~(Fundamental Theorem of
Indiscernibles) \textit{Suppose} $\mathcal{M}$ \textit{ is a structure with definable Skolem
functions, }$(I,<_{I})$ \textit{is a set of order indiscernibles in} $%
\mathcal{M}$, \textit{and} $\mathcal{L}=\mathcal{L}(\mathcal{M}).$
\medskip

\noindent \textbf{(a)} (Subset Theorem) \textit{For each subset} $I_{0}$
\textit{of} $I$, $\mathcal{M}_{I_{0}}\preceq \mathcal{M}_{I}\preceq \mathcal{%
M}.$ \textit{Moreover, if} $I$ \textit{is infinite and} $I_{0}\neq I$,
\textit{then} $\mathcal{M}_{I_{0}}\prec \mathcal{M}_{I}.$

\medskip

\noindent \textbf{(b)} (Stretching Theorem) \textit{If} $I$ \textit{is
infinite and }$(K,<_{K})$\textit{\ is a linear order, then there is a
model} $\mathcal{M}_{K}\equiv \mathcal{M}$\textit{\ in which} $(K,<_{K})$
\textit{forms a set of indiscernibles}, $K$ \textit{generates} $\mathcal{M}%
_{K},$ \textit{and for any} $\mathcal{L}$-\textit{formula} $\varphi
(x_{1},\cdot \cdot \cdot ,x_{n})$ \textit{we have}:

\begin{center}
$\forall \overline{i}\in \lbrack I]^{n}$ $\forall $ $\overline{k}\in \lbrack
K]^{n}$\ $\mathcal{M}\models \varphi (i_{1},\cdot \cdot \cdot ,i_{n})$ $%
\Longleftrightarrow \mathcal{M}_{K}\models \varphi (k_{1},\cdot \cdot \cdot
,k_{n}).$
\end{center}

\noindent \textbf{(c)} (Elementary Embedding Theorem) \textit{Let} $\mathcal{%
M}_{K}$ \textit{be as in} (b). \textit{Then each injective order-preserving
embedding }$e$ \textit{of} $(I,<_{I})$ \textit{into} $(K,<_{K})$ \textit{%
induces an elementary embedding} $\widehat{e}$ \textit{of} $\mathcal{M}_{I}$
\textit{into} $\mathcal{M}_{K}$, \textit{defined by}

\begin{center}
$\widehat{e}(f(i_{1},\cdot \cdot \cdot ,i_{n}))=f(e(i_{1}),\cdot \cdot \cdot
,e(i_{n})),$
\end{center}

\noindent \textit{where} $f$ \textit{is an} $\mathcal{M}$-\textit{definable
function}. \textit{Moreover, if} $e$ \textit{is surjective, then so is} $%
\widehat{e}$.\medskip

\noindent \textbf{2.2.3.~Remark.}~Since the moreover clause of part (a) of Theorem 2.2.2 is not included in Theorem 3.3.11 of \cite{Chang-Keisler}, we outline its proof here. It suffices to show that if $f(x_1,\cdot \cdot \cdot,x_n)$ is an $\mathcal{M}$-definable function, $(i_1,\cdot \cdot \cdot,i_{n})\in[I]^n$, and $j$ is an element of $I$ such that $j \notin  \{i_1,\cdot \cdot \cdot,i_{n} \}$, then $f(i_1,\cdot \cdot \cdot,i_{n}) \neq j$. By part (b) of Theorem 2.2.2 we may assume that the order-type of $I$ is $\mathbb{Q}$ (the rationals). There are three cases to consider:
\newline
\noindent Case A: $j$ is below $i_1$.
\newline
\noindent Case B: $j$ is between $i_k$ and $i_{k+1}$ , where $1 \leq k \leq n-1$.
\newline
\noindent Case C: $j$ is above $i_n$.
\newline We only consider Case B and leave the other two cases (which are handled similarly) to the reader. Suppose to the contrary that for $1 \leq k \leq n-1$:
\medskip
\newline
\noindent $(*)$~~~$i_{k} < j < i_{k+1}$, and $f(i_1,\cdot \cdot \cdot, i_n)=j$.
\medskip
\newline
\noindent Since the order-type of $I$ is assumed to be $\mathbb{Q}$, there is an element $j^{\prime}\neq j$ in $I$ such that $i_k < j^{\prime} < i_{k+1}$. By indiscernibility, $(*)$ implies that $f(i_1,\cdot \cdot \cdot,i_n)=j^{\prime}$, which is impossible since $f$ is a function.

\medskip

\begin{center}
\textbf{2.3.~Satisfaction classes}\medskip
\end{center}

Satisfaction classes are generalizations of the familiar model-theoretic
notion of `elementary diagram'. They play an important role in this paper;
the material below is the bare minimum that we will need.\medskip

\noindent \textbf{2.3.1.~Definition.}~Reasoning within \textrm{ZF}, for each
object $a$ in the universe of sets, let $c_{a\text{ }}$ be a constant symbol
denoting $a$ (where the map $a\mapsto c_{a}$ is $\Delta _{1}).$ For each
finite extension $\mathcal{L}$ of $\mathcal{L}_{\mathrm{Set}}$, let $\mathrm{%
Sent}_{\mathcal{L}^{+}}(x)$ be the $\mathcal{L}_{\mathrm{Set}}$-formula that
defines the class $\mathrm{Sent}_{\mathcal{L}^{+}}$ of sentences in the
language $\mathcal{L}^{+}=$ $\mathcal{L}\cup \{c_{a}:a\in \mathrm{V}\}$, and
let $\mathrm{Sent}_{\mathcal{L}^{+}}(i,x)$ be the $\mathcal{L}_{\mathrm{Set}%
} $-formula that expresses \textquotedblleft $i\in \omega ,$ $x\in \mathrm{%
Sent}_{\mathcal{L}^{+}},$ and $x$ is a $\Sigma _{i}(\mathcal{L}^{+})$%
-sentence\textquotedblright .\medskip

\noindent \textbf{2.3.2.~Definition.}~Suppose $\mathcal{L}$ is a finite
extension of $\mathcal{L}_{\mathrm{Set}}$, $\mathcal{M}\models \mathrm{ZF}%
\mathcal{(L)}$, $S\subseteq M$, and $k\in \omega ^{\mathcal{M}}$.\medskip

\noindent \textbf{(a)} $S$ is a $\Sigma _{k}$-\textit{satisfaction class for%
} $\mathcal{M}$ if $\left( \mathcal{M},S\right) \models \mathrm{Sat}(k,S)$,
where $\mathrm{Sat}(k,S)$ is the universal generalization of the conjunction
of the axioms $(I)$ through $(IV)$ below. We assume that first order logic
is formulated using only the logical constants $\left\{ \lnot ,\vee ,\exists
\right\} .$\medskip

\noindent $(I)\ \ \left[ \left( S\left( c_{x}=c_{y}\right) \leftrightarrow
x=y\right) \wedge \left( S\left( c_{x}\in c_{y}\right) \leftrightarrow x\in
y\right) \right] .$\medskip

\noindent $(II)\ \ \left[ \mathrm{Sent}_{\mathcal{L}^{+}}(k,\varphi )\wedge
\left( \varphi =\lnot \psi \right) \right] \rightarrow \left[ S(\varphi
)\leftrightarrow \lnot S\mathsf{(}\psi \mathsf{)}\right] \mathsf{.}$\medskip

\noindent $(III)$ $\ \left[ \mathrm{Sent}_{\mathcal{L}^{+}}(k,\varphi
)\wedge \left( \varphi =\psi _{1}\vee \psi _{2}\right) \right] \rightarrow %
\left[ S(\varphi )\leftrightarrow \left( S(\psi _{1})\vee S(\psi
_{2})\right) \right] \mathsf{.}$\medskip

\noindent $(IV)$ $\ \left[ \mathrm{Sent}_{\mathcal{L}^{+}}(k,\varphi )\wedge
\left( \varphi =\exists v\ \psi (v)\right) \right] \rightarrow \left[
S(\varphi )\leftrightarrow \exists x\ S\mathsf{(}\psi \mathsf{(}c_{x}\mathsf{%
))}\right] .$\medskip

\noindent \textbf{(b)} $S$ is a $\Sigma _{\omega }$\textit{-satisfaction
class for} $\mathcal{M}$ if for each $k\in \omega ,$ $S$ is a $\Sigma _{k}$%
-satisfaction class over $\mathcal{M}$. In other words, $S$ is a $\Sigma
_{\omega }$\textit{-}satisfaction class for $\mathcal{M}$ if $S$ agrees
with the usual Tarskian satisfaction class for $\mathcal{M}$ on all standard
$\mathcal{L}$-formulae. Note that if $\mathcal{M}$ is not $\omega $%
-standard, then such a satisfaction class $S$ does not necessarily satisfy
Tarski's compositional clauses for formulae of nonstandard length in $%
\mathcal{M}$. However, using a routine overspill argument, it can be readily
checked that if $S$ is a $\Sigma _{\omega }$\textit{-}satisfaction class
for $\mathcal{M}$ and $S$ is $\mathcal{M}$-amenable, then there is a
nonstandard $c\in \omega ^{\mathcal{M}}$ such that $S$ is a $\Sigma _{c}$%
-satisfaction class over $\mathcal{M}$; indeed, all that is needed for the
overspill argument is for $(\mathcal{M},S)$ to satisfy the scheme of
induction over $\omega ^{\mathcal{M}}$, a scheme that holds in $(\mathcal{M}%
,S)$ since $(\mathcal{M},S)$ satisfies the separation scheme $\mathrm{Sep}(%
\mathcal{L}).$\medskip

\noindent \textbf{(c)} $S$ is a \textit{full satisfaction class for} $%
\mathcal{M}$ if for each $k\in \omega ^{\mathcal{M}},$ $S$ is a $\Sigma _{k}$%
-satisfaction class over $\mathcal{M}$. In other words, $S$ is a full
satisfaction class for $\mathcal{M}$ if $S$ satisfies $(I)$, and the
strengthened versions of $(II)$, $(III)$, and $(IV)$ from part (a) in which
the conjunct $\mathrm{Sent}_{\mathcal{L}^{+}}(k,\varphi )$ is replaced by $%
\mathrm{Sent}_{\mathcal{L}^{+}}(\varphi ).$ Thus, in contrast with $\Sigma
_{\omega }$\textit{-}satisfaction classes which are only guaranteed to
satisfy Tarski's compositional clauses for standard formulae, full
satisfaction classes satisfy Tarski's compositional clauses for all formulae
in $\mathcal{M}$ (including the nonstandard ones, if any).\medskip

\noindent \textbf{(d)} Recall that given any language $\mathbb{L}$, $\mathbb{%
L}_{\infty ,\infty }$ is the union of logics $\mathbb{L}_{\kappa ,\lambda }$%
, where $\kappa $ and $\lambda $ are infinite cardinals and the logic $%
\mathbb{L}_{\kappa ,\lambda }$ is the extension of first order logic that
allows conjunctions and disjunctions of sets of formulae of cardinality less
than $\kappa $ and blocks of existential quantifiers and blocks of universal quantifiers of length less than $
\lambda$. Thus $\mathbb{L}_{\omega ,\omega }$ is none other than the usual
first order logic based on the language $\mathbb{L}$. $S$ is an $\mathcal{L}%
_{\infty ,\infty }$\textit{-satisfaction class for} $\mathcal{M}$ if $S$
satisfies $(I)$, the strengthened version of $(II)$ from part (a) in which
the conjunct $\mathrm{Sent}_{\mathcal{L}^{+}}(k,\varphi )$ is replaced by
the formula $\mathrm{Sent}_{\mathcal{L}_{\infty ,\infty }^{+}}(\varphi )$
that expresses \textquotedblleft $\varphi $ is a sentence of $\mathcal{L}%
_{\infty ,\infty }^{+}$\textquotedblright , as well as the following
stronger variants of $(III)$ and $(IV)$. Note that in $(III)^{\ast }$ below $%
\Psi $ ranges over \textit{sets} of formulae of $\mathcal{L}_{\infty ,\infty
}$ \medskip

\noindent $(III)^{\ast }$ $\ \left[ \mathrm{Sent}_{\mathcal{L}_{\infty
,\infty }^{+}}(\varphi )\wedge \left( \varphi =\bigvee \Psi \right) \right]
\rightarrow \left[ S(\varphi )\leftrightarrow \exists \psi \in \Psi \ S%
\mathsf{(}\psi \mathsf{)}\right] .$\medskip

\noindent $(IV)^{\ast }$ $\ \left[ \mathrm{Sent}_{\mathcal{L}_{\infty
,\infty }^{+}}(\varphi )\wedge \left( \varphi =\exists \left\langle
x_{\alpha }:\alpha <\lambda \right\rangle \ \psi \left( x_{\alpha }:\alpha
<\lambda \right) \right) \right] \rightarrow $

\begin{center}
$\left[ S(\varphi )\leftrightarrow \left( \exists \left\langle x_{\alpha
}:\alpha <\lambda \right\rangle \ S\mathsf{(}\psi \left( c_{x_{\alpha
}}:\alpha <\lambda \right) \right) \right] .$
\end{center}

\begin{itemize}
\item Given a satisfaction class $S$, in the interest of a lighter notation,
we will often write $\varphi \left( a_{1},\cdot \cdot \cdot ,a_{n}\right)
\in S$ instead of $\varphi \left( c_{a_{1}},\cdot \cdot \cdot
,c_{a_{n}}\right) \in S.$
\end{itemize}

\noindent \textbf{2.3.3.~Remark.}~It is a well-known result of Levy that if $%
\mathcal{M}\models \mathrm{ZF,}$ then there is a $\Sigma _{0}$-satisfaction
class for\textit{\ }$\mathcal{M}$\textit{\ }that is\textit{\ }definable in $%
\mathcal{M}$ by a $\Sigma _{1}$-formula (see \cite[p.~186]{Jechbook-2003}
for a proof). This makes it clear that for each $n\geq 1,$ there is a $%
\Sigma _{n}$-satisfaction class for\textit{\ }$\mathcal{M}$\textit{\ }that
is\textit{\ }definable in $\mathcal{M}$ by a $\Sigma _{n}$-formula. Levy's
result extends to models of $\mathrm{ZF(}\mathcal{L}\mathrm{)}$ if $\mathcal{%
L}$\ is finite. We use $\mathrm{Sat}_{\Delta _{0}}$ to refer to the
canonical $\Sigma _{0}$-satisfaction class (recall that by definition $%
\Delta _{0}=\Sigma _{0}$ in the Levy Hierarchy). \medskip

\begin{center}
\textbf{2.4.~The theory }$\mathbf{GBC\ }+$\ $``$\textbf{Ord is weakly
compact\textquotedblright }\medskip
\end{center}

The theory \textrm{GBC} $+$ \textrm{\textquotedblleft Ord is weakly
compact\textquotedblright }\textbf{\ }was first studied by McAloon and
Ressayre \cite{McAloon-Ressayre (French)}, and then later, using different
methods and motivations, by the author \cite{Ali NFUA}. Here we bring
together a number of results about this theory that are not only of
intrinsic foundational interest, but also play an essential role in the
proofs of the results in later sections.\medskip

\noindent \textbf{2.4.1.~Definition.}~$\mathrm{GBC}$\ is the G\"{o}%
del-Bernays theory of classes $\mathrm{GB}$ with global choice.\footnote{$%
\mathrm{GB}$\ is also referred to as $\mathrm{BG}$, $\mathrm{VNB}$ (\textit{%
von Neumann-Bernays) }and $\mathrm{NBG}$ (\textit{von Neumann-Bernays-G\"{o}%
del})\textit{\ }in the literature\textit{. }In some sources $\mathrm{GB}$
includes the global axiom of choice. It is well-known that GB is finitely
axiomatizable \cite[Exercise 13.5]{Jechbook-2003}.} Our set-up for $\mathrm{%
GB}$ is the standard one in which models of $\mathrm{GB}$ are viewed as
\textit{two-sorted} structures of the form $(\mathcal{M},\mathfrak{X})$,
where $\mathcal{M}\models \mathrm{ZF}$, and $\mathfrak{X}\subseteq \mathcal{P%
}\mathfrak{(}M)$. Thus, the language appropriate to $\mathrm{GB}$ (referred
to as the \textit{language of class theory}) is a two-sorted language: a
sort for sets (represented by lower case letters), a sort for classes
(represented by upper case letters), and a special membership relation
symbol $\in $ for indicating that a set $x$ is a member of a class $X$,
written $x\in X$. In the interest of a lighter notation, we use $\in $ both
as the formal symbol indicating membership between sets, and also for the
membership relation between sets and classes (since we use upper case
letters to symbolize classes, there is no risk of confusion). Also, since
coding of sequences is available in $\mathrm{GB}$, we shall use expressions
such as \textquotedblleft $F\in \mathfrak{X}$\textquotedblright , where $F$
is a function, as a substitute for the precise but lengthier expression
\textquotedblleft there is a class in $\mathfrak{X}$ that canonically codes $%
F$\textquotedblright . We will say $X\in \mathfrak{X}$ is a \textit{proper
class} if there is no $c\in M$ such that $\mathrm{Ext}_{\mathcal{M}}(c)=X$, otherwise we say that $X$ is coded as a set in $\mathcal{M}$.\medskip

\noindent \textbf{2.4.2.~Remark.}~It is well-known that for $\mathfrak{X}%
\subseteq \mathcal{P}(M),$ and $\mathcal{M}\models \mathrm{ZF}$, $(\mathcal{M%
},\mathfrak{X})\models \mathrm{GB}$ iff the following two conditions
hold:\medskip

\noindent \textbf{(a)} If $X_{1},\cdot \cdot \cdot ,X_{n}\in \mathfrak{X}$,
then $(\mathcal{M},X_{1},\cdot \cdot \cdot ,X_{n})\models \mathrm{ZF}%
(X_{1},\cdot \cdot \cdot ,X_{n})$.\medskip

\noindent \textbf{(b)} If $X_{1},\cdot \cdot \cdot ,X_{n}\in \mathfrak{X},$
and $Y$ is parametrically definable in $(\mathcal{M},X_{1},\cdot \cdot \cdot
,X_{n})$, then $Y\in \mathfrak{X}.$\medskip

\noindent \textbf{2.4.3.~Definition.}~\textquotedblleft $\mathrm{Ord}$ is
weakly compact\textquotedblright\ is the statement in the language of class
theory asserting that every \textrm{Ord}-tree has a branch, where \textrm{Ord%
}-trees are defined in analogy with the familiar notion of $\kappa $-trees
in infinite combinatorics: $(\tau ,<_{\tau })$ is an \textrm{Ord}-tree, if $%
(\tau ,<_{\tau })$ is a well-founded tree of height \textrm{Ord} such that
the collection of nodes of any prescribed ordinal rank is a set (as opposed
to a proper class). \medskip

The following result is the $\mathrm{GBC}$-adaptation of the $\mathrm{ZFC}$%
-formulation of the classical Erd\H{o}s-Hajnal-Rado Ramification Lemma. The
Ramification Lemma is a $\mathrm{ZFC}$-theorem with a parameter $\kappa $
that ranges over infinite cardinals $\kappa $; in the $\mathrm{GBC}$%
-adaptation below the class of ordinals $\mathrm{Ord}$ plays the role
typically played by $\kappa $. Lemma 2.4.4 shows that within each model $(%
\mathcal{M},\mathfrak{X})$ of $\mathrm{GBC}$, one can canonically associate
an $\mathrm{Ord}$-tree $\tau _{F}$ to each coloring of $F$ of $[\mathrm{Ord}%
]^{n+1}$, where $1\leq n\in \omega ^{\mathcal{M}}$ into set-many colors such
that the color associated by $F$ to each increasing chain of length $n+1$ in
$\tau _{F}$ is independent of the maximum element of the chain.\medskip

\noindent \textbf{2.4.4.~Lemma.}~\textit{Suppose} $(\mathcal{M},\mathfrak{X}%
)\models \mathrm{GBC}$, $1\leq n\in \omega ^{\mathcal{M}}$, \textit{and} $F:[%
\mathrm{Ord}]^{n+1}\rightarrow \lambda ,$ \textit{where} $F\in \mathfrak{X}$
\textit{and} $\lambda $ \textit{is a cardinal in} $\mathcal{M}$. \textit{%
There is a structure} $\tau _{F}=$ $\left( \mathrm{Ord}^{\mathcal{M}%
},<_{F}\right) $ \textit{coded in }$\mathfrak{X}$ \textit{such that the
following hold} \textit{in }$(\mathcal{M},F,\tau _{F})$:\medskip

\noindent \textbf{(a)} \textit{For all ordinals} $\alpha $ \textit{and} $%
\beta ,$ \textit{if} $\alpha <_{F}\beta $, \textit{then} $\alpha \in \beta .$
\textit{In particular}, $\tau _{F}$ \textit{is a well-founded tree.}\medskip

\noindent \textbf{(b)} $F(\alpha _{1},\alpha _{2},\cdot \cdot \cdot ,\alpha
_{n},\alpha _{n+1})=F(\alpha _{1},\alpha _{2},\cdot \cdot \cdot ,\alpha
_{n},\beta )$ \textit{whenever }

\begin{center}
$\alpha _{1}<_{F}\alpha _{2}<_{F}\cdot \cdot \cdot <_{F}\alpha
_{n-1}<_{F}\alpha _{n}<_{F}\alpha _{n+1}$, and $\alpha _{n}<_{F}\beta
.\medskip $
\end{center}

\noindent \textbf{(c) }\textit{For each ordinal }$\alpha $, \textit{the} $%
\alpha $-\textit{th level of the tree }$\tau _{F}$
\textit{has cardinality at most} $\lambda ^{\left\vert \omega +\alpha
\right\vert }$; \textit{in particular} $\tau _{F}$ \textit{is an }$\mathrm{%
Ord}$\textit{-tree}.$\medskip $

\noindent \textbf{Proof.}~The \textrm{ZFC}-proof presented in \cite[Lemma 7.2%
]{Kanamori book} can be readily adapted to the $\mathrm{GBC}$ context by
replacing the cardinal $\sigma $ in that proof with the proper class $%
\mathrm{Ord}$.\hfill $\square \medskip $

The above Lemma lies at the heart of the proof of the theorem below. In
parts (b) and (c) of the theorem, $\mathrm{Ord}\rightarrow \left( \mathrm{Ord%
}\right) _{\kappa }^{n}$ stands for the sentence in the language of class
theory that asserts that for every class function $F:[\mathrm{Ord}%
]^{n}\rightarrow \kappa $ (where $0<n\in \omega $, $\kappa $ is a finite or
infinite cardinal, and $[\mathrm{Ord}]^{n}$ is the class of all increasing
sequences of ordinals of length $n$) there is an unbounded $H\subseteq
\mathrm{Ord}$ such that $H$ is $F$-homogeneous, i.e., for any two increasing
$n$-tuples $\overline{x}$ and $\overline{y}$ from $H$, $F(\overline{x})=F(%
\overline{y})$.\medskip

\noindent \textbf{2.4.5.~Theorem.}~\textit{Suppose} $(\mathcal{M},\mathfrak{X%
})\models \mathrm{GBC}$ $+$ \textrm{\textquotedblleft Ord is weakly
compact\textquotedblright }. \textit{Then}:$\medskip $

\noindent \textbf{(a) }\textit{If} $1\leq n\in \omega ^{\mathcal{M}}$, $F\in
\mathfrak{X}$\textit{, }$\kappa \in \mathrm{Ord}^{\mathcal{M}}$ \textit{and}
$(\mathcal{M},F)\models F:[\mathrm{Ord}]^{n+1}\rightarrow \kappa $, \textit{%
then there is some proper class }$H\in \mathfrak{X}$ \textit{that is
`end-homogeneous', i.e.,} $(\mathcal{M},F,H)$ \textit{satisfies}:

\begin{center}
$\forall \overline{\alpha }\in \lbrack H]^{n+2}\ F(\alpha _{1},\cdot \cdot
\cdot ,\alpha _{n},\alpha _{n+1})=F(\alpha _{1},\cdot \cdot \cdot ,\alpha
_{n},\alpha _{n+2}).$
\end{center}

\noindent \textbf{(b) }\textit{For every cardinal }$\kappa $ \textit{in }$%
\mathcal{M}$, $(\mathcal{M},\mathfrak{X})\models \forall n\in \omega
\backslash \{0\}$ $\left( \varphi (n,\kappa )\rightarrow \varphi (n+1,\kappa
)\right) ,$ \textit{where}:

\begin{center}
$\varphi (n,\kappa ):=\left( \mathrm{Ord}\rightarrow \left( \mathrm{Ord}%
\right) _{\kappa }^{n}\right) .\medskip $
\end{center}

\noindent \textbf{(c)} \textit{If} $1\leq n\in \omega $, \textit{and} $%
\kappa $ \textit{is a cardinal of }$\mathcal{M}$, \textit{then} $(\mathcal{M}%
,\mathfrak{X})\models \mathrm{Ord}\rightarrow \left( \mathrm{Ord}\right)
_{\kappa }^{n}.$ $\footnote{%
As shown in Theorem 4.9 the statement\textit{\ }$\theta =\forall m,n\in
\omega \left( \mathrm{Ord}\rightarrow \left( \mathrm{Ord}\right)
_{m}^{n}\right) $ is not provable in $\mathrm{GBC}$ $+$ \textit{%
\textquotedblleft }\textrm{Ord is weakly compact}\textit{\textquotedblright
, }but part (b) of Theorem 2.4.5 shows that $\theta $ is provable in the
theory obtained by augmenting\textit{\ } GBC $+$ \textit{\textquotedblleft }%
\textrm{Ord is weakly compact}\textit{\textquotedblright } with $\Pi
_{2}^{1} $-induction (over the ambient $\omega ).$}\medskip $

\noindent \textbf{Proof.}~To verify (a), we argue in $(\mathcal{M},\mathfrak{%
X})$. Suppose $F:[\mathrm{Ord}]^{n+1}\rightarrow \kappa ,$ where $1\leq n\in
\omega ^{\mathcal{M}}$ and $F\in \mathfrak{X}$, and let $\tau _{F}$ be as in
Lemma 2.4.4. By weak compactness of Ord, there is some proper class $%
H\subseteq \mathrm{Ord}$ that is a cofinal branch of $\tau _{F}.$ Lemma
2.4.4 assures us that $F(\alpha _{1},\cdot \cdot \cdot ,\alpha _{n},\beta
)=F(\alpha _{1},\cdot \cdot \cdot ,\alpha _{n},\beta ^{\prime })$ if $%
\overline{\alpha }\in \lbrack H]^{n}$ and $\beta $ and $\beta ^{\prime }$
are any two elements of $H$ that are above $\alpha _{n}.$ Thus $H$ is
end-homogeneous, as desired.$\medskip $

To see that (b) holds, suppose $(\mathcal{M},\mathfrak{X})\models \varphi
(n,\kappa )$ for some nonzero $n\in \omega ^{\mathcal{M}}$ and some cardinal
$\kappa $ of $\mathcal{M}$. To verify that $(\mathcal{M},\mathfrak{X}%
)\models \varphi (n+1,\kappa ),$ suppose that for some $F\in \mathfrak{X},(%
\mathcal{M},F)\models F:[\mathrm{Ord}]^{n+1}\rightarrow \kappa $. By (a) we
can get hold of an end-homogeneous $H$ for $F$. Consider the function $%
G:[H]^{n}\rightarrow \kappa $ defined in $(\mathcal{M},F)$ by:

\begin{center}
$G(\alpha _{1},\cdot \cdot \cdot ,\alpha _{n}):=F(\alpha _{1},\cdot \cdot
\cdot ,\alpha _{n},\beta ),$ where $\beta \in H$ and $\beta >\alpha _{n}$.
\end{center}

\noindent The end-homogeneity of $H$\ assures us that $G$ is well-defined.
Hence by the assumption that $\varphi (n,\kappa )$ holds in $(\mathcal{M},%
\mathfrak{X})$, there is a proper class $H^{\prime }\subseteq H$ that is $G$%
-homogeneous. This makes it evident that $H^{\prime }$ is $F$-homogeneous,
thus completing the proof of (b).$\medskip $

(c) follows immediately from (b) by induction on metatheoretic natural
numbers $n.$\hfill $\square \medskip $

\begin{itemize}
\item Next we will describe a minor extension of another tree construction,
first introduced in \cite[Section 3]{Ali powerlike}, and later simplified in
\cite[Definition 2.2]{Ali-Joel}, where it was used to prove that models of $%
\mathrm{GBC}$ of the form $(\mathcal{M},\mathfrak{X})$, where $\mathfrak{X}$
is the collection of parametrically $\mathcal{M}$-definable subsets of $%
\mathcal{M}\models \mathrm{ZFC}$, never satisfy the axiom \textquotedblleft
\textrm{Ord} is weakly compact\textquotedblright .\medskip
\end{itemize}

\noindent \textbf{2.4.6.~Definition.}~Suppose $(\mathcal{M},\mathfrak{X}%
)\models \mathrm{GBC}$. Fix some ordering $<_{M}$ of $\mathcal{M}$ in $%
\mathfrak{X}$ such that $\left( \mathcal{M},<_{M}\right) \models \mathrm{GW.}
$ Within $\left( \mathcal{M},<_{M}\right) $, given ordinals $\alpha \in
\beta $ let:

\begin{center}
$\mathcal{V}_{\beta ,\alpha }=\left( \mathrm{V}(\beta ),\in ,<,a\right)
_{a\in \mathrm{V}(\alpha )}.$
\end{center}

\noindent Thus $\mathcal{V}_{\beta ,\alpha }$ is an $\mathcal{L}_{\alpha }$%
-structure, where $\mathcal{L}_{\alpha }$ is the result of augmenting $%
\mathcal{L}_{\mathrm{Set}}(<)$ with constant symbols $c_{a}$ for each $a\in
\mathrm{V}(\alpha ).$ Given $X\in \mathfrak{X}$\textit{\ }and\textit{\ }$%
n\in \omega ,$ within $\left( \mathcal{M},<_{M},X\right) ,$ let $\tau _{n}(X)$
be the tree whose elements are of the form:

\begin{center}
$T(X,\beta ,\alpha ,s):=\mathrm{Th}\left( \mathcal{V}_{\beta ,\alpha },X\cap
\mathrm{V}(\beta ),s\right) ,$
\end{center}

\noindent where $s\in \mathrm{V}(\beta )\backslash \mathrm{V}(\alpha ),$
with the additional requirement that:

\begin{center}
$\left( \mathrm{V}(\beta ),\in ,<,X\cap \mathrm{V}(\beta )\right) \prec
_{\Sigma _{n}(X)}\left( \mathrm{V},\in ,<,X\right) .$\medskip
\end{center}

\noindent Note that $T(X,\beta ,\alpha ,s)$ consists of $\mathcal{L}_{\alpha
}(X,c)$-sentences that hold in $\left( \mathcal{V}_{\beta ,\alpha },X\cap \mathrm{%
V}(\beta ),s\right) $, where:

\begin{center}
$\mathcal{L}_{\alpha }(X,c):=\{<, X, c\} \cup \{c_{a}:a\in \mathrm{V}(\alpha )\}.$
\end{center}

\noindent In the above, $X$ is a unary predicate (that is conflated with its
denotation), and $c$ is a new constant symbol whose denotation is $s$. The
ordering relation on $\tau _{n}(X)$ is set-inclusion.\medskip

\noindent \textbf{2.4.7.~Theorem.~}\textit{Suppose} $(\mathcal{M},\mathfrak{X%
})\models \mathrm{GBC}$, \textit{and let} $<_{M}$ \textit{be a member of} $%
\mathfrak{X}$ \textit{such that} $\left( \mathcal{M},<_{M}\right) \models
\mathrm{GW.}$ \textit{Then}:\medskip

\noindent \textbf{(a)}\textit{\ For each }$X\in \mathfrak{X},$ $(\mathcal{M},%
\mathfrak{X})\models \ ``\tau _{n}(X)$ \textit{is an} \textrm{Ord}-\textit{%
tree}\textquotedblright.\medskip

\noindent \textbf{(b) }\textit{If }$n\geq 1$ \textit{and the tree }$\tau
_{n}(X)$ \textit{as computed in} $\mathcal{M}^{+}:=\left( \mathcal{M},<_{M},X\right) $
\textit{has a branch }$B\in \mathfrak{X}$, \textit{then} \textit{there is an
}$\mathcal{L}_{\mathrm{Set}}(<,X)$-\textit{structure }$\mathcal{N}^{+}:=\left( \mathcal{N}%
,<_{N},X_{N}\right)$\textit{\ and
a proper }$\Sigma _{n}(<,X)$-\textit{elementary end} \textit{%
embedding }

\begin{center}
$j:\mathcal{M}^{+}\rightarrow \mathcal{N}^{+}.$
\end{center}

\noindent \textit{Moreover, }$\mathfrak{X}$ \textit{contains} \textit{both
the embedding} $j$, \textit{and a full satisfaction class for the structure} $\mathcal{N}^{+}$.

\medskip

\noindent \textbf{(c)} \textit{If} $(\mathcal{M},\mathfrak{X})\models ``%
\mathrm{Ord}$ $\mathrm{is}$ $\mathrm{weakly\ compact}$\textquotedblright $,$
\textit{then} \textit{for every} $X\in \mathfrak{X}$\textit{\ and every }$%
n\in \omega $,\textit{\ there is a }$\Sigma _{n}(<,X)$\textit{-e.e.e. }$
\mathcal{N}^{+}$ \textit{of }$\mathcal{M}^{+}:=\left( \mathcal{M},<_{M},X\right) $ \textit{such that}
\textrm{Ord}$^{\mathcal{N}}\backslash M$ \textit{has a minimum element.
Consequently, there is some }$S_{X}\in \mathfrak{X}$ \textit{that is a full
satisfaction class for} $\mathcal{M}^{+}$; \textit{indeed there is even some
}$S_{X,\infty }\in \mathfrak{X}$\textit{\ such that }$S_{X,\infty }$\textit{%
\ is an }$\mathcal{L}_{\infty ,\infty }$-\textit{satisfaction class} \textit{%
for} $\mathcal{M}^{+}.$\medskip

\noindent \textbf{Proof.}~The proofs of\textbf{\ }(a) and (b) are minor
variants of Lemmas 2.3 and 2.5 of \cite{Ali-Joel}, so we do not present them
here. To prove (c), given $X\in \mathfrak{X}$\textit{\ }and\textit{\ }$n\in
\omega $,\textit{\ }we first use (b) and the assumption that $(\mathcal{M},%
\mathfrak{X})\models \ ``\mathrm{Ord}$ $\mathrm{is}$ $\mathrm{weakly\ compact%
}$\textquotedblright\ to construct a $\Sigma _{n}(<,X)$-e.e.e.$~\mathcal{N}^{+}$ of $\mathcal{M}^{+}$. Then we use the following result to arrange for \textrm{Ord}$^{%
\mathcal{N}}\backslash M$ to have a minimum element. Note that this
immediately implies the existence of the satisfaction classes $S_{X}$ and $%
S_{X,\infty }$ as in the second assertion in (c) since if $n\geq 2$, then $\mathcal{N}^{+}$ is a model of a
substantial fragment of $\mathrm{ZF}$, including $\mathrm{KP}$
(Kripke-Platek set theory), and already $\mathrm{KP}$\ is sufficient for
defining the $\mathcal{L}_{\infty ,\infty }$-satisfaction predicate for
every $\mathcal{L}$-set-structure. \cite[%
III.2]{Barwise-book}.
\medskip

\noindent \textbf{2.4.8.~Theorem.}~\textit{Suppose} $(\mathcal{M},\mathfrak{X%
})\models \mathrm{GBC,}$ $X\in \mathfrak{X}$, \textit{and} $\left(\mathcal{M},<_{M}\right) \models
\mathrm{GW}$, \textit{where} $<_{M}$ \textit{is in} $\mathfrak{X}$. \textit{Suppose furthermore that }$\mathfrak{X}$
\textit{contains a full satisfaction class for }$\mathcal{M}^{+}:=\left( \mathcal{M}%
,<_{M},X_{M}\right) $ \textit{and also a full satisfaction class for some} $\Sigma
_{n+3}(<,X)$-\textit{e.e.e.~}$\mathcal{N}^{+}:=\left( \mathcal{N},<_{N},X_{N}\right) $ \textit{of} $\mathcal{M}^{+}$, \textit{where }$n\geq 1.$ \textit{Then there is some} $\mathcal{K}^{+}:=\left( \mathcal{K},<_{K},X_{K}\right) $ \textit{such that}:\medskip

\noindent \textbf{(a)} $\mathcal{K}^{+}$\textit{\ is a} $%
\Sigma _{n+1}(<,X)$-\textit{e.e.e.}$~$\textit{of} $\mathcal{M}^{+}.$%
\medskip

\noindent \textbf{(b)} $\mathfrak{X}$ \textit{contains a full} \textit{%
satisfaction class for} $\mathcal{K}^{+}$, \textit{and}
\medskip

\noindent \textbf{(c)} \textrm{Ord}$^{\mathcal{K}}\backslash M$ \textit{has
a minimum element.}$\medskip $

\noindent \textbf{Proof.}~The proof is similar to the proofs of \cite[%
Theorem 3.3]{Ali-Cons} and \cite[Theorem 2.1]{Ali-Joel}. Choose $S\in $ $%
\mathfrak{X}$ such that $S$ is a full satisfaction class for $\mathcal{M}^{+}$, where $\mathcal{M}^{+} \prec
_{\Sigma _{n+3}(<,X)}\mathcal{N}^{+}.$ For $n\in \omega $
consider the statement $\varphi _{n}$ that expresses the following instance
of the reflection theorem:

\begin{center}
$\forall \lambda \in \mathrm{Ord}\mathbf{\ }\exists \beta \in \mathrm{Ord}%
\left( \lambda \in \beta \wedge \ (\mathrm{V}(\beta ),<,\in ,X\cap \mathrm{V}%
(\beta ))\prec _{\Sigma _{n+1}(X)}\left( \mathrm{V},<,\in ,X\right) \right) .$
\end{center}

\noindent Note that $\varphi _{n}$ is a $\Pi _{n+3}(<,X)$-statement for each $1\leq k\in
\omega $ since the satisfaction predicate for $%
\Sigma _{k}(<,X)$-formulae is $\Sigma _{k}(<,X)$-definable. Therefore $\mathcal{N}^{+} \models \varphi $ since $\varphi $ holds in $%
\mathcal{M}^{+}$ by the reflection theorem. So we can fix some $\lambda \in
\mathrm{Ord}^{\mathcal{N}}\backslash \mathrm{Ord}^{\mathcal{M}}$ and some $%
\mathcal{N}$-ordinal $\beta >\lambda $ such that:

\begin{center}
$\mathcal{N}^{+} (\beta) \prec _{\Sigma
_{n+1}(<,X)}\mathcal{N}^{+} $,
\end{center}

\noindent where $\mathcal{N}^{+}(\beta ):=(\mathrm{V}(\beta ),<,\in )^{\mathcal{N}^{+}}. $ Note that this implies that $\mathcal{N}^{+}(\beta )$ can meaningfully
define the satisfaction predicate for every set-structure `living in'\ $%
\mathcal{N}^{+}(\beta )$. For any $\alpha \in \mathrm{Ord}^{\mathcal{M}}$ with $\alpha <\beta $,
\textit{within} $\mathcal{N}^{+} $ one can define the
submodel $\mathcal{K}^{+}_{\alpha}:= \left( \mathcal{K}_{\alpha },<_{K_{\alpha}},X_{K_{\alpha}} \right) $ of $\left(
\mathrm{V}(\beta),\in, <, X\cap \mathrm{V}(\beta )\right) $ whose universe $%
K_{\alpha }$ is defined via:

\begin{center}
$K_{\alpha }:=\{a\in \mathrm{V}(\beta ):a$ is first order definable in $(%
\mathrm{V}(\beta ),\in ,<,X\cap \mathrm{V}(\beta ),\lambda ,m)_{m\in
\mathrm{V}(\alpha )}\}$.
\end{center}

\noindent Clearly $M(\alpha)\cup \{\lambda \}\subsetneq K_{\alpha }$ and $\mathcal{K}^{+}_{\alpha} \prec_{\Sigma_{n+1}(<,X)} \mathcal{N}^{+} $, and of course $\mathcal{K}^{+}_{\alpha }$ is coded
in $\mathcal{N}^{+}$. Next let:

\begin{center}
$K:=\bigcup\limits_{\alpha \in \mathrm{Ord}^{\mathcal{M}}}K_{\alpha },$ and $%
X_{K}:=X_{N}\cap K,$
\end{center}

\noindent and let the $\mathcal{L}_{\mathrm{Set}}(<,X)$-structure $\mathcal{K}^{+}$ be the submodel of $\mathcal{N}^{+} (\beta) $ whose universe is $K$. Note that:

\begin{center}
$\mathcal{M}^{+} \prec _{\text{\textrm{end},}\Sigma
_{n+1}(<,X)}\mathcal{K}^{+},$
\end{center}

\noindent since:

\begin{center}
$\mathcal{M}^{+}\prec _{\mathrm{end}}\mathcal{K}^{+}\preceq \mathcal{N}^{+}(\beta) \prec
_{\Sigma _{n+1}(<,X)}\mathcal{N}^{+} .$
\end{center}

\begin{itemize}
\item Observe that if $S$ is a full satisfaction class for $\mathcal{N}^{+} $ such that $S\in \mathfrak{X}$, then there are full
satisfaction classes in $\mathfrak{X}$ for the structures $\mathcal{N}^{+}(\beta )$ and $\mathcal{K}^{+}$.
\end{itemize}

To prove that\textbf{\ }$\mathrm{Ord}^{\mathcal{K}}\backslash \mathrm{Ord}^{%
\mathcal{M}}$ has a least element, suppose to the contrary that $\mathrm{Ord}%
^{\mathcal{K}}\backslash \mathrm{Ord}^{\mathcal{M}}$ has no least element. Within $\mathcal{N}^{+} $ let $S_{\alpha}$ be the full satisfaction class for the
structure $$\left( \mathrm{V}(\beta ),\in ,<,X\cap \mathrm{V}(\beta ),\lambda
,m \right)_{m\in \mathrm{V}(\alpha )}.$$ Now let $\Phi :=\bigcup\limits_{\alpha \in \mathrm{Ord}^{\mathcal{M}}}\Phi
_{\alpha }$, where:

\begin{center}
$\Phi _{\alpha }:=\{\varphi (c,c_{m})\in M:\mathcal{N}^{+}
\models $ $\varphi (c_{\lambda },c_{m})\in S_{\alpha}\}\mathrm{,}$
\end{center}

\noindent  Note that $\Phi \in \mathfrak{X}$ since $%
S_{\alpha}\in \mathfrak{X}$, in particular $\Phi $ is $\mathcal{M}^{+}$-amenable. Also observe that in the above definition of $\Phi _{\alpha }$, $\varphi
(c,c_{m})$ ranges over formulae of the language $\mathcal{L}_{\alpha }(X,c)$ (in
the sense of $\mathcal{M)}$, where $\mathcal{L}_{\alpha }(X,c)$ is as in Definition 2.4.6.
 Also note that the constant $c$ is interpreted as $%
\lambda $ in the right-hand-side of the above definition of $\Phi _{\alpha
}. $ Thus $\Phi $ can be thought of as the \textit{type} of the element $%
\lambda $ in the structure $\mathcal{N}^{+}(\beta )$ over the parameter set $M$ (with the important provision that $%
\Phi $ includes nonstandard formulae if $\mathcal{M}$ is $\omega $%
-nonstandard). Now let:

\begin{center}
$\Gamma :=\left\{ t(c,c_{m})\in M:t(c,c_{m})\in \Phi \ \mathrm{and}\ \forall
\theta \in \mathrm{Ord\ }\left( t(c,c_{m})>c_{\theta }\right) \in \Phi
\right\} ,$
\end{center}

\noindent where $t$ is a definable function in the language $\mathcal{L}%
^{+}=\bigcup\limits_{\alpha \in \mathrm{Ord}^{\mathcal{M}}}\mathcal{L}%
_{\alpha }(X,c)$. So, officially speaking, $\Gamma $ consists of \textit{%
syntactic objects} $\varphi (c,c_{m},x)$ in $\mathcal{M}$ that satisfy the
following three conditions in $(\mathcal{M},\Phi )$:\medskip

\noindent (1) $\left[ \exists !x\varphi (c,c_{m},x)\right] \in \Phi .$%
\medskip

\noindent (2) $\left[ \forall x\left( \varphi (c,c_{m},x)\rightarrow x\in
\mathrm{Ord}\right) \right] \in \Phi .$\medskip

\noindent (3) $\forall \theta \in \mathrm{Ord}\ \left[ \forall x\left(
\varphi (c,c_{m},x)\rightarrow c_{\theta }\in x\right) \right] \in \Phi .$%
\medskip

\noindent Note that $\Gamma $ is definable in $\left( \mathcal{M}^{+},\Phi
\right) .$ Since we assumed that $\mathrm{Ord}^{\mathcal{K}}\backslash
\mathrm{Ord}^{\mathcal{M}}$ has no minimum element, $\left( \mathcal{M}^{+},\Phi
\right) \models \psi $, where:

\begin{center}
$\psi :=\forall t\left( t\in \Gamma \rightarrow \left( \exists t^{\prime
}\in \Gamma \wedge \left[ t^{\prime }\in t\right] \in \Phi \right) \right) .$
\end{center}

\noindent Choose $k\in \omega $ such that $\psi $ is a $\Sigma _{k}(<,X,\Phi
\mathrm{)}$-statement$,$ and use the reflection theorem in $\left( \mathcal{M}^{+},\Phi \right) $ to pick $\mu \in \mathrm{Ord}^{\mathcal{M}}$ such
that:

\begin{center}
$\left( \mathcal{M}^{+}(\mu ),\Phi \cap M\mathcal{(}\mu )\right)
\prec _{\Sigma _{k}(<,X,\Phi )}\left( \mathcal{M}^{+},\Phi \right) .$
\end{center}

\noindent Then $\psi $ holds in $\left( \mathcal{M}^{+}(\mu ),\Phi
\cap M(\mu )\right) $, so by DC (dependent choice, which holds in $%
\mathcal{M}$ since AC holds in $\mathcal{M}$), there is some function $f_{c}$
in $\mathcal{M}$ such that:

\begin{center}
$\left( \mathcal{M}^{+},\Phi \right) \models \forall n\in \omega \ \left[
f_{c}(n+1)\in f_{c}(n)\right] \in \Phi .$
\end{center}

\noindent Let $\alpha \in \mathrm{Ord}^{\mathcal{M}}$ be large enough so
that $M(\alpha)$ contains all constants $c_{m}$ that occur in any of the
terms in the range of $f$; let $f_{\lambda }(n)$ be defined in $\mathcal{N}$
as the result of replacing all occurrences of the constant $c$ with $%
c_{\lambda }$ in $f_{c}(n)$; and let $g(n)$ be defined in $\mathcal{N}^{+} $ as the interpretation of $f_{\lambda }(n)$ in:

\begin{center}
$\left( \mathrm{V}(\beta ),\in ,<,X\cap \mathrm{V}(\beta ),\lambda ,m\right)
_{m\in \mathrm{V}_{\alpha }}.$
\end{center}

\noindent Then $\mathcal{N}^{+} $ satisfies:

\begin{center}
$\forall n\in \omega \ \left( g(n)\in g(n+1)\right) $,
\end{center}

\noindent which contradicts the foundation axiom in $\mathcal{N}^{+}$. This
completes the proof of Theorem 2.4.8, which in turn concludes the proof of
part (c) of Theorem 2.4.7.\hfill $\square $\medskip

\noindent \textbf{2.4.9.~Theorem.}~(Different faces of weak compactness of
Ord) \textit{The following are equivalent for any model }$(\mathcal{M},%
\mathfrak{X})$\textit{\ of }$\mathrm{GBC}$:\medskip

\noindent $(i)$ (\textit{Tree property}) $(\mathcal{M},\mathfrak{X})\models
\psi _{1},$ \textit{where} $\psi _{1}$ \textit{expresses}: \textit{Every\ }$%
\mathit{\mathrm{Ord}}$\textit{-tree has a branch}.\medskip

\noindent $(ii)$ (\textit{Weak compactness})\textrm{\ }$(\mathcal{M},%
\mathfrak{X})\models \psi _{2}$, \textit{where} $\psi _{2}$ \textit{expresses%
}:\textit{\ For any language} $\mathbb{L}$\textrm{,\textit{\ if}\ }$T$%
\textit{\ is an}\textrm{\ }$\mathbb{L}_{\mathrm{\infty },\mathrm{\infty }}$%
\textrm{-\textit{theory of cardinality} Ord} \textit{such that every
set-sized subtheory of} $T$ \textit{has a model, then there is a full
satisfaction class for a model of} $T$. \medskip

\noindent $(iii)$ (\textit{Ramsey property for an arbitrary set of colors in}
$\mathcal{M}$ \textit{and an arbitrary metatheoretic exponent }$n\geq 2$) $(%
\mathcal{M},\mathfrak{X})\models \psi _{3,n}$, \textit{where} $n\geq 2$
\textit{and} $\psi _{3,n}$ \textit{expresses}: $\forall \kappa \left(
\mathrm{Ord}\rightarrow \left( \mathrm{Ord}\right) _{\kappa }^{n}\right) .$%
\footnote{%
As shown in Theorem 4.9, this result cannot be strengthened by quantifying
over $n$ \textit{within} the theory \textrm{GBC} $+$ \textit{%
\textquotedblleft }\textrm{Ord is weakly compact}\textit{\textquotedblright .%
}}\medskip

\noindent $(iv)$ (\textit{Ramsey property for exponent }$2$\textit{\ and }$2$%
\textit{\ colors}) $(\mathcal{M},\mathfrak{X})\models \psi _{4}$, \textit{%
where} $\psi _{4}$ \textit{expresses}: $\mathrm{Ord}\rightarrow \left(
\mathrm{Ord}\right) _{2}^{2}.$\medskip

\noindent $(v)$ (\textit{Keisler property}) $(\mathcal{M},\mathfrak{X}%
)\models \psi _{5},$ \textit{where} $\psi _{5}$ \textit{expresses}$:$
\textit{For all }$X$ \textit{there is some} $S$ \textit{such that} $S$
\textit{is an} $\mathbb{L}_{\infty ,\infty }$-\textit{satisfaction class for
an\ }$\mathbb{L}_{\infty ,\infty }$\textrm{-}\textit{e.e.e. of} $(\mathrm{V}%
,\in ,X)$.\medskip

\noindent $(vi)$ ($\Pi _{1}^{1}$-\textit{Reflection}) \textit{For every }$%
\mathcal{L}_{\mathrm{Set}}(X,Y)$-\textit{formula }$\varphi (X,Y,x)$, \textit{%
and for each} $m\in M$ \textit{and} $A\in \mathfrak{X}$, $(\mathcal{M},%
\mathfrak{X})$ \textit{satisfies the following sentence in which }$A_{\alpha
}:=A\cap \mathrm{V}(\alpha )$:

\begin{center}
$\left[ \forall X\ \varphi (X,A,m)\right] \longrightarrow $\medskip

$\left[ \exists \alpha \ \forall X\ \subseteq \mathrm{V}(\alpha )\ (\mathrm{V%
}(\alpha ),\in ,X,A_{\alpha })\models \varphi (X,A_{\alpha}
,m)\right] .$%
\medskip
\end{center}

\noindent \textbf{Proof.}~With the help of Theorem 2.4.5, the equivalence of
$(i)$, $(ii)$, $(iii)$, and $(iv)$ can be verified with the same strategy as
in the usual ZFC-proofs (e.g., as in \cite[Theorem 7.8]{Kanamori book}) of
the equivalence of various formulations of weak compactness of a cardinal.
It is easy to see that $(v)\Rightarrow (i)$. To show the equivalence of $(v)$
with any of $(i)$ through $(iv)$, however, takes much more effort in
contrast to the ZFC-setting, e.g., in order to show that $(ii)\Rightarrow
(v) $ one first needs to know that if $(\mathcal{M},\mathfrak{X})$ is a
model of $\mathrm{GBC}$ in which $(ii)$ holds, and $X\in \mathfrak{X}$, then
the\textrm{\ }$\mathcal{L}_{\infty ,\infty }$-elementary diagram\textrm{\ }%
of $(\mathcal{M},X)$ is available as a member of $\mathfrak{X}$ (where $%
\mathcal{L}=\mathcal{L}_{\mathrm{Set}}(X)).$ More officially, we need to
know that $\mathfrak{X}$ contains an $\mathcal{L}_{\infty ,\infty }$%
-satisfaction class for $(\mathcal{M},X)\mathfrak{\ }$(as defined in
Definition 2.3.2(d)). This is precisely where part (c) of Theorem 2.4.7%
\textbf{\ }comes to the rescue. With the equivalence of $(v)$ with each of $%
(i)$ through $(iv)$ at hand, the proof will be complete once we show that $%
(v)\Rightarrow (vi)\Rightarrow (i)$. To see that $(v)\Rightarrow (vi),$
suppose $(\mathcal{M},\mathfrak{X})$ is a model of $\mathrm{GBC}$ in which $%
(v)$ holds, and suppose $(\mathcal{M},\mathfrak{X})\models \forall X\
\varphi (X,A,m)$ for some $m\in M$ and $A\in \mathfrak{X}$. Let $\left(
\mathcal{N},B\right) $ be an $\mathcal{L}_{\infty ,\infty }$-elementary end
extension of $\left( \mathcal{M},A\right) $, where for some $S\in \mathfrak{X%
},$ $S$ is a $\mathcal{L}_{\infty ,\infty }$-satisfaction class for $\left(
\mathcal{M},A\right) $. Recall that in ZFC the well-foundedness of $\in $ is
expressible in $\mathcal{L}_{\omega _{1},\omega _{1}}$ via the sentence $%
\psi $ below:

\begin{center}
$\psi :=\lnot \exists \left\langle x_{n}:n\in \omega \right\rangle
\bigwedge\limits_{n\in \omega }x_{n+1}\in x_{n}.$
\end{center}

\noindent Therefore, since $\mathcal{M}$ satisfies $\psi $ and $\mathcal{N}$
is an $\mathcal{L}_{\infty ,\infty }$-elementary extension of $\mathcal{M},$
\textrm{Ord}$^{\mathcal{N}}\ \backslash \ \mathrm{Ord}^{\mathcal{M}}$ has a
minimum element $\kappa $, and thus $\left( \mathcal{N}(\kappa ),B\cap
N(\kappa )\right) =\left( \mathcal{M},A\right) $. Hence:

\begin{center}
$\left( \mathcal{N},B\right) \models \forall X\subseteq \mathrm{V}(\kappa )\
\varphi ^{\mathrm{V}(\kappa )}(X,A,m),$
\end{center}

\noindent where $\varphi ^{\mathrm{V}(\kappa )}$ is the result of
restricting the (set) quantifiers of $\varphi $ to $\mathrm{V}(\kappa ).$
Therefore since $\left( \mathcal{M},A\right) \prec \left( \mathcal{N}%
,B\right) $, we conclude:

\begin{center}
$\left( \mathcal{M},A\right) \models \exists \alpha \ \forall X\subseteq
\mathrm{V}(\alpha )\ \varphi ^{\mathrm{V}(\alpha )}(X,A,m)$,
\end{center}

\noindent thus completing the proof of $(v)\Rightarrow (vi)$. The proof of $%
(vi)\Rightarrow (i)$ is routine and uses the standard strategy of
showing within $\mathrm{ZFC}$ that the $\Pi _{1}^{1}$-Reflection property of
an inaccessible cardinal $\kappa $ implies that $\kappa $ has tree
property.\hfill $\square \medskip $

Recall that the notion \textquotedblleft $\kappa $ is $\alpha $%
-Mahlo\textquotedblright\ is defined recursively by decreeing that
\textquotedblleft $\kappa $ is $0$-Mahlo\textquotedblright\ means that $%
\kappa $ is strongly inaccessible, and for an ordinal $\alpha >0$
\textquotedblleft $\kappa $ is $\alpha $-Mahlo\textquotedblright\ means that
for all $\beta <\alpha $ the collection of cardinals that are $\beta $-Mahlo
are stationary in $\kappa .$ It is a classical fact that, provably in ZFC,
every weakly compact cardinal $\kappa $ is $\kappa $-Mahlo.\medskip

Theorem 2.4.12 below summarizes some well-known facts about the Levy scheme $%
\Lambda $; the statement of the theorem uses the following
definition.\medskip

\noindent \textbf{2.4.10.~Definition.}~In what follows $X$\ is a unary
predicate symbol (which will be conflated with its interpretation in a given
structure).\medskip

\noindent \textbf{(a) }$\Lambda =\{\lambda _{n}:n\in \omega \}$, where $%
\lambda _{n}$ is the $\mathcal{L}_{\mathrm{Set}}$-sentence asserting the
existence of an $n$-Mahlo cardinal $\kappa $ such that $\left( \mathrm{V}%
(\kappa ),\in \right) \prec _{\Sigma _{n}}\left( \mathrm{V},\in \right) .$
More generally, $\Lambda (X)=\{\lambda _{n}(X):n\in \omega \}$, and $\lambda
_{n}(X)$ is the $\mathcal{L}_{\mathrm{Set}}(X)$-sentence asserting the
existence of an $n$-Mahlo cardinal $\kappa $ such that $\left( \mathrm{V}%
(\kappa ),\in ,X\cap \mathrm{V}(\kappa )\right) \prec _{\Sigma
_{n}(X)}\left( \mathrm{V},\in ,X\right) .$ \medskip

\noindent \textbf{(b) }For $n\in \omega ,$ $\Lambda _{n}(X)=\{\lambda
_{n,i}(X):i\in \omega \}$, and $\lambda _{n,i}(X)$ is the sentence asserting
the existence of an $n$-Mahlo cardinal $\kappa $ such that $\left( \mathrm{V}%
(\kappa ),\in ,X\cap \mathrm{V}(\kappa )\right) \prec _{\Sigma
_{i}(X)}\left( \mathrm{V},\in ,X\right) .$\medskip

\noindent \textbf{(c) }$\Lambda ^{-}$ is the fragment of $\Lambda $
consisting of statements of the form \textquotedblleft there is an $n$-Mahlo
cardinal\textquotedblright, for $n\in \omega .$\medskip

\noindent \textbf{2.4.11.~Theorem.}~(Folklore).\medskip

\noindent \textbf{(a)} \textit{For }$n\in \omega ,$ $\kappa $\textit{\ is} $%
(n+1)$-\textit{Mahlo iff} \textit{for every} $X\subseteq \mathrm{V}(\kappa )$%
,$\ (\mathrm{V}(\kappa ),\in ,X)\models \Lambda _{n}(X).$\medskip

\noindent \textbf{(b) }$\kappa $\textit{\ is} $\omega $-\textit{Mahlo iff
for every} $X\subseteq \mathrm{V}(\kappa )$,$\ (\mathrm{V}(\kappa ),\in
,X)\models \Lambda (X).$\medskip

\noindent \textbf{(c)} $\mathrm{ZFC}~+~\Lambda ^{-}$ \textit{is mutually
interpretable with} $\mathrm{ZFC}~+~\Lambda .$\medskip

\noindent \textbf{(d)} \textit{Assuming the consistency of} $\mathrm{ZFC}+$%
\textquotedblleft \textit{there is an }$\omega $-\textit{Mahlo cardinal}%
\textquotedblright $,$ $\mathrm{ZFC}~+~\Lambda ^{-}\nvdash \Lambda .$\medskip

\noindent \textbf{(e)} \textit{If} $\mathcal{M}\models \mathrm{ZFC}+\Lambda $%
, \textit{then} $\mathrm{L}^{\mathcal{M}}\models \Lambda $ (\textit{where} $%
\mathrm{L}^{\mathcal{M}}$ \textit{is the constructible universe of} $%
\mathcal{M)}$. \medskip

\noindent \textbf{(f)} \textit{If} $\mathcal{M}\models \mathrm{ZFC}~+~\Lambda $%
, $\mathbb{P}$\textit{\ is a} \textit{set notion of forcing} $\mathbb{P}$
\textit{in} $\mathcal{M}$, \textit{and} $G$\textit{\ is} $\mathbb{P}$-%
\textit{generic over} $\mathcal{M}$, \textit{then} $\mathcal{M}[G]\models
\Lambda $.\medskip

\noindent \textbf{Proof.} Suppose $\kappa $ is a strongly inaccessible
cardinal and $X\subseteq \mathrm{V}(\kappa ).$ Let

\begin{center}
$C:=\left\{ \lambda \in \kappa :\left( \mathrm{V}(\lambda ),\in ,X\cap
\mathrm{V}(\lambda )\right) \prec \left( \mathrm{V}(\kappa ),\in ,X\right)
\right\} .$
\end{center}

\noindent A routine Skolem hull argument shows that $C$ is closed and
unbounded in $\kappa .$ This fact lies at the heart of the proofs of (a)
through (c); note that the proof of (c) uses Orey's Compactness Theorem 5.3. To verify (d), work in a model of $\mathrm{ZFC}$ $+\text{
\textquotedblleft }$there is an $\omega $-Mahlo cardinal$\text{%
\textquotedblright },$ and for each $n\in \omega $ let $\kappa _{n}$ be the
first $n$-Mahlo cardinal, and $\kappa _{\omega }:=\sup\limits_{n\in \omega
}\kappa _{n}$. Choose the first strongly inaccessible cardinal $\lambda
>\kappa _{\omega }$. Then $\mathrm{ZFC}+\Lambda ^{-}$ clearly holds in $%
\left( \mathrm{V}(\lambda ),\in \right) $. We will show that $\Lambda $
fails in $\left( \mathrm{V}(\lambda ),\in \right) .$ To see this, we first
note:\medskip

\noindent $(\ast )$ $\left( \mathrm{V}(\lambda ),\in \right) \models $
\textquotedblleft the collection of Mahlo cardinals is bounded in $\mathrm{%
Ord}$\textquotedblright .\medskip

\noindent On the other hand, ``there are unboundedly many $n$%
-Mahlo cardinals in the universe'' holds in $\left( \mathrm{%
V}(\kappa ),\in \right) $ for any $(n+1)$-Mahlo cardinal $\kappa $, and
therefore if $\left( \mathrm{V}(\kappa ),\in \right) \prec _{\Sigma _{m}}(%
\mathrm{V},\in )$, where the statement $\varphi =$ ``there are unboundedly
many $n$-Mahlo cardinals'' is a $\Sigma _{m}$-sentence, then $\varphi $ holds
in the universe. Together with $(\ast )$, this makes it clear that $\Lambda $
fails in $\left( \mathrm{V}(\lambda ),\in \right) .$ Part (e) follows from
routine absoluteness considerations, and part (f) is a consequence of the
preservation of both (1) the $n$-Mahlo property of a cardinal $\kappa $, and
(2) the property $\mathrm{V}(\kappa )\prec _{\Sigma _{n}}\mathrm{V}$, in $%
\mathbb{P}$-generic extensions satisfying $\mathbb{P} \in \mathrm{V}(\kappa )$. (1) is
established along the lines of the proof of \cite[Proposition 10.13]%
{Kanamori book}; (2) follows from a standard truth-and-forcing
argument.\hfill $\square \medskip $

The theorem below reveals the close relationship between the class theory $%
\mathrm{GBC}+$ \textquotedblleft $\mathrm{Ord}$\textrm{\ is weakly compact}%
\textquotedblright ,\ and the set theory $\mathrm{ZFC}+\Lambda .$\medskip

\noindent \textbf{2.4.12.~Theorem.}~\cite[Corollary 2.1.1]{Ali NFUA}~\textit{%
Let }$\varphi $\textit{\ be an }$\mathcal{L}_{\mathrm{Set}}$\textit{-sentence%
}$.$ \textit{The} \textit{following are equivalent}:\medskip

\noindent $(i)$ $\mathrm{GBC}\ +\ $\textrm{\textquotedblleft Ord}$\ \mathrm{%
is\ weakly\ compact}$\textquotedblright $\ \vdash \varphi .$\medskip

\noindent $(ii)$\textbf{\ }$\mathrm{ZFC}+\Lambda \vdash \varphi .$\bigskip

\begin{center}
\textbf{3.~BASIC\ FEATURES\ OF\ }$\mathbf{ZFI\ }$\textbf{AND }$\mathbf{ZFI}_{%
\mathrm{<}}$\bigskip
\end{center}

In this section we officially meet the principal characters of our paper,
namely the theory $\mathrm{ZFI}$, and its extension $\mathrm{ZFI}_{\mathrm{<}%
}.$ We establish two useful schemes (apartness and diagonal
indiscernibility) within $\mathrm{ZFI}$. These schemes are then used to
demonstrate some basic model-theoretic facts about \textrm{ZFI}\ and $%
\mathrm{ZFI}_{\mathrm{<}}$. In particular, we show that $\omega $%
-nonstandard models of $\mathrm{ZF}$ that have an expansion to $\mathrm{ZFI}$%
\ are recursively saturated, and $\omega $-standard models of ZF that have
an expansion to ZFI$_{<}$ satisfy \textquotedblleft 0$^{\#}$
exists\textquotedblright .\medskip

\noindent \textbf{3.1.~Definition.}~$\mathrm{ZFI}$ is the theory formulated
in the language $\mathcal{L}_{\mathrm{Set}}(I)$, where $I$ is a unary
predicate, whose axioms consist of the three groups below.

\begin{itemize}
\item Note that we often write $x\in I$ instead of $I(x)$.\medskip
\end{itemize}

\noindent $(1)$ $\mathrm{ZF(}I\mathrm{)}$. Recall from part (e) of
Definition 2.1.1 that $\mathrm{ZF(}I\mathrm{)}$ includes the separation
scheme $\mathrm{Sep(}I)$ and the collection scheme $\mathrm{Coll(}I)\mathrm{.%
}$\medskip

\noindent $(2)$ The sentence $\mathrm{Cof}(I)$ expressing \textquotedblleft $%
I$\ is a cofinal subclass \textrm{Ord}\textquotedblright .\medskip

\noindent $(3)$ The scheme $\mathrm{Indis}(I)=\{\mathrm{Indis}_{\varphi
}(I): $ $\varphi $ is a formula of $\mathcal{L}_{\mathrm{Set}}\}$ ensuring
that $I$ forms a class of order indiscernibles for the ambient model $(%
\mathrm{V},\in )$ of set theory$\mathrm{.}$ More explicitly, for each $n$%
-ary formula $\varphi (v_{1},\cdot \cdot \cdot ,v_{n})$ in the language $%
\{=,\in \},$ $\mathrm{Indis}_{\varphi }(I)$ is the sentence:\medskip

\begin{center}
$\forall x_{1}\in I\cdot \cdot \cdot \forall x_{n}\in I$ $\forall y_{1}\in
I\cdot \cdot \cdot \ \forall y_{n}\in I$ \medskip

$[(x_{1}\in \cdot \cdot \cdot \in x_{n})\wedge (y_{1}\in \cdot \cdot \cdot
\in y_{n})\rightarrow (\varphi (x_{1},\cdot \cdot \cdot
,x_{n})\leftrightarrow \varphi (y_{1},\cdot \cdot \cdot ,y_{n}))].$\medskip
\end{center}

\noindent The theory $\mathrm{ZFI}_{\mathrm{<}}$ is an extension of $\mathrm{%
ZFI}$; it is formulated in the language $\mathcal{L}_{\mathrm{Set}}(I,<)$,
whose axioms consist of $\mathrm{Cof}(I)$ above, together with the following
strengthenings of the axioms in (1) and (3) above:\medskip

\noindent $(1^{+})$ $\mathrm{ZF(}I,<\mathrm{)+GW}$.\medskip

\noindent $(3^{+})$ The scheme $\mathrm{Indis}_{<}(I)=\{\mathrm{Indis}%
_{\varphi }(I):$ $\varphi $ is a formula of $\mathcal{L}_{\mathrm{Set}}(<)\}$
ensuring that $I$ forms a class of order indiscernibles for $(\mathrm{V},\in
,<).$

\begin{itemize}
\item The above definition can be model-theoretically recast as follows: $%
\mathcal{M}\models \mathrm{ZF}$ has an expansion $(\mathcal{M},I)\models
\mathrm{ZFI}$ iff there is an $\mathcal{M}$-amenable cofinal subset $I$ of
\textrm{Ord}$^{\mathcal{M}}$ such that $(I,\in _{M})$ forms a class of
indiscernibles over $\mathcal{M}$. Similarly, a model $\left( \mathcal{M}%
,<_{M}\right) \models \mathrm{ZF(}\mathcal{<}\mathrm{)+GW}$ has an expansion
$\left( \mathcal{M},<_{M},I\right) \models \mathrm{ZFI}_{\mathrm{<}}$ iff
there is an $(\mathcal{M},<_{M})$-amenable cofinal subset $I$ of \textrm{Ord}%
$^{\mathcal{M}}$ such that $(I,\in _{M})$ forms a class of indiscernibles
over $\left( \mathcal{M},<_{M}\right) $. Therefore by Theorem 2.1.5 if $%
\mathcal{M}\models \mathrm{ZF}+\mathrm{V=HOD}$, and $\mathcal{M}$ has an
expansion to a model of $\mathrm{ZFI}$, then $\mathcal{M}$ is also
expandable to a model of $\mathrm{ZFI}_{<}$. Moreover, by Theorem 3.2(b)
below, the assumption that $\mathcal{M}\models \mathrm{V=HOD}$ can be
weakened to the assumption that $\mathcal{M}\models \exists p(\mathrm{V=HOD(}%
p\mathrm{)})$.
\end{itemize}

\noindent \textbf{3.2.~Theorem. }\textit{Let }\textrm{ZFI}$^{\ast }$ \textit{%
be the subsystem of }\textrm{ZFI} \textit{axiomatized by} $\mathrm{ZF+Coll%
\mathrm{(}}I\mathrm{)}+\mathrm{\mathrm{Cof}(}I\mathrm{)}+\mathrm{%
\mathrm{Indis}(}I\mathrm{)}$. \textit{The following schemes are provable in}
\textrm{ZFI}$^{\ast }$:\medskip

\noindent \textbf{(a) }\textit{The \textbf{apartness} scheme for }$\mathcal{L}_{%
\mathrm{Set}}$\textit{-formulae}:\medskip

\begin{center}
$\mathrm{Apart}(I) =\{ \mathrm{Apart}_{\varphi }(I):\varphi \in \mathrm{Form}_{n+1}\mathrm{(}%
\mathcal{L}_{\mathrm{Set}})$, $n\in \omega \mathrm{\},}$\medskip
\end{center}

\noindent \textit{where} $\mathrm{Form}_{n}\mathrm{(}\mathcal{L}_{\mathrm{Set%
}})$ \textit{is the collection of} $\mathcal{L}_{\mathrm{Set}}$-\textit{%
formulae whose free variables are} $x_{1},\cdot \cdot \cdot ,x_{n},$ \textit{%
and} $\mathrm{Apart}_{\varphi }(I)$ \textit{is the following formula}:\medskip

\begin{center}
$\forall i\in I\ \forall j\in I\left[ i<j\rightarrow \forall \overline{x}\in
\left( \mathrm{V}(i)\right) ^{n}\left( \exists y\varphi (\overline{x}%
,y)\rightarrow \exists y\in \mathrm{V}(j)\ \varphi (\overline{x},y)\right) %
\right] .$\medskip
\end{center}

\noindent \textbf{(b)} \textit{The \textbf{diagonal indiscernibility} scheme for }$%
\mathcal{L}_{\mathrm{Set}}$\textit{-formulae}:

\begin{center}
$\mathrm{Indis}^{+}(I)= \{\mathrm{Indis}_{\varphi }^{+}(I):$ $\varphi \in
\mathrm{Form}_{n+1+r}\mathrm{(}\mathcal{L}_{\mathrm{Set}}),\ n,r\in \omega
,\ r\geq 1\}$,
\end{center}

\noindent \textit{where} $\mathrm{Indis}_{\varphi }^{+}(I)$ \textit{is the
following formula}:

\begin{center}
$\forall i\in I\ \forall \overline{j}\in \lbrack I]^{r}\ \forall \overline{k}%
\in \lbrack I]^{r}\ \left[ \left( i<j_{1}\right) \wedge (i<k_{1})\right]
\longrightarrow \medskip $

$\left[ \forall \overline{x}\in \left( \mathrm{V}(i)\right) ^{n}\ \left(
\varphi (\overline{x},i,j_{1},\cdot \cdot \cdot ,j_{r})\leftrightarrow
\varphi (\overline{x},i,k_{1},\cdot \cdot \cdot ,k_{r})\right) \right] .$%
\medskip
\end{center}

\noindent \textit{Similarly, let} \textrm{ZFI}$_{<}^{\ast }$ \textit{be the
subsystem of }\textrm{ZFI}$_{<}$ \textit{axiomatized by} $\mathrm{ZF(<)}+%
\mathrm{GW}+\mathrm{Coll\mathrm{(}}I\mathrm{)}+\mathrm{\mathrm{Cof}(}I%
\mathrm{)}+\mathrm{Indis}_{<}(I)$. \textit{The following schemes are
provable in} \textrm{ZFI}$_{<}^{\ast }$:\medskip

\noindent \textbf{(c) }\textit{The \textbf{apartness} scheme for }$\mathcal{L%
}_{\mathrm{Set}}(<)$\textit{-formulae}:\medskip

\begin{center}
$\mathrm{Apart}_{<}(I) = \{\mathrm{Apart}_{\varphi }(I):\varphi \in \mathrm{Form}%
_{n+1}\mathrm{(}\mathcal{L}_{\mathrm{Set}}(<))$, $n\in \omega \mathrm{\},}$%
\medskip
\end{center}

\noindent \textit{where} $\mathrm{Form}_{n}\mathrm{(}\mathcal{L}_{\mathrm{Set%
}}(<))$ \textit{is the collection of} $\mathcal{L}_{\mathrm{Set}}(<)$-%
\textit{formulae whose free variables are} $x_{1},\cdot \cdot \cdot ,x_{n},$
\textit{and} $\mathrm{Apart}_{\varphi }(I)$ \textit{is the following formula}%
:\medskip

\begin{center}
$\forall i\in I\ \forall j\in I\left[ i<j\rightarrow \forall \overline{x}\in
\left( \mathrm{V}(i)\right) ^{n}\left( \exists y\varphi (\overline{x}%
,y)\rightarrow \exists y\in \mathrm{V}(j)\ \varphi (\overline{x},y)\right) %
\right] .$\medskip
\end{center}

\noindent \textbf{(d)} \textit{The \textbf{diagonal indiscernibility} scheme for}$%
\mathcal{L}_{\mathrm{Set}}(<)$\textit{-formulae}:

\begin{center}
$\mathrm{Indis}_{<}^{+}(I)=\{\mathrm{Indis}_{\varphi }^{+}(I):$ $\varphi \in
\mathrm{Form}_{n+1+r}\mathrm{(}\mathcal{L}_{\mathrm{Set}}(<)),\ n,r\in
\omega ,\ r\geq 1\}$,
\end{center}

\noindent \textit{where} $\mathrm{Indis}_{\varphi }^{+}(I)$ \textit{is the
following formula}:

\begin{center}
$\forall i\in I\ \forall \overline{j}\in \lbrack I]^{r}\ \forall \overline{k}%
\in \lbrack I]^{r}\ \left[ \left( i<j_{1}\right) \wedge (i<k_{1})\right]
\longrightarrow \medskip $

$\left[ \forall \overline{x}\in \left( \mathrm{V}(i)\right) ^{n}\ \left(
\varphi (\overline{x},i,j_{1},\cdot \cdot \cdot ,j_{r})\leftrightarrow
\varphi (\overline{x},i,k_{1},\cdot \cdot \cdot ,k_{r})\right) \right] .$
\end{center}

\noindent \textbf{Proof.}~We will only establish (a) and (b) since the proof
of (c) is similar to the proof of (a) and the proof of (d) is similar to the
proof of (b). Let $(\mathcal{M},I)\models \mathrm{ZFI}^{\ast }\mathrm{.}$ To
verify that the apartness scheme holds in $(\mathcal{M},I)$, fix some $%
i_{0}\in I$ and some $\varphi (\overline{x},y)\in \mathrm{Form}_{n+1}\mathrm{%
(}\mathcal{L}_{\mathrm{Set}}).$ Then, since the collection scheme $\mathrm{%
Coll}(I)$ holds in $(\mathcal{M},I),$ and $I$ is cofinal in \textrm{Ord}$^{%
\mathcal{M}}$, there is some $j_{0}\in I$ with $i_{0}<j_{0}$ such that:

\begin{center}
$(\mathcal{M},I)\models \forall \overline{x}\in \left( \mathrm{V}%
(i_{0})\right) ^{n}\ \left( \exists y\varphi (\overline{x},y)\rightarrow
\exists y\in \mathrm{V}(j_{0})\ \varphi (\overline{x},y)\right) .$
\end{center}

\noindent The above, together with the indiscernibility of $I$ in $\mathcal{M%
}$, makes it evident that $(\mathcal{M},I)\models \mathrm{Apart}_{\varphi }.$
\medskip

To verify that $\mathrm{Indis}_{\varphi }^{+}(I)$ holds in $(\mathcal{M},I)$%
, we will first establish a weaker form of diagonal indiscernibility of $I$
in which all $j_{n}<k_{1}$ (thus all the elements of $\overline{j}$ are less
than all the elements of $\overline{k}$)$.$ Fix some $\varphi \in \mathrm{%
Form}_{n+1+r}\mathrm{(}\mathcal{L}_{\mathrm{Set}})$ and $i_{0}\in I.$ Within
$\mathcal{M}$ consider the function $f:\left[ \mathrm{Ord}\right]
^{r}\rightarrow \mathcal{P}(\mathrm{V}(i_{0})^{n})$ by:

\begin{center}
$f(\overline{\gamma }):=\{\overline{a}\in \left( \mathrm{V}(i_{0})\right)
^{n}:\varphi (\overline{a},i_{0},\overline{\gamma })\}.$
\end{center}

\noindent Since $(\mathcal{M},I)$ satisfies the collection scheme $\mathrm{%
Coll}(I)$ and $I$ is cofinal in \textrm{Ord}$^{%
\mathcal{M}}$, there are $\mathcal{M}$-ordinals $\gamma
_{1}<\cdot \cdot \cdot <\gamma _{2r}$ in $I$ such that:

\begin{center}
$f(\gamma _{1},\cdot \cdot \cdot ,\gamma _{r})=f(\gamma _{r+1},\cdot \cdot
\cdot ,\gamma _{2r}).$
\end{center}

\noindent Thus we have:

\begin{center}
$(\mathcal{M},I)\models \left[ \forall \overline{x}\in \left( \mathrm{V}%
(i_{0})\right) ^{n}\ \left( \varphi (\overline{x},i_{0},\gamma _{1},\cdot
\cdot \cdot ,\gamma _{r})\leftrightarrow \varphi (\overline{x},i_{0},\gamma
_{r+1},\cdot \cdot \cdot ,\gamma _{2r})\right) \right] .$
\end{center}

\noindent By indiscernibility of $I$ in $\mathcal{M}$, the above implies the
following weaker form of $\mathrm{Indis}_{\varphi }^{+}(I)$:\medskip

\begin{center}
$\forall i\in I\ \forall \overline{j}\in \lbrack \mathrm{I}]^{r}\ \forall
\overline{k}\in \lbrack I]^{r}\ \left[ \left( i<j_{1}\right) \wedge
(j_{n}<k_{1})\right] \longrightarrow $\medskip

$\left[ \forall \overline{x}\in \left( \mathrm{V}(i)\right) ^{n}\ \left(
\varphi (\overline{x},i,j_{1},\cdot \cdot \cdot ,j_{r})\leftrightarrow
\varphi (\overline{x},i,k_{1},\cdot \cdot \cdot ,k_{r})\right) \right] .$%
\medskip
\end{center}

\noindent We will now show that the above weaker form of $\mathrm{Indis}%
_{\varphi }^{+}(I)$ implies $\mathrm{Indis}_{\varphi }^{+}(I).$ Given $i\in
I $, $\overline{\alpha }\in \lbrack I]^{r}\ $and $\overline{\beta }\in
\lbrack I]^{r},$ with $i<\alpha _{1}$ and $i<\beta _{1},$ choose $\overline{%
\gamma }\in \lbrack I]^{r}$ with $\gamma _{1}>\max \left\{ \alpha _{n},\beta
_{n}\right\} .$ Then by the above we have:

\begin{center}
$\mathcal{M}\models \left[ \forall \overline{x}\in \left( \mathrm{V}%
(i)\right) ^{n}\ \left( \varphi (\overline{x},i,\alpha _{1},\cdot \cdot
\cdot ,\alpha _{r})\leftrightarrow \varphi (\overline{x},i,\gamma _{1},\cdot
\cdot \cdot ,\gamma _{r})\right) \right] ,$
\end{center}

\noindent and

\begin{center}
$\mathcal{M}\models \left[ \forall \overline{x}\in \left( \mathrm{V}%
(i)\right) ^{n}\ \left( \varphi (\overline{x},i,\beta _{1},\cdot \cdot \cdot
,\beta _{r})\leftrightarrow \varphi (\overline{x},i,\gamma _{1},\cdot \cdot
\cdot ,\gamma _{r})\right) \right] ,$
\end{center}

\noindent which together imply:

\begin{center}
$\mathcal{M}\models \left[ \forall \overline{x}\in \left( \mathrm{V}%
(i)\right) ^{n}\ \left( \varphi (\overline{x},i,\alpha _{1},\cdot \cdot
\cdot ,\alpha _{r})\leftrightarrow \varphi (\overline{x},i,\beta _{1},\cdot
\cdot \cdot ,\beta _{r})\right) \right] .$

\hfill $\square $\medskip
\end{center}

\begin{itemize}
\item Note that the diagonal indiscernibility scheme for $\mathcal{L}_{%
\mathrm{Set}}$-formulae ensures that if $(\mathcal{M},I)\models \mathrm{ZFI}$
and $i\in I,$ then $I^{\geq i}$ is a set of indiscernibles over the expanded
structure $(\mathcal{M},m)_{m\in \mathrm{V}(i)}$, where $I^{\geq i}=\{j\in
I:j\geq i\}.$ Similarly, the diagonal indiscernibility scheme for $\mathcal{L}_{%
\mathrm{Set}}(<)$-formulae ensures that if $(\mathcal{M},<,I)\models \mathrm{%
ZFI}_{<}$ and $i\in I,$ then $I^{\geq i}$ is a set of indiscernibles over
the expanded structure $(\mathcal{M},<_{M},m)_{m\in \mathrm{V}(i)}.$
\end{itemize}

The fact that the apartness scheme holds in $\mathrm{ZFI}$ and $\mathrm{ZFI}%
_{<}$ will be employed in the following theorem to show that \textrm{ZFI} is
able to define a $\Sigma _{\omega }$-satisfaction predicate over the ambient
model of $\mathrm{ZF}$, and \textrm{ZFI}$_{<}$ is able to define a $\Sigma
_{\omega }$-satisfaction predicate over the ambient model of $\mathrm{%
ZF(<)+GW}$ (in the sense of part (b) of Definition 2.3.2).\medskip

\noindent \textbf{3.3.~Theorem.}~\textit{There is a formula }$\sigma (x)$
\textit{in the language} $\mathcal{L}_{\mathrm{Set}}(I)$ \textit{such that
for all models} $(\mathcal{M},I)$\textit{\ of }$\mathrm{ZFI,}$ $\sigma ^{%
\mathcal{M}}$ \textit{is a }$\Sigma _{\omega }$-\textit{satisfaction class}
\textit{for} $\mathcal{M}$. \textit{In particular}:\medskip

\noindent \textbf{(a)}
\textit{If} $(\mathcal{M},I)\models \mathrm{ZFI}$, \textit{then} $\sigma ^{\mathcal{M}}$ \textit{is an amenable $\Sigma _{\omega }$-satisfaction class for} $%
\mathcal{M}$.

\medskip

\noindent \textbf{(b)} \textit{If} $(\mathcal{M},I)\models \mathrm{ZFI}$ ,
\textit{and }$\mathcal{M}$\textit{\ is }$\omega $-\textit{standard, then} $%
\sigma ^{\mathcal{M}}$ \textit{is an amenable full satisfaction class for} $%
\mathcal{M}$.\medskip

\noindent \textit{Similarly}, \textit{there is a formula }$\sigma _{<}(x)$
\textit{in the language} $\mathcal{L}_{\mathrm{Set}}(<,I)$ \textit{such that
for all models} $(\mathcal{M},<,I)$\textit{\ of }$\mathrm{ZFI}_{<}\mathrm{,}$
$\sigma _{<}^{\mathcal{M}}$ \textit{is a }$\Sigma _{\omega }$-\textit{%
satisfaction class} \textit{for} $\left( \mathcal{M},<_{M}\right) $. \textit{%
In particular}:\medskip

\noindent \textbf{(c)} \textit{If} $(\mathcal{M},<_{M},I)\models \mathrm{ZFI}_{<}$ , \textit{then} $\sigma _{<}^{\mathcal{M}}$ \textit{is an amenable $\Sigma _{\omega }$-satisfaction class for} $\left( \mathcal{M},<_{M}\right) $.

\medskip

\noindent \textbf{(d)} \textit{If} $(\mathcal{M},<_{M},I)\models \mathrm{ZFI}%
_{<}$ , \textit{and }$\mathcal{M}$\textit{\ is }$\omega $-\textit{standard,
then} $\sigma _{<}^{\mathcal{M}}$ \textit{is an amenable full satisfaction
class for} $\left( \mathcal{M},<_{M}\right) $.\medskip

\noindent \textbf{Proof.}~(a) and (b) are immediate consequences of the
first assertion of the theorem, which we will establish. The proofs of (c)
and (d) are similar and will not be presented. The following definition
takes place in $(\mathcal{M},I)$: Given any $\varphi (\overline{x})\in
\mathrm{Form}_{k}(\mathcal{L}_{\mathrm{Set}})$ and any $k$-tuple $\overline{a%
}$, let $i_{0}$ be the first element of $I$ such that\textit{\ }$\overline{a}%
\in \mathrm{V}(i_{0}),$ and for each $n\in \omega ,$ let $i_{n+1}$ be the
first element of $I$ that exceeds $i_{n}.$ Then let $\alpha :=\sup_{n\in
\omega }i_{n}.$ \medskip

It is easy to see, by Tarski's test (for elementarity) and the veracity of
the Apartness scheme in $(\mathcal{M},I)$ that $\mathcal{M}(\alpha )\prec
\mathcal{M}.$ Therefore, if $S$ is defined in $\mathcal{M}$ by:

\begin{center}
$\varphi (\overline{a})\in S$ iff $(\mathrm{V}(\alpha ),\in )\models \varphi
(\overline{a}),$
\end{center}

\noindent then $S$ is a $\Sigma _{\omega }$-satisfaction class for $%
\mathcal{M}.$ Our description of $S$ makes it clear that $S$ is definable in
$\mathcal{M}$ by a parameter-free formula $\sigma (x)$ in the language $%
\mathcal{L}_{\mathrm{Set}}(I).$\hfill $\square $ \medskip

\noindent \textbf{3.4.~Remark.}~Theorem 2.4.7(c) together with the proof of
the $(iii)\Rightarrow (i)$ direction of Theorem 4.1 shows that if\textit{\ }$(%
\mathcal{M},I)\models \mathrm{ZFI}_{\mathrm{<}}$, then $\mathcal{M}_{I}$
carries an amenable full satisfaction class $S$ (but $S$ need not be definable in $\mathcal{M}_{I}$). On the other hand, it is known
\cite[Theorem 6.3]{Ali-Matt-Zach-Aut} that if $\mathcal{M}$ and $%
\mathcal{N}$ are models of $\mathrm{ZFC}$ such that $\mathcal{M}$\ is a
cofinal elementary submodel of $\mathcal{N}$, then for any $\mathcal{M}$%
-amenable subset $X_{M}$ of $M$, there is a (unique) subset $X_{N}$ of $N$ such that $(\mathcal{M},X_{M})\prec (\mathcal{N},X_{N}).$ Thus,
if\textit{\ }$(%
\mathcal{M},I)\models \mathrm{ZFI}_{\mathrm{<}}$, then $\mathcal{M}$
carries an amenable full satisfaction class.

\medskip

\noindent \textbf{3.5.~Corollary}. \textit{Suppose }$\mathcal{M}\models
\mathrm{ZF}$. \textit{There is no} \textit{parametrically }$\mathcal{M}$-%
\textit{definable subset} $I$ \textit{of} \textrm{Ord}$^{\mathcal{M}}$%
\textit{\ such that} $(\mathcal{M},I)\models \mathrm{ZFI}$. \textit{%
Similarly, if} $\mathcal{M}$ \textit{has an expansion} $\left( \mathcal{M}%
,<_{M}\right) \models \mathrm{ZF(<)}+\mathrm{GW,}$ \textit{then there is no parametrically }$\left( \mathcal{M},<_{M}\right) $-\textit{definable
subset} $I$ \textit{of} \textrm{Ord}$^{\mathcal{M}}$\textit{\ such that} $(%
\mathcal{M},<_{M},I)\models \mathrm{ZFI}_{<}.$\medskip

\noindent \textbf{Proof.} Put Theorem 3.3 together with Tarski's theorem on
undefinability of truth. Alternatively, one can take advantage of diagonal
indiscernibility.\hfill $\square $\medskip

\noindent \textbf{3.6.~Corollary}.~If $(\mathcal{M},I)\models \mathrm{ZFI},$
\textit{and} $\mathcal{M}$ \textit{is} $\omega $-\textit{nonstandard, then }$%
\mathcal{M}$ \textit{is recursively saturated}. \textit{Similarly, if }$(%
\mathcal{M},<_{M},I)\models \mathrm{ZFI}_{<},$ \textit{and} $\mathcal{M}$
\textit{is} $\omega $-\textit{nonstandard, then }$\left( \mathcal{M}%
,<_{M}\right) $ \textit{is recursively saturated.}\medskip

\noindent \textbf{Proof.} We will only verify the $\mathrm{ZFI}$ case; a
similar strategy works for $\mathrm{ZFI}_{<}\mathrm{.}$ This is established
using a well-known overspill argument using the fact that induction over $%
\omega ^{\mathcal{M}}$ holds in $(\mathcal{M},S)$, where $S$ is a $\Sigma
_{\omega }$-satisfaction class given by Theorem 3.3. More specifically,
since $S$ satisfies Tarski's compositional conditions for each $\Sigma _{n}$%
-formula (where $n\in \omega ),$ by overspill we can fix some nonstandard $%
c\in \omega ^{\mathcal{M}}$ such that $S$ satisfies Tarski conditions for $%
\Sigma _{c}$-formulae$.$ Next let $\left\langle \varphi _{i}(x):i\in \omega
\right\rangle $ be a recursive enumeration in the real world of the formulae
of a recursive type $p(x)$ (involving finitely many parameters from $%
\mathcal{M}$), where $p(x)$ is finitely realizable in $\mathcal{M}$. This
enumeration can be extended to some enumeration $\left\langle \varphi
_{i}(x):i\in \omega ^{\mathcal{M}}\right\rangle $ in $\mathcal{M}$. For each
$i\in \omega ^{\mathcal{M}}$ let

\begin{center}
$\psi _{i}:=\exists x\bigwedge\limits_{j\leq i}\varphi _{j}(x).$
\end{center}

\noindent Then for every $n\in \omega $, $(\mathcal{M},S)\models \theta (n)$%
, where $\theta (i):=\left( \psi _{i}\in \Sigma _{c}\right) \wedge S(\psi
_{i})$, and therefore by overspill, there is some nonstandard $d\in \omega ^{%
\mathcal{M}}$ such that $(\mathcal{M},S)\models \theta (d).$ It is now easy
to see (using the fact that $S$ satisfies Tarski's compositional clauses for
all $\Sigma _{c}$-formuale) that $p(x)$ is realized in $\mathcal{M}$.\hfill $%
\square $ \medskip

\noindent \textbf{3.7.~Corollary.~}\textit{A countable} $\omega $-\textit{nonstandard model} $\mathcal{M}%
\models \mathrm{ZFC}$ \textit{has an expansion to} \textit{a model of} $%
\mathrm{ZFI}_{\mathrm{<}}$\ \textit{iff} $\mathcal{M}$ \textit{is
recursively saturated and} $\mathcal{M}\models \Lambda$.

\medskip

\noindent \textbf{Proof.} The left-to-right direction follows from Corollary 3.6 and Theorem 4.1. The right-to-left direction follows from Theorem 4.1 and the resplendence property of countable recursively saturated models \cite[Theorem 15.7]{Richard-book}.\footnote{Recall that if $\mathcal{M}$ is resplendent, and $\mathcal{M}$ has an elementary extension to a recursive (computable) theory $T$ (such as $\mathrm{ZFI}_{\mathrm{<}}$) formulated in a language extending the language of $\mathcal{M}$, then $\mathcal{M}$ has an expansion to a model of $T$.}\hfill $%
\square $ \medskip

\begin{itemize}
\item In what follows $\mathcal{M}_{X}$ is the elementary submodel of $\left( \mathcal{M},<_{M}\right)$ generated by $X$, as
in part (j) of Definition 2.1.1, thus the universe $M_X$ of $\mathcal{M}_{X}$ consists of the elements of $M$ that are
definable in $\left( \mathcal{M},<_{M}\right) $ with parameters from $X$.
\end{itemize}

\noindent \textbf{3.8.~Theorem.}~\textit{Suppose} $(\mathcal{M},<_{M},I)$
\textit{is an} $\omega $-\textit{standard model of} $\mathrm{ZFI}_{\mathrm{<}%
}$\textit{. Then}:\medskip

\noindent \textbf{(a)} \textit{For each subset} $X$ \textit{of} $M$ \textit{%
that is definable in} $(\mathcal{M},<_{M},I)$, $M_{X}$ \textit{is definable
in} $\left( \mathcal{M},<_{M},I \right) .$\medskip

\noindent \textbf{(b) }$\mathcal{M}_{I_{1}}\cong \mathcal{M}_{I_{2}}$
\textit{for any cofinal subsets} $I_{1}$ \textit{and} $I_{2}$ \textit{of} $I$
\textit{that are definable in} $(\mathcal{M},<_{M},I).$ \textit{Moreover,
the isomorphism between} $\mathcal{M}_{I_{1}}$ \textit{and} $\mathcal{M}%
_{I_{2}}$ \textit{is definable in} $(\mathcal{M},<_{M},I).$\medskip

\noindent \textbf{(c)} \textit{There is a nontrivial elementary embedding} $%
j:\mathcal{M}_{I}\rightarrow \mathcal{M}_{I}$ \textit{such that} $j$ \textit{%
is definable in} $(\mathcal{M},<_{M},I).$\medskip

\noindent \textbf{(d)} $M_{I}$ \textit{is a proper subset of }$M$.\medskip

\noindent \textbf{(e)} $\mathcal{M}\models \ ``0^{\#}\ \mathrm{exists}$%
\textquotedblright ,\textit{\ in particular} $\mathcal{M}\models \mathrm{V}%
\neq \mathrm{L}$.\medskip

\noindent \textbf{(f)}\textit{\ The core model} $\mathrm{K}^{\mathcal{M}}$
\textit{of} $\mathcal{M}$ \textit{satisfies \textquotedblleft there is a
proper class of almost Ramsey cardinals}\textquotedblright\ (\textit{in the
sense of} \cite{Vickers-Welch}).\medskip

\noindent \textbf{Proof.}~(a) can be easily verified with the help of
Theorem 3.3.\medskip

To prove (b), first we observe that within $\mathrm{ZF}(\mathcal{L})$ (for
any $\mathcal{L}$) one can prove that if $I_{1}$ and $I_{2}$ are definable
cofinal subsets of the class of ordinals, then there is a definable
isomorphism $g:I_{1}\rightarrow I_{2}$. By Theorem 2.2.2, $g$ lifts to an
isomorphism $\widehat{g}:$ $\mathcal{M}_{I_{1}}\rightarrow \mathcal{M}%
_{I_{2}}.$ Let $\widehat{g}:$ $\mathcal{M}_{I_{1}}\rightarrow \mathcal{M}%
_{I_{2}}$ be given by

\begin{center}
$\widehat{g}(f(i_{1},\cdot \cdot \cdot ,i_{n}))=f(g(i_{1}),\cdot \cdot \cdot
,g(i_{n})),$
\end{center}

\noindent where $f$ is an $\mathcal{M}$-definable function. On the other
hand, by Theorem 3.3(d), there is a full satisfaction predicate over $(%
\mathcal{M},<_{M})$ that is definable in\textit{\ }$(\mathcal{M},<_{M},I),$
which together with (a) makes it clear that the proof of the fact that $%
\widehat{g}$ is an\ isomorphism of $\mathcal{M}_{I_{1}}$ and $\mathcal{M}%
_{I_{2}}$ can be carried out within \textit{\ }$(\mathcal{M},<_{M},I)$%
.\medskip

To see that (c)\textbf{\ }holds, we first observe that, reasoning in \textrm{%
ZFI}, there is definable order-isomorphism $f:\mathrm{Ord}\rightarrow I$,
and thus there is a definable enumeration $\left\langle i_{\xi }:\xi \in
\mathrm{Ord}\right\rangle $ of $I$, where $f(\xi )=i_{_{\xi }}$. Therefore
the map $h:I\rightarrow I$ given by $h(i_{\xi })=i_{\xi +1}$ is an $(%
\mathcal{M},I)$-definable order-preserving map whose range $I_{0}$ is a
proper subset of $I$. By part (a) of Theorem 2.2.2, $h$ induces an
elementary embedding $\widehat{h}$ of $\mathcal{M}_{I}$ onto $\mathcal{M}%
_{I_{0}},$ where $\mathcal{M}_{I_{0}}$ is a \textit{proper} elementary
submodel of $\mathcal{M}_{I}.$ Note that by Theorem 3.3, $\widehat{h}$ is
definable in $(\mathcal{M},<_{M},I)$. Thus $\widehat{h}$ is the desired
nontrivial elementary self-embedding $j$ of $\mathcal{M}_{I}$. \medskip

To verify (d)\textbf{, }suppose $M_{I}=M.$ Then by (c) there is a nontrivial
elementary embedding $j:\mathcal{M}\rightarrow \mathcal{M}$ such that $j$ is
$\mathcal{M}$-amenable. But Kunen's venerable theorem \cite[Theorem 17.7]%
{Jechbook-2003} bars the existence of such an embedding j. Thus $%
M_{I}\subsetneq M$. \medskip

Next we establish (e)\textbf{.} The fact that there is an $\mathcal{M}$%
-amenable satisfaction class over $\mathcal{M}$ makes it clear that there is
a cofinal subset $X\subseteq \mathrm{Ord}^{\mathcal{M}}$ such that $\mathcal{%
M}(\alpha )\prec \mathcal{M}$ for each $\alpha \in X$. Therefore for each $%
\alpha \in X$ the statement:

\begin{center}
\textquotedblleft $I\cap \mathrm{L}(\alpha )$ is a set of indiscernibles
over $\left( \mathrm{L}(\alpha ),\in \right) $\textquotedblright
\end{center}

\noindent holds in $\mathcal{M}$. So by picking an element $\alpha \in X$
such that $\mathcal{M}$ satisfies \textquotedblleft $I\cap \mathrm{L}(\alpha
)$ is uncountable\textquotedblright , we can deduce that $\mathcal{M}$
satisfies that $0^{\#}$ exists by a classical theorem of Silver \cite[%
Corollary 18.18]{Jechbook-2003}. Alternatively, one can put (c) together with
Kunen's theorem \cite[Theorem 18.20]{Jechbook-2003} that says that $0^{\#}$
exists iff the constructible universe admits a nontrivial elementary
self-embedding. This is because within $\left(\mathcal{M},<_M,I \right)$, there is an
isomorphism between \textrm{L}$^{\mathcal{M}}$ and \textrm{L}$^{\mathcal{M}%
_{I}}$, and therefore if $j:\mathcal{M}_{I}\rightarrow \mathcal{M}_{I}$ is a
nontrivial elementary embedding such that $j$ is definable in $\left(\mathcal{M}%
,<_{M},I \right)$, then $j$ induces a nontrivial elementary self-embedding of
\textrm{L}$^{\mathcal{M}}$ that is $\mathcal{M}$-amenable\textrm{.}\medskip

The proof of (f) is based on a key result of Vickers and Welch \cite%
{Vickers-Welch}, which states that if there is an inner model $\mathcal{M}_{0}$
of a model $\mathcal{M}$ of ZFC, and an $\mathcal{M}$-amenable nontrivial
elementary embedding $j:\mathcal{M}_{0}\rightarrow \mathcal{M}$,\ then the
core model $\mathrm{K}^{\mathcal{M}}$ of $\mathcal{M}$ satisfies
\textquotedblleft there is a proper class of almost Ramsey
cardinals\textquotedblright . Note that by (d), $\mathcal{M}_{I}$ is a
proper elementary submodel of $\mathcal{M}$, and by (a) its universe $M_{I}$
is definable in $(\mathcal{M},<_{M},I),$ therefore if $c:M_{I}\rightarrow
M_{0}$ is the collapsing map of $\mathcal{M}_{I}$ onto an inner model $%
\mathcal{M}_{0}$ of $\mathcal{M}$, then $c^{-1}:\mathcal{M}_{0}\rightarrow
\mathcal{M}$ is a nontrivial elementary embedding that is clearly $\mathcal{M%
}$-amenable.\hfill $\square $ \medskip

\noindent \textbf{3.9.~Remark.}~Parts (d), (e), and (f) of Theorem 3.8 can
be strengthened, as explained below. \medskip

Part (d) holds also when $\mathcal{M}$ is $\omega $-nonstandard. To see
this, suppose $M_{I}=M$ for $(\mathcal{M},<_{M},I)\models \mathrm{ZFI}_{%
\mathrm{<}}$, where $\mathcal{M}$ is $\omega $-nonstandard. By Theorem 3.3
there is a $\Sigma _{\omega }$-satisfaction class $S$ on $\mathcal{M}$ that
is definable in $(\mathcal{M},I)$. Consider the function

\begin{center}
$h:M\rightarrow \mathrm{Ext}_{\mathcal{M}}(\omega ^{\mathcal{M}}),$
\end{center}

\noindent where $h$ is defined in $(\mathcal{M},I)$ by $h(m):=$ the (G\"{o}%
del number of) the least $\mathcal{L}_{\mathrm{Set}}$-formula $\varphi (x,%
\overline{y})$ such that, as deemed by $S$, $m$ is defined by $\varphi (x,%
\overline{i})$ for some tuple $\overline{i}$ of parameters from $I,$ i.e., $%
S $ contains the sentences $\varphi (m,\overline{i})$ and $\exists !x\varphi
(x,\overline{i}).$ Note that the set of \textit{standard} elements of $\omega ^{%
\mathcal{M}}$ are definable in $(\mathcal{M},I)$ as the set of $i\in \omega $
such that $i<j\in \omega $ for some $j$ in the range of $h$. Thus $(\mathcal{%
M},I)$ is an $\omega $-nonstandard model of $\mathrm{ZFI}_{\mathrm{<}}$, in which the
set of standard elements of $\omega ^{\mathcal{M}}$ is definable, which is
impossible. \medskip

A straightforward modification of the proof of part (e) shows that the
statement \textquotedblleft $r^{\#}$ exists for all $r\subseteq \omega $%
\textquotedblright\ holds in every $\omega $-standard model of $\mathrm{ZFI}%
_{\mathrm{<}}$. \medskip

Finally, by taking advantage of the diagonal indiscernibility property of $I$,
the proof of part (f) can be modified to show that if $(\mathcal{M},<_{M},I)$
is an $\omega $-standard model of $\mathrm{ZFI}_{\mathrm{<}}$, then for any $%
m\in M$ there is an inner model $\mathcal{M}_{0}$ of $\mathcal{M}$ such that
$m\in M_{0}$, and for some $\mathcal{M}$-amenable nontrivial elementary
embedding $j:\mathcal{M}_{0}\rightarrow \mathcal{M}$, $j(m)=m$. This shows
that if $(\mathcal{M},<_{M},I)$ is an $\omega $-standard model of $\mathrm{%
ZFI}_{\mathrm{<}}$, then $\mathcal{M}$ exhibits \textquotedblleft inner
model reflection\textquotedblright\ in the sense of \cite{Barton et al},
i.e., any first order property of $\mathcal{M}$ (parameters allowed)
reflects to a proper inner model of $\mathcal{M}$. This result is a variant
of a theorem of Vickers and Welch \cite[Theorem 2.3$(i)$]{Vickers-Welch}
that derives inner model reflection from the existence of a proper class $I$
of \textquotedblleft good indiscernibles\textquotedblright\ for $(\mathrm{V}%
,\in )$.\footnote{$I$ is a set of good indiscernibles over a model $\left(
\mathcal{M},<_{M}\right) \models \mathrm{ZF}(<)+\mathrm{GW}$ if (1) $I$ is a
cofinal subset of $\mathrm{Ord}^{\mathcal{M}}$ that is $\mathcal{M}$%
-amenable, (2) for each $\alpha \in I,$ $(\mathcal{M}(\alpha ),<_{M_{\alpha}})\prec (\mathcal{M%
},<_{M})$, and (3) $I$ satisfies the diagonal indiscernibility scheme. In light of
Theorem 3.2, $I$ is a set of good indiscernibles over $\left( \mathcal{M}%
,<_{M}\right) \models \mathrm{ZF}(<)+\mathrm{GW}$ iff $\left( \mathcal{M}%
,<_{M},I\right) \models \mathrm{ZFI}_{<}$ and (2) holds. Thus $I$ is a set of good indiscernibles over $\left( \mathcal{M}%
,<_{M}\right) \models \mathrm{ZF}(<)+\mathrm{GW}$ iff $\left( \mathcal{M}%
,<_{M},I\right) \models \mathrm{ZFI}^{\mathrm{Good}}_{<}$, where $\mathrm{ZFI}^{\mathrm{Good}}_{<}$ is as in part $(iii)$ of Theorem 4.1. } \medskip

\noindent \textbf{3.10.~Corollary.}~\textit{No well-founded model} $\mathcal{%
M}$ of \textrm{ZF} \textit{that satisfies any of the conditions below has an
expansion to a model of }$\mathrm{ZFI}_{<}$\textrm{.}\medskip

\noindent \textbf{(a)} $\mathcal{M}\models $ $\mathrm{V=L.}$\medskip

\noindent \textbf{(b)} $\mathcal{M}=(\mathrm{V}(\kappa ),\in )$, \textit{%
where} $\kappa $ \textit{is the first cardinal satisfying} $P(\kappa)$, \textit{%
and} $P(\kappa )$ \textit{is a large cardinal property consistent with} $%
\mathrm{V=L}$, \textit{e.g.,} $P(\kappa )=$ \textquotedblleft $\kappa $%
\textit{\ is inaccessible/Mahlo/weakly compact/ineffable}\textquotedblright .%
$\bigskip $

\begin{center}
\textbf{4.~WHAT }$\mathbf{ZFI}_{<}$\textbf{\ KNOWS ABOUT SET THEORY}\bigskip
\end{center}

\noindent In contrast to the previous section whose main focus was on the
model-theoretic behavior of the theories $\mathrm{ZFI}$ and $\mathrm{ZFI}%
_{<} $, the main focus of this section is to use model-theoretic methods to
gauge the \textit{proof-theoretic strength} of these theories. As mentioned
in the introduction, a simple compactness argument shows that $\mathrm{ZFI}_{%
\mathrm{<}}$ is consistent if there is a weakly compact cardinal. The main
result of this section is Theorem 4.1, which pinpoints the set-theoretical
strength of $\mathrm{ZFI}_{<}$. Note that Theorem 4.1 shows that the
consistency strength of $\mathrm{ZFI}_{\mathrm{<}}$ is roughly the
consistency strength of the existence of an $\omega $-Mahlo cardinal, which
is considerably below the consistency strength of the existence of a weakly
compact cardinal. This calibration of the consistency strength of $\mathrm{%
ZFI}_{\mathrm{<}}$ also follows from part (b) of Theorem 5.5. \emph{As explained in Remark 4.8, Theorem 4.1 can be strengthened by adding two
additional equivalent conditions to the five equivalent conditions of the
theorem.}\medskip

\noindent \textbf{4.1.~Theorem}.~\textit{The following are equivalent for an
}$\mathcal{L}_{\mathrm{Set}}$-\textit{sentence} $\varphi $:\textit{\ }%
\medskip

\noindent $(i)$ $\mathrm{ZFI}_{<}^{\ast }\vdash \varphi$, where $\mathrm{ZFI}_{<}^{\ast }$ is the subsystem of $\mathrm{ZFI}_{<}$ axiomatized by: $$\mathrm{ZF(<)+ GW}+\mathrm{Coll}(<,I)+\mathrm{Cof}(I)+\mathrm{Indis}_{<}(I).$$

\medskip

\noindent $(ii)$ $\mathrm{ZFI}_{<}\vdash \varphi .$\medskip

\noindent $(iii)$ $\mathrm{ZFI}^{\mathrm{Good}}_{<} \vdash \varphi$, where $\mathrm{ZFI}^{\mathrm{Good}}_{<}:= \mathrm{ZFI}_{<} + \psi$, where $\psi$ is the single sentence expressing: $$\forall \alpha, \beta\in I [\alpha \in \beta \rightarrow (\mathrm{V}(\alpha),\in,<) \prec (\mathrm{V}(\beta),\in,<)].$$

\medskip

\noindent $(iv)$ $\mathrm{ZFC}+\Lambda \vdash \varphi .$\medskip

\noindent $(v)$ $\mathrm{GBC\ }+\mathrm{\ }$\textquotedblleft $\mathrm{Ord\
is\ weakly\ compact}$\textquotedblright\ $\vdash \varphi .$\medskip

\noindent \textbf{Proof.}~Recall that Theorem 2.4.12 assures us of the
equivalence of $(iv)$ and $(v)$. Since $(i)\Rightarrow (ii)$ and $(ii)\Rightarrow (iii)$ are both trivial,
the proof of the theorem will be complete once we establish:
 $$(iii)\Rightarrow (iv)\Rightarrow (i).$$

\medskip

\noindent To prove\textbf{\ }$(iii)\Rightarrow (iv),$ suppose
that for some $\mathcal{L}_{\mathrm{Set}}$-sentence $\varphi $ we
have:\medskip

\noindent (1) $\mathrm{ZFI}^{\mathrm{Good}}_{<}\vdash \varphi.$

\medskip

\noindent Assume on the contrary that $\mathrm{ZFC}+\Lambda +\lnot \varphi $
is consistent. By Theorem 2.4.12 and the completeness theorem for first
order logic, there is a model $\left( \mathcal{M}_{0},\mathfrak{X}%
_{0}\right) \models \mathrm{GBC\ }+\mathrm{\ }$\textquotedblleft $\mathrm{%
Ord\ is\ weakly\ compact}$\textquotedblright\ such that:\medskip

\noindent (2) $\mathcal{M}_{0}\models \lnot \varphi $. \medskip

\noindent Since by Theorem 2.4.9 for each metatheoretic natural number $%
n\geq 2$,

\begin{center}
$(\mathcal{M}_{0},\mathfrak{X}_{0})\models \mathrm{Ord}\rightarrow \left(
\mathrm{Ord}\right) _{2^{n}}^{n}$,
\end{center}

\noindent there is an elementary extension $(\mathcal{M},\mathfrak{X})$ of $(%
\mathcal{M}_{0},\mathfrak{X}_{0})$ such that for some nonstandard $c\in
\omega ^{\mathcal{M}}$ we have:$\medskip $

\noindent (3) $(\mathcal{M},\mathfrak{X})\models \mathrm{Ord}\rightarrow
\left( \mathrm{Ord}\right) _{2^{c}}^{c}.$ $\medskip $

\noindent Let $<_{M}$ be a member of $\mathfrak{X}_{0}$ such that $\left(
\mathcal{M},<_{M}\right) \models \mathrm{GW}.$ By Theorem 2.4.7(c) we can
get hold of a full satisfaction class $S\in \mathfrak{X}$ for $\left(
\mathcal{M},<_{M}\right) $. Since $S$ is $\mathcal{M}$-amenable, by using the reflection theorem within $(\mathcal{M},S)$, there is an $(\mathcal{M},S)$-definable unbounded subset $I_0$ of $\mathrm{Ord}^{\mathcal{M}}$ (in particular, $I_0 \in\mathfrak{X}$) such that:

\medskip

\noindent (4) For each $\alpha$ in $I_0$ $(\mathcal{M}(\alpha),<_{M_{\alpha}}, S_{\alpha}) \prec_{\Sigma_2(<,S)} (\mathcal{M},<_{M},S)$, where $S_{\alpha}:=S\cap M(\alpha)$.

\medskip

\noindent Since the predicate ``$S$ is a full satisfaction class for $(V,\in,<)$" is $\Pi_2(<,S)$, by (4) $S_{\alpha}$ is a full satisfaction class for $(\mathcal{M}(\alpha),<_{M_{\alpha}})$ for each $\alpha$ in $I_0$. This fact, in turn, readily implies:

\medskip

\noindent (5) $(\mathcal{M},I_0) \models \forall \alpha,\beta \in I_0[\alpha \in \beta \rightarrow (\mathrm{V}(\alpha), \in,<) \prec (\mathrm{V}(\beta), \in,<)].$
\medskip

\noindent On the other hand, it is easy to construct a recursive list of $%
\mathcal{L}_{\mathrm{Set}}(\mathrm{<)}$-formulae $\left\langle \varphi
_{i}(x_{1},\cdot \cdot \cdot ,x_{i}):i<\omega \right\rangle $ such that the
free variables of $\varphi _{i}$ are among $x_{1},\cdot \cdot \cdot ,x_{i}.$
This enumeration can be naturally prolonged within $\mathcal{M}$ so as to
obtain an enumeration $\left\langle \varphi _{i}(x_{1},\cdot \cdot \cdot
,x_{i}):i<c\right\rangle $ of $\mathcal{L}_{\mathrm{Set}}(\mathrm{<)}$%
-formulae in the sense of $\mathcal{M}.$ Next we define the following
evaluation function $e_S:[I_0]^{c}\rightarrow \{0,1\}^{c}$ within $(%
\mathcal{M},S)$ by:$\medskip $

\begin{center}
$e_{S}(\alpha _{1},\cdot \cdot \cdot ,\alpha _{c})=\left\langle \left\Vert
\varphi _{i}(\alpha _{1},\cdot \cdot \cdot ,\alpha _{i})\right\Vert
_{S}:i<c\right\rangle $,

\medskip

where $\alpha _{1}<\cdot \cdot \cdot <\alpha _{i}$ and for all $i<c$

\medskip

$\left\Vert \varphi _{i}(\alpha _{1},\cdot \cdot \cdot ,\alpha
_{i})\right\Vert _{S}=1$ iff $\varphi _{i}(\alpha _{1},\cdot \cdot \cdot
,\alpha _{i})\in S$.
\end{center}

\noindent By (3) there is some $I\in \mathfrak{X}$, such that $I \subseteq I_0$, $I$ is homogeneous for $%
e_S $ and $I$ is unbounded in $\mathrm{Ord}^{\mathcal{M}}$. It is evident that $I$ is
a cofinal set of indiscernibles over $\left( \mathcal{M},<_{M}\right) $ that
is $\mathcal{M}$-amenable. Thus, in light of (5) and the fact that $I_0 \subseteq I$, $(\mathcal{M},<_{M},I)\models
\mathrm{ZFI}^{\mathrm{Good}}_{<}$, so by (1) $\mathcal{M}\models \varphi $%
, which contradicts (2). This contradiction concludes the proof of $%
(iii)\Rightarrow (iv).$

\begin{itemize}
\item The proof of $(iv)\Rightarrow (i)$ of Theorem 4.1 relies on the
following lemma, in which $\mathcal{M}$ is a
model of $\mathrm{ZF}(<)+\mathrm{GW}$, $\mathcal{M}_{I}$ is the elementary submodel of $%
\mathcal{M}$ generated by $\mathcal{M}$-definable functions (as in part (j)
of Definition 2.1.1) and $\mathcal{M}_{I+J}$ is the elementary extension of%
\textit{\ }$\mathcal{M}_{I}$ resulting from stretching\textit{\ }$I$ to the
linear order\textit{\ }$I+J$, as in Theorem 2.2.2(b). Here $I+J$ is the
linear order on $I\cup J$ in which the elements of $J$ all exceed the
elements of $I$ (where $I$ and $J$ are disjoint). Thus, $\mathcal{M}$, $\mathcal{M}_{I}$, and $\mathcal{M}_{I+J}$ are $\mathcal{L}_{\mathrm{Set}}(<)$-structures that satisfy $\mathrm{ZF}(<)+\mathrm{GW}$.
\end{itemize}

\noindent \textbf{4.2.~Lemma}.~\textit{Suppose} $\mathcal{M}$ \textit{is a
model of} $\mathrm{ZF}(<)+\mathrm{GW}$ \textit{that has an expansion} $(\mathcal{M}%
,I)\models \mathrm{ZFI}_{<}^{\ast }$. \textit{Let} $\left(
J,<_{J}\right) $ \textit{be a linear order without a minimum element that is
disjoint from }$M$\textit{, and let} $\mathfrak{X:}=\mathrm{Cod}_{M_{I}}(%
\mathcal{M}_{I+J}).$ \textit{Then the following hold}:\medskip

\noindent \textbf{(a)} $\mathcal{M}_{I}\prec _{\mathrm{end}}\mathcal{M}%
_{I+J}.$\medskip

\noindent \textbf{(b) }$J$ \textit{is downward cofinal in}\textbf{\ }$%
M_{I+J}\backslash M_{I}$,\textit{\ i.e.}, $\forall x\in M_{I+J}\backslash
M_{I}$ $\exists j\in J$ $(j<_{M_{I+J}}x).$ \medskip

\noindent \textbf{(c) }$\left( \mathcal{M}_{I},\mathfrak{X}\right) \models
\mathrm{GBC.}$\medskip

\noindent \textbf{(d) }$\left( \mathcal{M}_{I},\mathfrak{X}\right) \models \
``\mathrm{Ord\ is\ weakly\ compact}$''.\medskip

\noindent \textbf{Proof.}~To prove (a), we note that by the Stretching
Theorem 2.2.2(b), $\mathcal{M}_{I+J}$ is an elementary extension of $%
\mathcal{M}_{I}$, so the proof of (a) is complete once we verify that $%
\mathcal{M}_{I+J}$ end extends $\mathcal{M}_{I}$. For this purpose, since $I$
is cofinal in the ordinals of $\mathcal{M}_{I}$ it suffices to show that if $%
f$ is an $\mathcal{M}$-definable function, where $f$ is $\left( n+s\right) $%
-ary, $\overline{i}\in \lbrack I]^{n}$, and $\overline{j}\in \lbrack J]^{s}$%
, then the following statement $(\nabla )$ holds for any $i\in I$:

\begin{center}
$(\nabla )$\qquad $\left[ \mathcal{M}_{I+J}\models f(\overline{i},\overline{j%
})\in \mathrm{V}(i)\right] \Longrightarrow f(\overline{i},\overline{j})\in
M_{I}$.
\end{center}

\noindent To establish $(\nabla )$, suppose: $\medskip $

\noindent (1) $\mathcal{M}_{I+J}\models f(\overline{i},\overline{j})\in
\mathrm{V}(i)$ for some $\overline{i}\in \lbrack I]^{n}$, $\overline{j}\in
\lbrack J]^{s}$, and $i\in I.$ \medskip

\noindent Let $i_{n}=\max (\overline{i}).$ Putting (1) together with the
assumption that $\mathcal{M}_{I+J}$ is obtained by stretching $I$ to\textit{%
\ }$I+J$ implies:\medskip

\noindent (2) $(\mathcal{M},I)\models \forall \overline{x}\in \lbrack
I]^{s}\ \left[ i_{n}<x_{1}<\cdot \cdot \cdot <x_{s}\rightarrow f(\overline{i}%
,x_{1},\cdot \cdot \cdot ,x_{s})\in \mathrm{V}(i)\right] .$\medskip

\noindent By $\mathcal{M}$-amenability of $I$, the Collection Scheme $%
\mathrm{Coll}(I)$ holds in $(\mathcal{M}%
,I)$, which coupled with (2) yields:\medskip

\noindent (3) $(\mathcal{M},I)\models \exists y\in \mathrm{V}(i)\ \forall
\alpha \in \mathrm{Ord}\ \exists \overline{x}\in \lbrack I]^{s}\ \left(
\alpha <x_{1}<\cdot \cdot \cdot <x_{s}\right) \wedge $ $f(\overline{i}%
,x_{1},\cdot \cdot \cdot ,x_{s})=y.$\medskip

\noindent Since $\mathcal{M}_{I}\prec \mathcal{M}$, by (3), we can find $%
\overline{k},\overline{l}\in \lbrack I]^{2}$ with $i_{n}<_{M}k_{1}<_{M}\cdot
\cdot \cdot <_{M}k_{s}<_{M}l_{1}<_{M}\cdot \cdot \cdot <_{M}l_{s}$ such
that:\medskip

\noindent (4) $\mathcal{M}_{I}\models f(\overline{i},\overline{k})=f(%
\overline{i},\overline{l}).$\medskip

\noindent By combining (4) with the assumption that $\mathcal{M}_{I+J}$ is
obtained by stretching $I$ to\textit{\ }$I+J$ we can conclude that $f(%
\overline{i},\overline{j})=f(\overline{i},\overline{k})\in M_{I}$, which
shows that $(\nabla )$ holds, thus completing the proof of (a). Note that
the assumption that $J$ has no minimum element was not invoked in the proof
of (a).\medskip

We next establish (b). In light of (a) it is sufficient to show:\medskip

\noindent $(\heartsuit )$ If $c\in M_{I+J}$ and $c<_{M_{I+J}}j $ for each $j\in J,$ then $c\in M_{I}.$ \medskip

\noindent We will establish the following stronger form $(\heartsuit ^{+})$
of $(\heartsuit )$. In what follows $f$ is an $(n+s)$-ary $\mathcal{M}$%
-definable function, $\overline{i}\in \lbrack I]^{n}$, $\overline{j}\in
\lbrack J]^{s}$ and $i_{n}=\max (\overline{i}).$\medskip

\noindent $(\heartsuit ^{+})$ If $f(\overline{i},\overline{j})\in M_{I+J}$,
and for all $j\in J$ $\mathcal{M}_{I+J}\models f(\overline{i},\overline{j}%
)<j,$ then $f(\overline{i},\overline{j})<k^{\prime }$ for any $k^{\prime
}\in I$ such that $i_{n}<_{M}k^{\prime }.$\medskip

\noindent To establish $(\heartsuit ^{+})$, suppose that for all $j\in J,$ $%
\mathcal{M}_{I+J}\models f(\overline{i},\overline{j})<j$. Since $J$ has no
minimum element, there is some $j^{\prime }\in J$ that is below $j_{1}=\min (%
\overline{j}),$ therefore: \medskip

\noindent (5) For all $j\in J,$ $\mathcal{M}_{I+J}\models f(\overline{i},%
\overline{j})<j^{\prime }<j_{1}$. \medskip

\noindent The fact that $\mathcal{M}_{I+J}$ is the elementary extension of $%
\mathcal{M}_{I}$ resulting from stretching $I$ to $I+J$ assures us that if
we choose $k^{\prime }>i_{n}$ with $k^{\prime }\in I$ and some $\overline{k}%
\in \lbrack I]^{s}$ with $k^{\prime }<k_{1}=\min (\overline{k}),$
then:\medskip

\noindent (6) $\mathcal{M}_{I}\models f(\overline{i},\overline{k})<k^{\prime
}$.\medskip

\noindent Thanks to (6), we can conclude that $\mathcal{M}_{I+J}\models f(%
\overline{i},\overline{j})<k^{\prime },$ thus $(\heartsuit ^{+})$ holds.
This concludes our verification of (b). Note that (b) implies that $%
M_{I+J}\backslash M_{I}$ has no $<_{M_{I+J}}$-minimum element.\medskip

To establish (c), we first claim:

\begin{center}
$(\diamondsuit )$\qquad If $X\in \mathfrak{X}$, then $X=D\cap M_{I}$ for
some $D\subseteq M$ such that $D$ is parametrically $(\mathcal{M},I)$%
-definable.
\end{center}

\noindent To demonstrate $(\diamondsuit )$ let $X\in \mathfrak{X,}$ and
choose $a\in M_{I+J}$ such that $X=\mathrm{Ext}_{\mathcal{M}_{I+J}}(a).$
Then $a=f(\overline{i},\overline{j})$ for some $\mathcal{M}$-definable $%
(n+s) $-ary function $f$, where $\overline{i}\in \lbrack I]^{n}$, and $%
\overline{j}\in \lbrack J]^{s}.$ Thus

\begin{center}
$X=\{m\in M_{I}:\mathcal{M}_{I+J}\models m\in f(\overline{i},\overline{j}%
)\}. $
\end{center}

\noindent Note that the veracity of the diagonal indiscernibility scheme in $%
\mathcal{M}_{I+J}$ implies:\medskip

\noindent (7) For $m\in M,$ $\mathcal{M}_{I+J}\models m\in f(\overline{i},%
\overline{j})$ iff there is some \textquotedblleft sufficiently
large\textquotedblright\ $\overline{k}\in \lbrack I]^{s}$, $\mathcal{M}%
\models m\in f(\overline{i},\overline{k}),$ \medskip

\noindent where \textquotedblleft sufficiently large\textquotedblright\
means that there is some $u\in I$ such that $\mathcal{M}(u)$ (i.e., $\mathrm{%
V}^{\mathcal{M}}(u)$) contains $\overline{i}$ and $m$ and $u<k_{1}$ (recall
that $k_{1}=\min (\overline{k})).$ Thus (7) makes it clear that $X=D\cap
M_{I}$, where:\medskip

\begin{center}
$D:=\left\{ m\in M:(M,I)\models \exists u\in I\ \exists \overline{k}\in
\lbrack I]^{s}\ [u<k_{1}\wedge \{\overline{i},m\}\subseteq \mathrm{V}^{%
\mathcal{M}}(u)\wedge m\in f(\overline{i},\overline{k})]\right\} $.
\end{center}

\noindent This concludes the verification of $(\diamondsuit )$. Note that $%
(\diamondsuit )$ readily implies that for each $X \in \mathfrak{X}$ $(\mathcal{M}_{I},X)$ satisfies the following weak form of $\mathrm{Coll}(X)$:

\medskip
\noindent (8)  $(\mathcal{M}_{I},X)\models \forall v \left[  \left(\forall x\in v\text{ }\exists y\text{\ } (\langle x,y \rangle \in X
)\right) \rightarrow \left( \exists w\text{ }\forall x\in v\text{ }\exists
y\in w\text{ }(\langle x,y \rangle\in X)\right) \right].$

\medskip

\noindent We will next verify $(\triangle )$ below, which together with (8) will allow us to conclude that $(\mathcal{M}_{I},X)$ satisfies the full scheme $\mathrm{Coll}(X)$.

\begin{center}
$(\triangle )$\qquad\ If $X\in \mathfrak{X}$ and $Y=\{m\in M_{I}:(\mathcal{M}_{I},p,X)\models\varphi
(m,p,X)\}$ for some $\mathcal{L}_{\mathrm{Set}}(<,X)$-formula $\varphi (x,y,X)$
and some parameter $p\in M,$ then $Y\in \mathfrak{X.}$
\end{center}

\noindent The proof of $(\triangle )$ is carried out
by induction on the complexity of $\varphi .$ We may assume that the logical
connectives consist of $\{\lnot ,\vee ,\exists \}.$ The atomic case and the
Boolean cases go through smoothly (since $\mathfrak{X}$ is readily seen to
be closed under complements and unions), but the existential case requires a
nontrivial argument. To handle the existential case, we need to show:

\begin{center}
$(\ast )$\qquad If $Y:=\{x\in M_{I}:\exists y\in M_{I}\ \langle x,y \rangle\in X\}$, then $%
Y\in \mathfrak{X.}$
\end{center}

\noindent Let $r\in M_{I+J}$ such that $X=\mathrm{Ext}_{\mathcal{M}%
_{I+J}}(r)\cap M_{I}.$ To verify $(\ast )$, it is sufficient to show $(\ast
\ast )$\ below:

\begin{center}
$(\ast \ast )$\qquad There is some $c\in M_{I+J}$ such that $\forall x,y\in
M_{I+J}\ (\mathcal{M}_{I+J}\models \lbrack \langle x,y \rangle\in r$ and $y<c]$ iff $y\in
M_{I}).$
\end{center}

\noindent To see that $(\ast \ast )\Rightarrow (\ast )$, choose $d$ in $%
\mathcal{M}_{I+J}$ such that $\mathcal{M}_{I+J}\models d=\left\{ x:\exists
y<c\ \langle x,y \rangle\in r\right\} $, thus:

\begin{center}
$Y=\mathrm{Ext}_{\mathcal{M}_{I+J}}(d)\cap M_{I}$,
\end{center}

\noindent which makes it clear that $Y\in \mathfrak{X}.$ In order to
establish $(\ast \ast )$, choose a function $g$ in $\mathcal{M}_{I+J}$ such
that $\mathcal{M}_{I+J}$ thinks that the domain of $g$ is the same as the
domain $\mathrm{Dom}(r)$ of $r$, where $\mathrm{Dom}(r)=\left\{ x:\exists
y\langle x,y \rangle \in r\right\} $ and

\begin{center}
$\mathcal{M}_{I+J}\models \forall x\in \mathrm{Dom}(r)$ $[g(x)=$ $<$-least $%
y $ such that $\langle x,y \rangle\in r].$
\end{center}

\noindent Choose an $\mathcal{M}$-definable function $f$, where $\overline{i}%
\in \lbrack I]^{n}$, and $\overline{j}\in \lbrack J]^{s}$ such that $g=f(%
\overline{i},\overline{j})$, and let

\begin{center}
$G:=\mathrm{Ext}_{\mathcal{M}_{I+J}}(g)\cap M_{I}.$
\end{center}

\noindent Note that $Y=\{x\in M_{I}:g(x)\in M_{I}\}.$ We will establish $%
(\ast \ast )$ by showing that there is a lower bound $c\in M_{I+J}$ for $%
\left\{ g(x):x\in M_{I},\ g(x)\notin M_{I}\right\} $ (in the sense of $%
<_{M_{I+J}}).$ We may assume that $I$ is cofinal in $(M_{I},<_{M})$ (by
replacing $<_{M}$, if necessary, with $<^{*}_{M}$, as in part (k) of Definition 2.1.1). Coupled with (8), we may conclude:
\medskip

\noindent (9) $\forall k\in I\ \exists k^{\prime }\in I$ such that $(%
\mathcal{M}_{I},G)\models $ $\forall x<k\ (x\in \mathrm{Dom}(G)\rightarrow G(x)<k^{\prime }).$\medskip

\noindent Given elements $k$ and $k^{\prime }$ of $I$, let

\begin{center}
$Z_{k,k^{\prime }}:=\left\{ y\in M_{I+J}:\mathcal{M}_{I+J}\models \exists
x<k,\ y=g(x)>k^{\prime }\right\} .$
\end{center}

\noindent Choose any $k\in I.$ By (9) there is some $k^{\prime }\in I$ such
that $Z_{k,k^{\prime }}\subseteq M_{I+J}\backslash M_{I}$. Reasoning in $%
\mathcal{M}_{I+J}$, let $u\in M_{I+J}$ be the $<$-least element of $%
Z_{k,k^{\prime }}$. By part(b) of Lemma 4.2, there is some $%
j^{\prime }\in J$ such that $j^{\prime }<u$ and therefore $j^{\prime }$ is
a strict lower bound for $Z_{k,k^{\prime }}.$ Thus:

\medskip

\noindent (10) $\mathcal{M}_{I+J}\models \forall x<k$ $\forall y\ \left[
\left( y=g(x)>k^{\prime }\right) \rightarrow j^{\prime }<y\right] $.

\medskip

\noindent Recall that $g=f(\overline{i},\overline{j}).$ So (10) states that $\mathcal{M}_{I+J}$ satisfies a particular first order statement with parameters $k$, $k^{\prime}$, $j^{\prime}$, $\overline{i}$, and $\overline{j}$ (all of which are in $I$), which coupled with the indiscernibility property of $I$, shows that (11) holds for any $%
k<k^{\prime }$ in $I$, as long as $k$ is above $i_{n}$. This shows that $j^{\prime }$ serves as the element $c$
in $(\ast \ast )$, thus concluding the verification of ($\triangle )$. \medskip

\noindent Thanks to (8) and ($\triangle )$ we have:

\medskip

\noindent (11) $(\mathcal{M}_{I},X)\models \mathrm{Coll}(X)$.

\medskip

\noindent On the other hand, each member of $\mathfrak{X}$ is clearly
piecewise coded in $\mathcal{M}_{I}$ since $\mathcal{M}_{I}$ is a
rank-extension of $\mathcal{M}_{I+J}$ (thanks to (a) and the fact that
elementary end extensions of models of $\mathrm{ZF}$ are rank extensions).
As pointed out in 2.2.2(g), the piecewise codability of $X$ together with
(11) allows us to conclude that $(\mathcal{M}_{I},X)\models \mathrm{ZF}%
(X)$. Therefore thanks to the fact that finitely many members of $%
\mathfrak{X}$ can be coded by a single member of $\mathfrak{X,}$ we
have:\medskip

\noindent (12) $(\mathcal{M}_{I},X_{1},...,X_{n})\models \mathrm{ZF}%
(X_{1},...,X_{n})$ for any finite subset $\{X_{1},...,X_{n}\}$ of $\mathfrak{X.}$\medskip

\noindent In light of Remark 2.4.2, ($\triangle )$ and (12) make it clear that $\left( \mathcal{M}_{I},%
\mathfrak{X}\right) \models \mathrm{GBC}$, thus concluding the proof of (c).

\medskip

Finally, we turn to establishing (d). Suppose $\tau $ is an \textrm{Ord}%
-tree (in the sense of $\mathcal{M}_{I}$) coded in $\mathfrak{X}$. Thanks to the
existence of the global well-ordering $<_{M}$ in $\mathfrak{X}$ we may
assume without loss of generality that $\tau =(M_{I},\vartriangleleft )$ for
some tree-ordering relation $\vartriangleleft $ coded in $\mathfrak{X}$. Fix
$r\in M_{I+J}$ such that

\begin{center}
$\mathrm{Ext}_{\mathcal{M}_{I+J}}(r)\cap M=\left\{ (x,y)\in
M_{I}:x\vartriangleleft y\right\} .$
\end{center}

\noindent Within $\mathcal{M}_{I+J}$ let $k$ be the field of $r$, i.e., the
set of elements that occur as the first or second coordinates of an ordered
pair in $r$. Without loss of generality we may assume that every element of $%
r$ is an ordered pair from the point of view of $\mathcal{M}_{I+J}.$
Consider the relational structure $\tau ^{\ast }:=(k,r)\in M_{I+J}$. Within $%
\mathcal{M}_{I+J}$ for each ordinal $\alpha ,$ let $\tau ^{\ast }(\alpha )$
be the initial segment of $\tau $ consisting of elements of $\tau ^{\ast }$
whose rank (in the tree $\tau ^{\ast }$) is at most $\alpha ,$ and let $%
\theta (x)$ be the following formula that expresses:

\begin{center}
\textquotedblleft $x\in \mathrm{Ord}$ and $\tau ^{\ast }(x)$ is a
well-founded tree\textquotedblright $.$
\end{center}

\noindent The assumption that $(\mathcal{M}_{I},\tau )\models $
\textquotedblleft $\tau $ is a well-founded tree of height
Ord\textquotedblright\ implies:\medskip

\noindent (12) $\mathcal{M}_{I+J}\models \theta (\alpha )$ for all $\alpha
\in \mathrm{Ord}^{\mathcal{M}_{I}}$.\medskip

\noindent Recall that (c) implies that \textrm{Ord}$^{\mathcal{M}%
_{I+J}}\backslash M_{I}$ has no least element. Therefore (12) assures us via
a simple overspill argument there is some $\beta \in \mathrm{Ord}^{\mathcal{M%
}_{I+J}}$ such that $\mathcal{M}_{I+J}\models \theta (\beta ).$ This shows
that the initial segment $\tau (\beta )$ of $\tau ^{\ast }$ in $\mathcal{M}%
_{I+J}$ properly end extends $\tau $, i.e., $\tau (\beta )$ does not contain
any new elements $\vartriangleleft $-below the elements of $\tau $. So we
can construct a branch $B$ of $\tau $ such that $B\in \mathfrak{X}$ by
considering the elements below a member of $\tau (\beta )$ whose height is
above Ord$^{M}.$ More specifically, choose $t\in \mathrm{Ext}_{\mathcal{M}%
_{I+J}}(\tau (\beta ))\backslash M_{I}$, and define the desired branch $B\in
\mathfrak{X}$ of $\tau $ by

\begin{center}
$B:=\{m\in M_{I}:\left( m,t\right) \in r\}.$
\end{center}

\hfill $\square $ (Lemma 4.2)\medskip

With Lemma 4.2 at hand, we are now in a position to smoothly verify the
direction $(iv)\Rightarrow (i)$ of Theorem 4.1. Suppose $\mathrm{ZFC}%
+\Lambda \vdash \varphi ,$ and assume on the contrary that $\mathrm{ZFI}%
_{<}^{\ast }+\lnot \varphi $ is consistent, and therefore there is a
countable $(\mathcal{M},I)$ of $\mathrm{ZFI}_{<}^{\ast }$ such that $%
\mathcal{M}\models \lnot \varphi .$ Let $J$ be any linear order with no
minimum element that is disjoint from $M$, and let $\mathcal{M}_{I+J}$ be
the elementary extension of\textit{\ }$\mathcal{M}_{I}$ resulting from
stretching\textit{\ }$I$ to the linear order\textit{\ }$I+J.$ By Lemma 4.2 $%
\mathcal{M}_{I+J}$ is an elementary end extension of $\mathcal{M}_{I}$ and $%
\left( \mathcal{M}_{I},\mathfrak{X}\right) \models \mathrm{GBC}+$
\textquotedblleft $\mathrm{Ord\ is\ weakly\ compact}$\textquotedblright . So
by Theorem 2.4.12 $\mathcal{M}_{I}$ satisfies $\mathrm{ZFC}+\Lambda $, which
in light of the fact that $\mathcal{M}_{I}\prec \mathcal{M}$ implies that $%
\varphi $ holds in $\mathcal{M}$, contradiction.\hfill $\square $ (Theorem
4.1)\medskip

\noindent \textbf{4.3.~Remark.}~It is not clear whether the scheme $\mathrm{%
\Lambda }$ is provable in $\mathrm{ZFCI}$ (i.e., \textrm{ZFI} plus the axiom
of choice)$\mathrm{.}$ However, note that by part (b) of Theorem 3.2 for
any $(\mathcal{M},I)\models \mathrm{ZFI}$, and any $p\in M$, a tail of $I$\
is indiscernible in $\mathrm{HOD}^{\mathcal{M}}(p).$ Together with the fact that
there is a well-ordering of $\mathrm{HOD}^{\mathcal{M}}(p)$ that is
parametrically definable in $\mathcal{M}$, one can use the strategy of the $%
(iv)\Rightarrow (i)$ direction of the proof of Theorem 4.1 so as to show
that if $(\mathcal{M},I)\models \mathrm{ZFI}$, and $m\in M$, then $\mathrm{%
HOD}^{\mathcal{M}}(m)\models \mathrm{\Lambda .}$ \textbf{\medskip }

\noindent \textbf{4.4.~Remark.}~The proof of Theorem 4.1 makes it clear that
the following hold:\textbf{\medskip }

\noindent \textbf{(a)} If $(\mathcal{M},\mathfrak{X})\models \mathrm{GBC}+$
\textquotedblleft $\mathrm{Ord\ is\ weakly\ compact}$\textquotedblright $,$
and $\mathcal{M}$ is $\omega $-nonstandard, then $\mathcal{M}$ has an
expansion to a model of $\mathrm{ZFI}^{+}_{<}.$
\medskip

\noindent \textbf{(b)} If $(\mathcal{M},I)\models \mathrm{ZFI}_{\mathrm{<}}$, then
the elementary submodel $\mathcal{M}_{I}$ of $\mathcal{M}$ has an expansion
to a model of $\mathrm{GBC}+$ \textquotedblleft $\mathrm{Ord\ is\ weakly\
compact}$\textquotedblright .\textbf{\medskip }

Next we define the extensions $\mathrm{ZFI}_{\mathrm{<}}^{k}$ and $\mathrm{%
ZFI}_{\mathrm{<}}^{\omega }$ of $\mathrm{ZFI}_{<}$, which despite their
powerful appearance, turn out to be rather mild extensions of $\mathrm{%
ZFI}_{<}.$\textbf{\medskip }

\noindent \textbf{4.5.~Definition.}~The theory\textbf{\ }$\mathrm{ZFI}_{%
\mathrm{<}}^{\omega }$ is the union of the theories $\mathrm{ZFI}_{\mathrm{<}%
}^{k}$ for $1\leq k\in \omega \mathrm{,}$ where $\mathrm{ZFI}_{\mathrm{<}%
}^{k}$ is formulated in the language $\mathcal{L}_{k}=\mathcal{L}_{\mathrm{%
Set}}\cup \{I_{j}(x):j<k\}$, and each $I_{j}\mathrm{(}x\mathrm{)}$ is a
unary predicate. The axioms of $\mathrm{ZFI}_{\mathrm{<}}^{1}$ are obtained
from the axioms of $\mathrm{ZFI}_{\mathrm{<}}$ simply by renaming $I$ as $%
I_{0}.$ The axioms of $\mathrm{ZFI}_{\mathrm{<}}^{k+1}$ consist of the union
of the axioms of $\mathrm{ZFI}_{\mathrm{<}}^{k}$ with the following four
groups of sentences:\medskip

\noindent $(1_{k+1})$ $\mathrm{ZFC(}\mathcal{L}_{k+1});$\medskip

\noindent $(2_{k+1})$ The sentence $\mathrm{Cof}(I_{k})$ expressing
\textquotedblleft $I_{k}$\ is a cofinal subclass of the class of
ordinals\textquotedblright ; and\medskip

\noindent $(3_{k+1})$ The scheme $\mathrm{Indis}_{k}(I_{k})=\{\mathrm{Indis}%
_{\varphi }(I_{k}):$ $\varphi $ is a formula of $\mathcal{L}_{k}\}$ ensuring
that $I_{k}$ is a class of order indiscernibles for the structure $\mathrm{%
(V,\in ,<,}I_{j}\mathrm{)}_{j<k}.$ More explicitly, for each $n$-ary formula
$\varphi (v_{1},\cdot \cdot \cdot ,v_{n})$ in the language $\mathcal{L}_{k},$
$\mathrm{Indis}_{\varphi }(I_{k})$ is the following sentence:\medskip

\begin{center}
$\forall x_{1}\in I_{k}\cdot \cdot \cdot \forall x_{n}\in I_{k}$ $\forall
y_{1}\in I_{k}\cdot \cdot \cdot \forall y_{n}\in I_{k}$ \textbf{\medskip }

$[(x_{1}<\cdot \cdot \cdot <x_{n})\wedge (y_{1}<\cdot \cdot \cdot
<y_{n})\rightarrow (\varphi (x_{1},\cdot \cdot \cdot ,x_{n})\leftrightarrow
\varphi (y_{1},\cdot \cdot \cdot ,y_{n}))].$\medskip
\end{center}

\noindent $(4_{k+1})$ The sentence asserting that $I_{k}$ is subclass of $%
I_{k-1}$ (for $k\geq 1).$\medskip

\begin{itemize}
\item Thus $\mathrm{ZFI}_{\mathrm{<}}^{k+1}$ bears the same relation to $%
\mathrm{ZFI}_{\mathrm{<}}^{k}$ that $\mathrm{ZFI}_{\mathrm{<}}$ bears to $%
\mathrm{ZF+GW,}$ i.e., for $1\leq k\in \omega $, a model $\left( \mathcal{M}%
,<_{M}\right) $ of $\mathrm{ZF+GW}$ has an expansion to a model of $\mathrm{%
ZFI}_{\mathrm{<}}^{k+1}$ iff there is a nested sequence $I_{0}\supseteq
\cdot \cdot \cdot \supseteq I_{k}$ of cofinal subsets of $\mathrm{Ord}^{%
\mathcal{M}}$ such that $\left( I_{0},\in \right) $ is indiscernible over $%
\left( \mathcal{M},<_{M}\right) $, $\left( I_{1},\in \right) $ is
indiscernible over $\left( \mathcal{M},<_{M},I_{0}\right) $, $\cdot \cdot
\cdot $, and $\left( I_{k},\in \right) $ is indiscernible over $\left(
\mathcal{M},<_{M},I_{0},\cdot \cdot \cdot ,I_{k-1}\right) .$\medskip
\end{itemize}

\noindent \textbf{4.6.~Theorem.}~\textit{Suppose} $\varphi $\textit{\ is a
sentence in the language }$\mathcal{L}_{\mathrm{Set}}$, \textit{then}%
:\medskip

\begin{center}
$\mathrm{ZFI}_{\mathrm{<}}^{\omega }\vdash \varphi $ iff $\mathrm{GBC}+$
\textquotedblleft $\mathrm{Ord\ is\ weakly\ compact}$\textquotedblright\ $%
\vdash \varphi .$
\end{center}

\noindent \textbf{Proof.}~Note that the right-to-left direction of the above
equivalence is an immediate consequence of $(iii)\Rightarrow (i)$ of Theorem
4.1. The left-to-right direction of the above equivalence is an elaboration
of the proof of $(ii)\Rightarrow (iii)$ of Theorem 4.1. More explicitly, it
suffices to show that for any nonzero $k\in \omega $, if $(\mathcal{M},%
\mathfrak{X})\models $ $\mathrm{GBC}+$ \textquotedblleft $\mathrm{Ord\ is\
weakly\ compact}$\textquotedblright , then a sufficient condition for $%
\mathcal{M}$ to have an expansion to a model of $\mathrm{ZFI}_{\mathrm{<}%
}^{k}$ is that there is a nonstandard $c\in \omega ^{\mathcal{M}}$ such that:%
$\medskip $

\noindent (1) $(\mathcal{M},\mathfrak{X})\models \mathrm{Ord}\rightarrow
\left( \mathrm{Ord}\right) _{2^{c}}^{c}.$ $\medskip $

\noindent By the reasoning of the proof of $(iii)\Rightarrow (iv)$ of Theorem
4.1 using (1) we can find some $<_{M}$ in $\mathfrak{X}$ and $I_{0}\in
\mathfrak{X}$ such that $(\mathcal{M},<_{M},I_{0})\models \mathrm{ZFI}_{%
\mathrm{<}}^{1}.$ Let $\left\langle \varphi _{i}(x_{1},\cdot \cdot \cdot
,x_{i}):i<\omega \right\rangle $ be a recursive list of $\mathcal{L}_{%
\mathrm{Set}}(\mathrm{<},~I_{0}\mathrm{)}$-formulae $\left\langle \varphi
_{i}(x_{1},\cdot \cdot \cdot ,x_{i}):i<\omega \right\rangle $ such that the
free variables of $\varphi _{i}$ are among $x_{1},\cdot \cdot \cdot ,x_{i}$,
and let $\left\langle \varphi _{i}(x_{1},\cdot \cdot \cdot
,x_{i}):i<c\right\rangle $ be an extension of this enumeration in $\mathcal{M%
}$. Fix a full satisfaction class $S\in \mathfrak{X}$ for $(M,<_{M},I_{0})$
and let $e_{S}:[\mathrm{Ord}]^{c}\rightarrow \{0,1\}^{c}$ within $(\mathcal{M}%
,S) $ by:$\medskip $

\begin{center}
$e_{S}(\alpha _{1},\cdot \cdot \cdot ,\alpha _{c})=\left\langle \left\Vert
\varphi _{i}(\alpha _{1},\cdot \cdot \cdot ,\alpha _{i})\right\Vert
_{S}:i<c\right\rangle $, $\medskip $

where $\alpha _{1}<\cdot \cdot \cdot <\alpha _{i}$ and for all $i<c$

\medskip

$\left\Vert \varphi _{i}(\alpha _{1},\cdot \cdot \cdot ,\alpha
_{i})\right\Vert _{S}=1$ iff $\varphi _{i}(\alpha _{1},\cdot \cdot \cdot
,\alpha _{i})\in S$.
\end{center}

\noindent By (1) there is some $I_{1}\in \mathfrak{X}$ with $I_{1}\subseteq
I_{0}$ such that $I_{1}$ is homogeneous for $e$ and unbounded in $\mathrm{Ord%
}^{\mathcal{M}}$. It is evident that $\left( \mathcal{M},<_{M},I_{0},I_{1}%
\right) \models \mathrm{ZFI}_{\mathrm{<}}^{2}$. By repeating this argument
we can thus obtain an expansion of $\mathcal{M}$ that satisfies $\mathrm{ZFI}%
_{\mathrm{<}}^{k}$ for any desired nonzero $k\in \omega .$\hfill $\square $
\medskip

\noindent \textbf{4.7.~Theorem.}~If $\widehat{\mathcal{M}}=\left( \mathcal{M}%
,<_{M},I_{k}\right) _{1\leq k\in \omega }$ \textit{is a model of} $\mathrm{%
ZFI}_{\mathrm{<}}^{\omega }$ \textit{and} $\mathfrak{X}$ \textit{is the
collection of parametrically} $\widehat{\mathcal{M}}$-\textit{definable
subsets of} $M$, \textit{then} $(\mathcal{M},\mathfrak{X})\models \mathrm{GBC%
}+$ \textquotedblleft $\mathrm{Ord\ is\ weakly\ compact}$\textquotedblright .%
$\medskip $

\noindent \textbf{Proof.}~It should be clear that $(\mathcal{M},\mathfrak{X}%
)\models \mathrm{GBC}$. To verify that \textquotedblleft $\mathrm{Ord\ is\
weakly\ compact}$\textquotedblright\ holds in $(\mathcal{M},\mathfrak{X})$,
by Theorem 2.4.9 it suffices to verify that the partition relation $\mathrm{%
Ord}\rightarrow \left( \mathrm{Ord}\right) _{2}^{2}$ holds in $(\mathcal{M},%
\mathfrak{X}).$ But this is easy, since if for some $F\in \mathfrak{X}$ we
have:

\begin{center}
$(\mathcal{M},F)\models F:[\mathrm{Ord}]^{2}\rightarrow \{0,1\},$
\end{center}

\noindent then $F$ is definable in $\left( \mathcal{M},<_{M},I_{k}\right)
_{1\leq k\leq m}$ for some $m\in \omega $, and therefore $I_{m+1}$ is proper
class that is an unbounded $F$-homogeneous member of $\mathfrak{X}$, as
desired.\hfill $\square $ \medskip

\noindent \textbf{4.8.~Remark.}~If the model $\widehat{\mathcal{M}}$ as in
Theorem 4.7 is recursively saturated, then the proof of Theorem 4.1 of \cite%
{Ali BS for set theory} shows that the model $(\mathcal{M},\mathfrak{X})$ in
the statement of Theorem 4.7 also satisfies the scheme $\Sigma _{1}^{1}$-$%
\mathrm{AC}$ (and therefore the scheme $\Delta _{1}^{1}$-\textrm{CA}).
Together with Theorem 4.6, this shows that Theorem 4.1 can be strengthened
by asserting that the following two conditions $(vi)$ and $(vii)$ on an $%
\mathcal{L}_{\mathrm{Set}}$-sentence $\varphi $ are also equivalent to
conditions $(i)$ through $(v)$ of that theorem:\medskip

\noindent $(vi)$ $\mathrm{ZFI}_{\mathrm{<}}^{\omega }\vdash \varphi $\medskip

\noindent $(vii)$ $\mathrm{GBC}+$ \textquotedblleft $\mathrm{Ord\ is\ weakly\
compact}$\textquotedblright\ + $\Sigma _{1}^{1}$-$\mathrm{AC}\vdash \varphi $%
\medskip

The next result shows that a statement that one might expect to be provable
in the theory \textrm{GBC} $+$ \textrm{\textquotedblleft Ord is weakly
compact\textquotedblright } is actually unprovable in that theory. \medskip

\noindent \textbf{4.9.~Theorem.}~\textit{If the theory} \textrm{GBC} $+$
\textrm{\textquotedblleft Ord is weakly compact}\textit{\textquotedblright\
is consistent, then it does not prove the statement} $\theta =\forall m,n\in
\omega \left( \mathrm{Ord}\rightarrow \left( \mathrm{Ord}\right)
_{m}^{n}\right) .\medskip $

\noindent \textbf{Proof.} Let $\left\langle \varphi _{i}(x_{1},\cdot \cdot
\cdot ,x_{i}):i<\omega \right\rangle $ be a recursive list of $\mathcal{L}_{%
\mathrm{Set}}(\mathrm{<)}$-formulae such that the free variables of $\varphi
_{i}$ are among $x_{1},\cdot \cdot \cdot ,x_{i}.$ For each\textbf{\ }$n\in
\omega $ let $T_{n}$ be the fragment of $\mathrm{ZFI}_{\mathrm{<}}$ whose
axioms consist of $\mathrm{ZF}(\mathcal{L})+\mathrm{GW}$ for $\mathcal{L}=%
\mathcal{L}_{\mathrm{Set}}(I,<,)$ and $\mathrm{Cof}(I)$ and sentences
\textrm{Indisc}$_{\varphi _{i}}$ for $i\leq n.$ We next prove a key
lemma.\medskip

\noindent \textbf{4.10.~Lemma}.~$\mathrm{GBC}+\theta \vdash \forall n\in
\omega \ \mathrm{Con}(T_{n}).$\medskip

\noindent \textbf{Proof.}~Let $(\mathcal{M},\mathfrak{X})\models \mathrm{GBC}%
+$ $\theta .$ Then in particular $(\mathcal{M},\mathfrak{X})$ satisfies $%
\forall n\in \omega \left( \mathrm{Ord}\rightarrow \left( \mathrm{Ord}%
\right) _{2^{n}}^{n}\right) .$ Given any fixed $n\in \omega ^{\mathcal{M}}$
and arguing in $(\mathcal{M},\mathfrak{X}),$ we will show the consistency of
$T_{n}.$ By Theorem 2.4.7(c) there is a full satisfaction class $S\in
\mathfrak{X}$ for $\mathcal{M}$, which we can use together with $\mathrm{Ord}%
\rightarrow \left( \mathrm{Ord}\right) _{2^{n}}^{n}$ to get hold of an
unbounded homogeneous set $I\in \mathfrak{X}$ for the map $e_{S}^{n}:[%
\mathrm{Ord}]^{n}\rightarrow \{0,1\}^{n}$ that is defined within $(\mathcal{M%
},S)$ by:$\medskip $

\begin{center}
$e_{S}^{n}(\alpha _{1},\cdot \cdot \cdot ,\alpha _{n})=\left\langle
\left\Vert \varphi _{i}(\alpha _{1},\cdot \cdot \cdot ,\alpha
_{i})\right\Vert _{S}:i<n\right\rangle $, $\medskip $

where $\alpha _{1}<\cdot \cdot \cdot <\alpha _{i}$ and for all $i<n$

\medskip

$\left\Vert \varphi _{i}(\alpha _{1},\cdot \cdot \cdot ,\alpha
_{i})\right\Vert _{S}=1$ iff $\varphi _{i}(\alpha _{1},\cdot \cdot \cdot
,\alpha _{i})\in S$.$\medskip $
\end{center}

\noindent Clearly $(M,I)\models T_{n}.$ By Theorem 2.4.7(c) there is a full
satisfaction predicate $S$ for $(\mathcal{M},I)$\ such that $S\in \mathfrak{X%
}$, which shows that $\mathrm{Con}(T_{n})$ holds in $\mathcal{M}$. \hfill $%
\square $ (Lemma 4.10)\medskip

By Lemma 4.10 and compactness, $\mathrm{Con}(\mathrm{ZFI}_{<}\mathrm{)}$ is
provable in $\mathrm{GBC}+\theta .$ Since Theorem 4.1 is readily verifiable
in ZFC, the formal consistency of \textrm{GBC} $+$ \textit{\textquotedblleft
}\textrm{Ord is weakly compact}\textit{\textquotedblright }\textrm{\ }is
provable in $\mathrm{GBC}+\theta .$ In light of G\"{o}del's second
incompleteness theorem, the proof is complete. \hfill $\square $ (Theorem
4.9)

\bigskip

\begin{center}
\textbf{5. ~INTERPRETABILITY\ ANALYSIS OF }$\mathbf{ZFI}_{<}$\textbf{\
\bigskip }
\end{center}

In this section we study $\mathrm{ZFI}$ and $\mathrm{ZFI}_{<}$ through the
lens of interpretability theory, a lens that brings both the semantic and
syntactic features of the theories under its scope into a finer focus. We
review some relevant interpretability-theoretic preliminaries before presenting our
results.\medskip

\noindent \textbf{5.1.~Definitions.~}Suppose $U$ and $V$ are first order
theories, and for the sake of notational simplicity, let us assume that $U$
and $V$ are theories that \textit{support a definable pairing function}. We
use $\mathcal{L}_{U}$ and $\mathcal{L}_{V}$ to respectively designate the
languages of $U$ and $V$.\medskip

\noindent \textbf{(a) }An interpretation $\mathcal{I}$ of $U$ in $V$,
written:

\begin{center}
$\mathcal{I}:U\rightarrow V$,
\end{center}

\noindent is given by a translation $\tau $ of each $\mathcal{L}_{U}$%
-formula $\varphi $ into an $\mathcal{L}_{V}$-formula $\varphi ^{\tau }$
with the requirement that $V\vdash \varphi ^{\tau }$ for each $\varphi \in U$%
, where $\tau $ is determined by an $\mathcal{L}_{V}$-formula $\delta (x)$
(referred to as a \textit{domain formula}), and a mapping $P\mapsto _{\tau
}A_{P}$ that translates each $n$-ary $\mathcal{L}_{U}$-predicate $P$ into
some $n$-ary $\mathcal{L}_{V}$-formula $A_{P}$. The translation is then
lifted to the full first order language in the obvious way by making it
commute with propositional connectives, and subject to:

\begin{center}
$\left( \forall x\varphi \right) ^{\tau }=\forall x(\delta (x)\rightarrow
\varphi ^{\tau })$ and $\left( \exists x\varphi \right) ^{\tau }=\exists
x(\delta (x)\wedge \varphi ^{\tau }).$
\end{center}

\begin{itemize}
\item Note that each interpretation $\mathcal{I}:U\rightarrow V$\/gives rise
to an \textit{inner model construction} \textit{that \textbf{uniformly}
builds a model} $\mathcal{M}^{\mathcal{I}}\models U$\/ \textit{for any} $%
\mathcal{M}\models V$. \medskip
\end{itemize}

\noindent \textbf{(b) }$U$ is \textit{interpretable} in $V$ (equivalently: $%
V $\textit{\ interprets} $U$), written $U\trianglelefteq V$, iff there is an
interpretation $\mathcal{I}:U\rightarrow V.$ $U$ is \textit{locally} \textit{%
interpretable} in $V$, written $U\trianglelefteq _{\mathrm{loc}}V$ if $%
U_{0}\trianglelefteq V$ for every finitely axiomatizable subtheory $U_{0%
\text{ }}$of $U$.\medskip

\noindent \textbf{(c) }$U$\textbf{\ }and $V$ are \textit{mutually
interpretable} when $U\trianglelefteq V$ and $V\trianglelefteq U.$\medskip

\noindent \textbf{(d) }$U$ is a \textit{retract} of $V$ iff there are
interpretations $\mathcal{I}$ and $\mathcal{J}$ with $\mathcal{I}%
:U\rightarrow V$ and $\mathcal{J}:V\rightarrow U$, and a binary $U$-formula $%
F$ such that $F$\ is, $U$-verifiably, an isomorphism between \textrm{id}$%
_{U} $ (the identity interpretation on $U$) and $\mathcal{J}\circ \mathcal{I}
$\textsf{.} In model-theoretic terms, this translates to the requirement
that the following holds for every $\mathcal{M}\models U$:

\begin{center}
$F^{\mathcal{M}}:\mathcal{M}\overset{\cong }{\longrightarrow }\mathcal{M}%
^{\ast }:=\left( \mathcal{M}^{\mathcal{J}}\right) ^{\mathcal{I}}.$
\end{center}

\noindent \textbf{(e) }$U$ and $V$ are \textit{bi-interpretable} iff there
are interpretations $\mathcal{I}$ and $\mathcal{J}$ as above that witness
that $U$ is a retract of $V$, and additionally, there is a $V$-formula $G,$
such that $G$\ is, $V$-verifiably, an isomorphism between $\mathrm{id}_{V}$
and $\mathcal{I}\circ \mathcal{J}.$ In particular, if $U$ and $V$ are
bi-interpretable, then given $\mathcal{M}\models U$ and $\mathcal{N}\models
V $, we have

\begin{center}
$F^{\mathcal{M}}:\mathcal{M}\overset{\cong }{\longrightarrow }\mathcal{M}%
^{\ast }:=\left( \mathcal{M}^{\mathcal{J}}\right) ^{\mathcal{I}}$ and $G^{%
\mathcal{N}}:\mathcal{N}\overset{\cong }{\longrightarrow }\mathcal{N}^{\ast
}:=\left( \mathcal{N}^{\mathcal{I}}\right) ^{\mathcal{J}}.$
\end{center}

\noindent \textbf{(f) }The above notions can also be localized at a pair of
models. Suppose $\mathcal{N}$ is an $\mathcal{L}_{U}$-structure and $%
\mathcal{M}$ is an $\mathcal{L}_{V}$-structure. We say that $\mathcal{N}$ is
\textit{parametrically} \textit{interpretable} in $\mathcal{M}$, written $%
\mathcal{N}\trianglelefteq _{\mathrm{par}}\mathcal{M}$ (equivalently: $%
\mathcal{M}\trianglerighteq _{\mathrm{par}}\mathcal{N}$) iff the universe of
discourse of $\mathcal{N}$, as well as all the $\mathcal{N}$-interpretations
of $\mathcal{L}_{U}$-predicates are $\mathcal{M}$-definable. Similarly, we
say that $\mathcal{M}$ and $\mathcal{N}$ are \textit{parametrically
bi-interpretable }if there are parametric interpretations\textit{\ }$%
\mathcal{I}$\textit{\ }and\textit{\ }$\mathcal{J}$\textit{, }together with
an $\mathcal{M}$-definable $F$ and an $\mathcal{N}$-definable map $G$ such
that:

\begin{center}
$F^{\mathcal{M}}:\mathcal{M}\overset{\cong }{\longrightarrow }\mathcal{M}%
^{\ast }:=\left( \mathcal{M}^{\mathcal{J}}\right) ^{\mathcal{I}}$ and $G^{%
\mathcal{N}}:\mathcal{N}\overset{\cong }{\longrightarrow }\mathcal{N}^{\ast
}:=\left( \mathcal{N}^{\mathcal{I}}\right) ^{\mathcal{J}}.$\medskip
\end{center}

\noindent \textbf{(g) }A \textit{sequential theory} is a theory equipped
with a\textit{\ `}$\beta $-function' for handling finite sequences of
objects in the domain of discourse.$\medskip $

The following theorems are classical. Theorem 5.2 was first proved for
\textrm{PA} by Mostowski. His argument was later generalized by Montague as
in Theorem 5.2 below. In part (b) of the theorem, $\mathcal{L}_{\mathrm{Arith%
}}$ is the usual language of arithmetic $\{+,\cdot ,<,0,1\},$ and for $%
\mathcal{L\supseteq L}_{\mathrm{Arith}}$, \textrm{PA(}$\mathcal{L}$\textrm{)
}is the natural extension of $\mathrm{PA}$ in which $\mathcal{L}$-formulae
can appear in the scheme of induction.\medskip

\noindent \textbf{5.2.~Theorem.~}(Montague) \textit{If }$T$ \textit{is a
sequential theory and} $T$ \textit{can prove the induction scheme over its
ambient set of natural numbers, then} $T$ \textit{is a reflexive theory, i.e.%
}, $T$ \textit{proves the formal consistency of each of its finite} \textit{%
subtheories. In particular}:\medskip

\noindent \textbf{(a)} \textit{For all }$\mathcal{L}\supseteq \mathcal{L}_{%
\mathrm{Arith}}$, \textit{every extension} (\textit{in the same language})
\textit{of}\textbf{\ }\textrm{PA(}$\mathcal{L}$\textrm{)} \textit{is
reflexive.}\medskip

\noindent \textbf{(b)} \textit{For all }$\mathcal{L}\supseteq \mathcal{L}_{%
\mathrm{Set}}$, \textit{every extension} (\textit{in the same language})
\textit{of}\textbf{\ }\textrm{Z(}$\mathcal{L}$\textrm{)} \textit{is
reflexive, where }\textrm{Z(}$\mathcal{L}$\textrm{)}\textit{\ is Zermelo set
theory, as in Definition }$2.1.1(e)$\textit{.}\medskip

\noindent \textbf{5.3.~Theorem.~}(Orey's Compactness Theorem)\textbf{\ }%
\textit{If} $U$ \textit{is a reflexive sequential theory, and} $V\trianglelefteq _{\mathrm{loc}%
}U $\textit{\ for some recursively enumerable theory }$V$\textit{, then }$%
V\trianglelefteq U$.\medskip

\noindent \textbf{5.4.~Theorem.~}$\mathrm{GB}\ntrianglelefteq \mathrm{ZF}$%
\textrm{. \textit{Similarly,} GBC} \textit{is not interpretable in} $\mathrm{%
ZF}(<)+\mathrm{GW}$.\medskip

We are now ready to present the new results of this section. In part (b) of
Theorem 5.5, $\Lambda ^{-}$ is the subset of $\Lambda $ consisting of
sentences of the form \textquotedblleft there is an $n$-Mahlo
cardinal\textquotedblright\ for each metatheoretic $n\in \omega ,$ as in
part (c) of Definition 2.4.12.\medskip

\noindent \textbf{5.5.~Theorem.~}(Relative interpretability results)\medskip

\noindent \textbf{(a)} \textrm{GBC} $+$ \textrm{\textquotedblleft Ord is
weakly compact\textquotedblright\ }\textit{is not interpretable in} $\mathrm{%
ZFI}_{<}.$\medskip

\noindent \textbf{(b) }\textit{The theories} $\mathrm{ZFC}+\Lambda ^{-},$ $%
\mathrm{ZFC}+\Lambda $ \textit{and} $\mathrm{ZFI}_{\mathrm{<}}$ \textit{are
pairwise mutually interpretable.}\medskip

\noindent \textbf{(c)} $\mathrm{ZFI}_{<}$ \textit{is interpretable in }$%
\mathrm{GBC}$ $+$ \textrm{\textquotedblleft Ord is weakly
compact\textquotedblright }\textit{.}\medskip

\noindent \textbf{Proof}. ~The proof of (a) combines Theorem 4.1 together
with the strategy that proves Theorem 5.4. More specifically, since $\mathrm{%
ZFI}_{<}$ is a reflexive theory (by Theorem 5.2(a)), and the theory $\mathrm{%
GBC}$ $+$ \textrm{\textquotedblleft Ord is weakly compact\textquotedblright\
}is finitely axiomatizable, the interpretability of \textrm{GBC} $+$ \textrm{%
\textquotedblleft Ord is weakly compact\textquotedblright }\textit{\ }in $%
\mathrm{ZFI}_{<}$ would imply that $\mathrm{ZFI}_{<}$ proves the $\mathcal{L}%
_{\mathrm{Set}}$-sentence expressing the formal consistency of \textrm{GBC} $%
+$ \textrm{\textquotedblleft Ord is weakly compact\textquotedblright }$,$
which in light of Theorem 4.1 contradicts G\"{o}del's second incompleteness
theorem.\medskip

To prove (b), first recall that by part (c) of Theorem 2.4.11 $\mathrm{ZFC}%
+\Lambda ^{-}$ and $\mathrm{ZFC}+\Lambda $ are mutually interpretable\textit{%
. }Also note that since Theorem 4.1 assures us that $\mathrm{ZFI}_{\mathrm{<}%
}\vdash \Lambda $, the identity interpretation serves as a witness to the
interpretability of $\mathrm{ZFC}+\Lambda $ within $\mathrm{ZFI}_{\mathrm{<}%
}.$ So the proof of (b)\ will be complete once we establish the
interpretability of $\mathrm{ZFI}_{\mathrm{<}}$ within $\mathrm{ZFC}+\Lambda
.$ Towards this goal, thanks to Orey's Compactness Theorem 5.3, it will
suffice to show that every finite subtheory of $\mathrm{ZFI}_{\mathrm{<}}$
is interpretable in \textrm{ZFC + }$\Lambda .$ Indeed we will show that for
each $n\in \omega ,$ \textrm{ZFC + }$\Lambda $ can interpret the subtheory $%
T_{n}$ of $\mathrm{ZFI}_{<}$, where $T_{n}$ is the same theory as in the
proof of Theorem 4.9. Fix some nonzero $n\in \omega $, and reasoning within
\textrm{ZFC + }$\Lambda ,$ we consider the class function $F:\left[ \mathrm{%
Ord}\right] ^{n}\rightarrow \{0,1\}^{n}$, where

\begin{center}
$F(\alpha _{1},\cdot \cdot \cdot ,\alpha _{n})=\left\langle \left\Vert
\varphi _{i}(\alpha _{1},\cdot \cdot \cdot ,\alpha _{i})\right\Vert :1\leq
i\leq n\right\rangle ,$ where for all $i<n$

\medskip

$\left\Vert \varphi _{i}(\alpha _{1},\cdot \cdot \cdot ,\alpha
_{i})\right\Vert =1$ iff $\varphi _{i}(\alpha _{1},\cdot \cdot \cdot ,\alpha
_{i}).$
\end{center}

\noindent Note that if $F(\alpha _{1},\cdot \cdot \cdot ,\alpha
_{n})=F(\beta _{1},\cdot \cdot \cdot ,\beta _{n}),$ then $\varphi
_{i}(\alpha _{1},\cdot \cdot \cdot ,\alpha _{i})\leftrightarrow \varphi
_{i}(\beta _{1},\cdot \cdot \cdot ,\beta _{i})$ whenever $1\leq i\leq n.$
Let $\tau _{F}$ be the $\mathrm{Ord}$-tree as in Lemma 2.4.4 whose cofinal
branches are end-homogeneous proper classes for $F,$ i.e., not dependent on
the $n$-th coordinate of any increasing chain of length $n$. Let $\kappa
_{1} $ be an $m$-Mahlo cardinal such that $\mathrm{V}(\kappa _{1})$ is a $%
\Sigma _{m}$-elementary submodel of the universe, where $m\geq n$, each $\varphi_i$ is in $\Sigma_m$, and also $%
m$ is large enough so that the following statement is in $\Sigma _{m}$:

\begin{center}
\textquotedblleft $\tau _{F}$ is an Ord-tree, and the value of $F$ on any
increasing chain in $\tau _{F}$ of length $n$ is independent of its $n$-th
component\textquotedblright .
\end{center}

\noindent Choose any ordinal $\lambda $ above $\kappa _{1}$ and let $H_{1}$
be the intersection of $\kappa _{1}$ with the collection of ordinals that
are below $\lambda $ in the sense of the ordering $<_{F}$ of $\tau _{F}.$
Then by the choice of $\kappa _{1}$, $H_{1}$ is a cofinal branch of the tree
computed in $\left( \mathrm{V}(\kappa _{1}),\in \right) $ via the defining
formula of $\tau _{F}.$ Hence:

\begin{center}
$\left( \mathrm{V}(\kappa _{1}),\in ,H_{1}\right) \models H_{1}$ is end-homogeneous for $F$,
\end{center}

\noindent i.e., $F(\alpha _{1},\cdot \cdot \cdot ,\alpha _{n-1},\alpha
_{n})=F(\alpha _{1},\cdot \cdot \cdot ,\alpha _{n-1},\beta )$ for any
increasing elements $\alpha _{1}<\cdot \cdot \cdot <\alpha _{n-1}<\alpha
_{n} $ from $H_{1}$, and any $\beta \in H_{1}$ that is greater than $\alpha
_{n}.$

\medskip

In the next step we consider the restriction of the function $F$ to $[H_1]$ within $\left( \mathrm{V}(\kappa _{1}),\in ,H_{1}\right)$ and obtain an $\left( m-1\right) $-Mahlo cardinal $\kappa _{2}<\kappa _{1}$
that satisfies the following two properties:\medskip

\noindent (1) $\left(\mathrm{V}(\kappa _{2}),\in ,H_{1}\cap \mathrm{V}(\kappa _{2})\right) \prec  \left( \mathrm{V}(\kappa _{1}),\in ,H_{1}\right).$%

\medskip

\noindent (2) There is some cofinal subset $H_{2}$ of $\kappa _{2}$ such
that the value of $F$ on any increasing chain of length $n$ from $H_{2}$ is
independent of the choices of the $(n-1)$-th and the $n$-th components of
the chain. \medskip

Recall that $m\geq n,$ so by the $m$-Mahlo property of $\kappa _{1},$
we can repeat this process so as to obtain an $\left( m-n\right) $%
-Mahlo cardinal $\kappa _{n}$ and for some
cofinal subset $H_{n-1}$ of $\kappa _{n}$ such that $ \left( \mathrm{V}(\kappa _{n}),\in \right) \prec \left( \mathrm{V}(\kappa _{1}),\in
\right)$  such that for all increasing $n$-tuples $\overline{\alpha} = \left(\alpha_1,\cdot \cdot \cdot, \alpha_n\right)$ and $\overline{\beta}= \left(\beta_1,\cdot \cdot \cdot, \beta_n \right)$ from $H_{n-1}$, $\alpha_1 = \beta_1$  implies that $F(\overline{\alpha})=F(\overline{\beta})$. With such a set $H_{n-1}$ at hand, it is then easy to build $H_n \subseteq H_{n-1}$ (as in the proof of Theorem 2.4.5(b)) such that:
\begin{center}
$H_{n}$ is cofinal in $\kappa_n$ and $\left( \mathrm{V}(\kappa _{n}),\in ,H_{n}\right) \models H_{n}$ is homogeneous for $F$.%
\end{center}

\noindent This makes it clear that $\left( \mathrm{V}(\kappa _{n}),\in
,H_{n}\right) $ is our desired model of the subtheory $T_{n}$ of $\mathrm{ZFI%
}_{<}.$ This concludes the proof of (b).\medskip

Finally, to demonstrate (c), we can simply put part (b) together with
Theorem 2.4.12 that assures us that $\mathrm{ZFC}+\Lambda $ is provable in
the theory $\mathrm{GBC}$ $+$ \textrm{\textquotedblleft Ord is weakly
compact\textquotedblright .}\hfill $\square $\textbf{\medskip }

\noindent \textbf{5.6.~Remark.~}By a slight modification of the proof
strategy of part (b) of Theorem 5.5, one could also show that $\mathrm{ZFC}%
+\Lambda $ is mutually interpretable with the extension $\mathrm{ZFI}%
_{<}^{\omega }$ of $\mathrm{ZFI}_{<}$ studied in the previous section. This
modified proof can be combined with Theorem 4.7 to give a new proof of $%
(i)\Rightarrow (ii)$ of Theorem 2.4.12.\textbf{\medskip }

The following definition is motivated by the work of Albert Visser \cite%
{Albert-Tehran}; it was introduced in \cite{Ali-Visserian}.\textbf{\medskip }

\noindent \textbf{5.7.~Definition.~}Suppose $T$ is a first order theory. $T$
is \textit{solid} iff the following property $(\nabla )$ holds for all
models $\mathcal{M}$, $\mathcal{M}^{\ast },$ and $\mathcal{N}$ of $T$%
:\medskip

\noindent $(\nabla )$\qquad If $\mathcal{M}\trianglerighteq _{\mathrm{par}}%
\mathcal{N}\trianglerighteq _{\mathrm{par}}\mathcal{M}^{\ast }$ and there is
a parametrically $\mathcal{M}$-definable isomorphism\textit{\ }$i_{0}:%
\mathcal{M}\rightarrow \mathcal{M}^{\ast }$, then there is a parametrically $%
\mathcal{M}$-definable isomorphism\textit{\ }$i:\mathcal{M}\rightarrow
\mathcal{N}$.\medskip

Visser showed that PA is a solid theory, a result that was extended to ZF
and Kelley-Morse theory of classes in \cite{Ali-Visserian}. An examination
of the proof of solidity of \textrm{ZF} presented in \cite{Ali-Visserian}
shows a slightly more general result that plays a crucial role in the proof
of Theorem 5.9 below, namely:\medskip

\noindent \textbf{5.8.~Theorem.~}\textit{Suppose} $\mathcal{M}$ \textit{and}
$\mathcal{M}^{\ast }$ \textit{are models of} \textrm{ZF}, \textit{and} $%
\left( \mathcal{N},X\right) \models \mathrm{ZF(}X\mathrm{)}$. \textit{Then} $%
(\nabla ^{+})$ \textit{below holds}:\medskip

\noindent $(\nabla ^{+})$\qquad \textit{If }$\mathcal{M}\trianglerighteq _{%
\mathrm{par}}\left( \mathcal{N},X\right) \trianglerighteq _{\mathrm{par}}%
\mathcal{M}^{\ast }$ \textit{and there is a parametrically }$\mathcal{M}$-%
\textit{definable isomorphism\ }$i_{0}:\mathcal{M}\rightarrow \mathcal{M}%
^{\ast }$, \textit{then there is a parametrically }$\mathcal{M}$-\textit{%
definable isomorphism\ }$i:\mathcal{M}\rightarrow \mathcal{N}$.\medskip

The following general result shows that in contrast with Theorem 5.5(b), the
theories $\mathrm{ZFC}+\Lambda $ and $\mathrm{ZFI}_{\mathrm{<}}$\textit{\ }%
are not bi-interpretable\textit{.}\medskip

\noindent \textbf{5.9.~Theorem.~}\textit{No model of} $\mathrm{ZF}$\textit{\
is parametrically bi-interpretable with a model of} $\mathrm{ZFI}$. \medskip

\noindent \textbf{Proof.}~Suppose to the contrary that there are
interpretations $\mathcal{I}$ and $\mathcal{J}$\ that witness that some
model $\mathcal{M}$ of ZFC is parametrically bi-interpretable with a model
of $\mathrm{ZFI}$. Then by Theorem 5.8, $\mathcal{M}$ can parametrically
define a class $I$ of indiscernibles for itself. But this contradicts
Corollary 3.5.\hfill $\square $

\bigskip

\begin{center}
\textbf{6.~SOME\ VARIANTS\ OF }$\mathbf{ZFI}_{<}$\textbf{\bigskip }
\end{center}

In this section we discuss four variants of \textrm{ZFI}$_{<}$. We begin
with presenting two of these variants that turn out to be conservative over
\textrm{ZFC}. The first such system \textrm{ZFI}$_{<}^{\circ}$ below can be
intuitively thought of as weakening the stipulation in \textrm{ZFC}$_{<}$
that there is a proper class of indiscernibles over the universe to the
stipulation that there are arbitrarily large sets of indiscernibles over the
universe.\medskip

\noindent \textbf{6.1.~Definition.~}$\mathrm{ZFI}_{<}^{\circ}$ is a theory
formulated in the language $\mathcal{L}_{\mathrm{Set}}\cup \{<,I(x,y)\},$
where $I(x,y)$ is a binary predicate, whose axioms consist of the following
three groups of axioms:

\begin{itemize}
\item We will write $I(x,\alpha )$ as $x\in I_{\alpha }$ for better
readability. \medskip
\end{itemize}

\noindent $(1)$ $\mathrm{ZF}(<,I)+\mathrm{GW}$.\medskip

\noindent $(2)$ The conjunction of $\forall \alpha \in \mathrm{Ord\ }\forall
x\mathrm{(}x\in I_{\alpha }\rightarrow \left( x\in \mathrm{Ord}\wedge \alpha
\in \mathrm{Ord}\right) \mathrm{)}$ with $\forall \alpha \in \mathrm{Ord\ }%
\left\vert \{x:x\in I_{\alpha }\}\right\vert \geq \aleph _{\alpha }.$\medskip

\noindent $(3)$ A scheme consisting of sentences of the form $\forall \alpha
\in \mathrm{Ord}\left( \mathrm{Indis}_{\varphi }(I_{\alpha })\right) $, for
each formula $\varphi $ in the language $\mathcal{L}_{\mathrm{Set}}(\mathrm{%
<)}$. This scheme ensures that $\left( I_{\alpha },\in \right) $ is a set of
order indiscernibles for the ambient model $(\mathrm{V},\in )$ of set theory
for each ordinal $\alpha \mathrm{.}$ More explicitly, if $\varphi =\varphi
(v_{1},\cdot \cdot \cdot ,v_{n}),$ then $\mathrm{Indis}_{\varphi }(I_{\alpha
})$ is the formula below:\medskip

\begin{center}
$\forall x_{1}\in I_{\alpha }\cdot \cdot \cdot \forall x_{n}\in I_{\alpha }$
$\forall y_{1}\in I_{\alpha }\cdot \cdot \cdot \ \forall y_{n}\in I_{\alpha
} $ \medskip

$[(x_{1}\in \cdot \cdot \cdot \in x_{n})\wedge (y_{1}\in \cdot \cdot \cdot
\in y_{n})\rightarrow (\varphi (x_{1},\cdot \cdot \cdot
,x_{n})\leftrightarrow \varphi (y_{1},\cdot \cdot \cdot ,y_{n}))].$\medskip
\end{center}

\begin{itemize}
\item Thus $\mathrm{ZFI}_{<}^{\circ }$ is a theory that ensures that for each
ambient infinite cardinal $\aleph _{\alpha }$, there is a set of
indiscernibles for $(\mathrm{V},\in ,<)$ of size at least $\aleph _{\alpha }$%
.
\end{itemize}

\noindent \textbf{6.2.~Theorem.~}$\mathrm{ZFI}_{<}^{\circ}$ \textit{is a}
\textit{conservative extension of }\textrm{ZFC.}\medskip

\noindent \textbf{Proof.} To show the conservativity of $\mathrm{ZFI}%
_{<}^{\circ}$ \ over $\mathrm{ZFC}$, it suffices to show that every
countable model $\mathcal{M}$\ of $\mathrm{ZFC}$ has an elementary extension
to a model $\mathcal{M}^{*}$ which has an expansion to \textrm{ZFI}$%
_{<}^{\circ }.$ So let $\mathcal{M}$ be a countable model of $\mathrm{ZFC}$.
By Theorem 2.1.3, there is an expansion $(\mathcal{M},<_{M})$ of $\mathcal{M}
$ that satisfies $\mathrm{ZF}(<)+\ \mathrm{GW}$. By compactness, to show the existence of
the desired elementary extension $\mathcal{M}^{*}$ of $\mathcal{M}$, it suffices to show
that the elementary diagram of $\mathcal{M}$ is consistent with $\mathrm{ZFCI%
}_{<}^{\circ }.$ Towards this goal, fix some list $\left\langle \varphi
_{i}(x_{1},\cdot \cdot \cdot ,x_{i}):i<\omega \right\rangle $ of $\mathcal{L}%
_{\mathrm{Set}}(\mathrm{<)}$-formulae such that the free variables of $%
\varphi _{i}$ are among $x_{1},\cdot \cdot \cdot ,x_{i}$, and let:

\begin{center}
$T:=\mathrm{Th}(\mathcal{M},<_{M},m)_{m\in M}\cup \mathrm{ZFI}_{<}^{\circ},$
\end{center}

\noindent and let $T_{0}$ be a finite subset of $T$. Then there is some $%
j\in \omega $ such that if an axiom of the form $\forall \alpha \in \mathrm{%
Ord}\left( \mathrm{Indis}_{\varphi }(I_{\alpha })\right) $ is included in $%
T_{0},$ then $\varphi $ is among $\left\{ \varphi _{0},\cdot \cdot \cdot
,\varphi _{j}\right\} .$ Recall that by the classical Erd\H{o}s-Rado theorem
\cite[Theorem 7.3]{Kanamori book}, ZFC proves:

\begin{center}
$\beth _{n}(\kappa )\rightarrow (\kappa ^{+})_{\kappa }^{n+1}$ for every
infinite cardinal $\kappa $ and every $n\in \omega ,$
\end{center}

\noindent where $\beth _{n}(\kappa )$ is the Beth function, defined by: $%
\beth _{0}(\kappa )=$ $\kappa $ and $\beth _{n+1}=2^{\beth _{n}(\kappa )}.$
The Erd\H{o}s-Rado theorem, together with a global well-ordering $<_{M}$,
then allows us to define within $(\mathcal{M},<_{M})$ a function $F:\mathrm{%
Ord}\rightarrow \mathrm{V}$ such that for each $\alpha \in \mathrm{Ord,}$ $%
F(\alpha )$ is a set of ordinals of cardinality at least $\aleph _{\alpha },$
and $\left( F(\alpha ),\in \right) $ is homogeneous for $\left\{ \varphi
_{0},\cdot \cdot \cdot ,\varphi _{j}\right\} .$ This makes it clear that if $%
I_{M}(x,y)$ is defined in $(\mathcal{M},<_{M})$ as $\left[ y\in \mathrm{Ord}%
\wedge x\in F(y)\right] ,$ then:

\begin{center}
$\left( \mathcal{M},<_{M},I_{M}\right) \models T_{0}.$
\end{center}

\noindent Thus every finite subset of $T$, and therefore $T$ itself, is
consistent, as promised.\hfill $\square \medskip $

The second variant of \textrm{ZFI}$_{<}$ we consider, denoted \textrm{ZFI}$%
_{<}^{-}$ is obtained from $\mathrm{ZFI}_{<}\ $by weakening the demand that $%
I$ is amenable to the demand that it satisfies $\mathrm{Sep}(\mathcal{L}_{%
\mathrm{Set}}(<,I))$. Note that \textrm{ZFI}$_{<}^{-}$ does include $\mathrm{%
Coll}(\mathcal{L}_{\mathrm{Set}}).\medskip $

\noindent \textbf{6.3.~Definition. }Let \textrm{ZFI}$_{<}^{-}$ be the
subsystem of \textrm{ZFI}$_{<}$ whose axioms consist of the following:$%
\medskip $

\noindent $(1)$ $\mathrm{ZF(<)}+\mathrm{GW}+$ $\mathrm{Sep}(<,I).\medskip $

\noindent $(2)$ The sentence expressing that $I$ is cofinal in $\mathrm{Ord}$%
.$\medskip $

\noindent $(3)$ The scheme \textrm{Indis}$_{\mathcal{L}_{\mathrm{Set}}}(I)$
(as in Definition 3.1).\medskip

\noindent \textbf{6.4.~Theorem.~}$\mathrm{ZFI}_{<}^{-}$ \textit{is a}
\textit{conservative extension of }\textrm{ZFC.}\medskip

\noindent \textbf{Proof.~}It suffices to show that every countable model $%
\mathcal{M}$\ of $\mathrm{ZFC}$ has an elementary extension to a model $\mathcal{M}%
^{\ast }$ which has an expansion to \textrm{ZFI}$_{<}^{-}.$ So let $\mathcal{%
M}$ be a countable model of $\mathrm{ZFC}$, and $(\mathcal{M},<_{M})$ be an
expansion of $\mathcal{M}$ that satisfies $\mathrm{ZF}(<)+\ \mathrm{GW.}$
Then let:

\begin{center}
$T:=\mathrm{Th}(\mathcal{M},<_{M},m)_{m\in M}\cup \mathrm{ZFI}_{<}^{-},$
\end{center}

\noindent and let $T_{0}$ be a finite subset of $T$. Let $n\in \omega $ be
large enough so that any sentence in $T_{0}$ that belongs to the
elementary diagram of $\mathcal{M}$ is $\Sigma _{n}$, and let $j\in \omega $
be large enough so that if the sentence $\mathrm{Indis}_{\varphi }(I)$ is in
$T_{0},$ then $\varphi $ is among $\left\{ \varphi _{0},\cdot \cdot \cdot
,\varphi _{j}\right\} .$ By the reflection theorem for $\mathrm{ZF}(<)$,
there is some $\alpha \in \mathrm{Ord}^{M}$ such that:

\begin{center}
$\left( \mathcal{M}(\alpha ),\ <_{M(\alpha )}\right) \prec _{\Sigma
_{n}}\left( \mathcal{M},<_{M}\right) $ and $\mathcal{M}\models \mathrm{cf}%
(\alpha )=\omega ,$
\end{center}

\noindent where $<_{M(\alpha )}$ is the restriction of $<_{M}$ to $M(\alpha
).$ Since $\alpha $ has countable cofinality in $\mathcal{M}$, by the
infinite Ramsey Theorem applied within $\mathcal{M}$, there is a cofinal
subset $I$\ of $\alpha $ that is $\varphi _{i}$-indiscernible for each $%
i\leq j.$ Since any expansion of $\mathcal{M}(\alpha )$ within $\mathcal{M}$
satisfies Zermelo set theory in the extended language, this makes it clear
that:

\begin{center}
$\left( \mathcal{M}(\alpha ),<_{M(\alpha )},I\right) \models T_{0},$
\end{center}

\noindent which completes the proof of consistency of $T$.\hfill $\square
\medskip $

Finally, in the remarks below, we briefly discuss two natural strengthenings of $\mathrm{ZFI}_{<}$ whose purely set-theoretical consequences go beyond $\mathrm{ZFC} + \Lambda$.

\medskip

\noindent \textbf{6.5.~Remark.~}Recall the classical fact of large cardinal
theory that the Silver indiscernibles (of the constructible universe) are
closed and unbounded in the ordinals, and satisfy the so-called
remarkability condition \cite[Lemma 9.10]{Kanamori book}. A moment's
reflection reveals that the axiom \textquotedblleft $I$ is closed and
unbounded in\textrm{\ Ord}\textquotedblright\ is inconsistent with $\mathrm{%
ZFI}$ based on cofinality considerations and indiscernibility. More specifically, $%
\mathrm{ZFI}$ implies that either all limit ordinals in $I$ have cofinality $\omega $, or they are all of uncountable cofinality; each
of which is inconsistent with $I$ being closed and unbounded. On the other
hand, if $(\mathcal{M},<, I)\models \mathrm{ZFI}_{<}$, using a class-theoretic adaptation of
Baumgartner's characterization of $n$-ineffable cardinals in terms of
regressive partition relations \cite{Baumgartner (ineffable)}, the
remarkability condition of $I$ can be recast as asserting that $I$ is \textquotedblleft
definably stationary\textquotedblright\ in $\mathcal{M}$, i.e., $I$
intersects every closed unbounded subset of $\mathrm{Ord}^{\mathcal{M}}$
that is parametrically definable in $(\mathcal{M},<)$. Using this equivalence one can readily show that the remarkability condition can be
consistently added to $\mathrm{ZFI}_{<}$, assuming that $\mathrm{ZFC}$ +
\textquotedblleft there is a cardinal $\kappa $ that is $n$-ineffable for
each $n\in \omega $\textquotedblright\ is consistent. Moreover, the techniques of this paper can be extended to show that $%
\mathcal{L}_{\mathrm{Set}}$-consequences of the strengthening of $\mathrm{ZFI%
}_{<}$ by an axiom scheme expressing the remarkability of $I$ turn out to
coincide with the theorems of $\mathrm{ZFC}$ + $\Theta $ where $\Theta
=\{\theta _{n}:n\in \omega \}$ and $\theta _{n}$ is the $\mathcal{L}_{%
\mathrm{Set}}$-sentence asserting that there is an $n$-ineffable cardinal $%
\kappa $ such that $\mathrm{V}(\kappa )$ is a $\Sigma _{n}$-elementary
submodel of the universe $\mathrm{V}$. Another axiomatization for $\Theta ,$
in the presence of ZFC, is $\Theta ^{\prime }=\{\theta _{n}^{\prime }:n\in
\omega \}$, where $\theta _{n}^{\prime }$ is the $\mathcal{L}_{\mathrm{Set}}$%
-sentence asserting that there is an $n$-subtle cardinal $\kappa $ such that
$\mathrm{V}(\kappa )$ is a $\Sigma _{n}$-elementary submodel of the universe
$\mathrm{V}$. It is worth mentioning that the proof of 2.4.12 can be
modified to show that ZFC + $\Theta $ axiomatizes the purely set-theoretical
consequences of the class theory $\mathrm{GBC}$ + \{\textrm{Ord} is $n$%
-ineffable: $n<\omega \}$; this class theory can also be axiomatized by $%
\mathrm{GBC}$ + \{\textrm{Ord} is $n$-subtle: $n<\omega \}.$
\medskip

\pagebreak

\noindent \textbf{6.6.~Remark.~}Recall from part $(iii)$ of Theorem 4.1 that $\mathrm{ZFI}^{\mathrm{Good}}_{<}$ is the result of augmenting the theory $\mathrm{ZFI}_{<}$ with the sentence $\psi$ that expresses:
\begin{center}
$\forall \alpha, \beta\in I [\alpha \in \beta \rightarrow (\mathrm{V}(\alpha),\in,<) \prec (\mathrm{V}(\beta),\in,<)].$
\end{center}

\noindent Note that thanks to $\psi$ within $\mathrm{ZFI}^{\mathrm{Good}}_{<}$ there is a definable \emph{full} satisfaction predicate $S$ for $(\mathrm{V},\in,<)$, since $S$ can simply be defined as the union of the Tarskian satisfaction predicates for structures of the form $(\mathrm{V}(\alpha),\in,<)$ as $\alpha$ ranges in $I$. Using $S$ we can formulate the following axiom $\sigma$ that expresses that elements of $I$ are order indiscernible in the sense of $S$ for the structure $(\mathrm{V},\in,<)$:

$$\forall n\in \omega ~\forall\varphi(\overline{x})\in \mathrm{Form}_{n}(\mathcal{L}_{\mathrm{Set}}(<))~\forall \overline{i} \in [I]^{n} ~\forall \overline{j}\in [I]^{n}~ (\varphi(\overline{i})\in S) \leftrightarrow (\varphi(\overline{j})\in S).$$

\noindent  It is easy to see that if $\kappa$ is a Ramsey cardinal, then $(\mathrm{V}(\kappa),\in)$ has an expansion to $\mathrm{ZFI}^{\mathrm{Good}}_{<}+\sigma$. In contrast to the consistency of $\mathrm{ZFI}^{\mathrm{Good}}_{<}+\mathrm{V} = \mathrm{L}$, $\mathrm{V} \neq \mathrm{L}$ is provable in $\mathrm{ZFI}^{\mathrm{Good}}_{<}+\sigma$. Indeed $\mathrm{ZFI}^{\mathrm{Good}}_{<}+\sigma$ proves the existence of $0^{\#}$, as well as the  existence of a proper class of almost Ramsey cardinals in the core model. This is because the proof strategy of Theorem 3.8(c,d,e) can be carried out \emph{within} $\mathrm{ZFI}^{\mathrm{Good}}_{<}+\sigma $ to obtain a nontrivial elementary self-embedding of $(\mathrm{V}_I,\in,<)$ by shifting the indiscernibles, where $\mathrm{V}_I$ is the proper class consisting of sets that are definable (in the sense of $S$) in $(\mathrm{V},\in,<)$ with parameters from $I$.

\medskip

 It is not clear whether the purely set-theoretical consequences of $\mathrm{ZFI}^{\mathrm{Good}}_{<}+\sigma$ can be axiomatized by a `natural' extension of ZFC + $\Lambda$. However, it is noteworthy that the purely set-theoretical consequences of $\mathrm{ZFI}^{\mathrm{Good}}_{<}+\sigma$ coincide with the purely set-theoretical consequences of $\mathrm{GBC}$ + ``Ord is Ramsey'', where ``Ord is Ramsey'' expresses: $\mathrm{Ord}\rightarrow \mathrm{Ord}^{<\omega}_{2}$. This follows from the fact that a model $\mathcal{M}$ of ZFC has an expansion to $\mathrm{ZFI}^{\mathrm{Good}}_{<}+\sigma$ iff $\mathcal{M}$ has an expansion to $\mathrm{GBC}$ + ``Ord is Ramsey''. The right-to-left direction of this equivalence is handled by a routine argument; the left-to-right direction is established by noting that if $(\mathcal{M},<_M,I)$ is a model of $\mathrm{ZFI}^{\mathrm{Good}}_{<}+\sigma$, and $\mathfrak{X}$ is chosen as the collection of subsets of $M$ that are parametrically definable in $(\mathcal{M},<_M)$ in the sense of the aforementioned $(\mathcal{M},<_M)$-definable full satisfaction predicate $S$, then $(\mathcal{M},\mathfrak{X})$ satisfies $\mathrm{GBC}$ + ``Ord is Ramsey''.
\medskip

\begin{center}
\textbf{7.~OPEN QUESTIONS\bigskip }
\end{center}

Here we draw attention to some natural questions that arise from the results
of the paper. $\medskip $

\noindent \textbf{7.1.~Question. }\textit{Does} $\mathrm{ZFCI}\vdash \Lambda
$?\medskip

\begin{itemize}
\item One would expect that by the use of a generic global well-ordering one
could show that $\mathrm{ZFI}_{<}$ is a conservative extension of $\mathrm{%
ZFCI}$, but our attempts in this direction have been unsuccessful. See also Remark 4.3.
\end{itemize}

\noindent \textbf{7.2.~Question. }Can Theorem 4.9 be improved by weakening
the statement $\theta $ of that theorem to the statement $\theta
^{-}=\forall n\in \omega \left( \mathrm{Ord}\rightarrow \left( \mathrm{Ord}%
\right) _{2}^{n}\right) $?

\begin{itemize}
\item We conjecture that the answer to Question 7.2 is in the positive, in
analogy with the well-known unprovability of the statement $\psi =\forall n\in \omega
\left( \omega \rightarrow \left( \omega \right) _{2}^{n}\right) $ in \textrm{%
ACA}$_{0}$.
\end{itemize}

\noindent \textbf{7.3.~Question. }\textit{Let} $T\vdash _{\pi }\varphi $
\textit{indicate that} $\pi $ \textit{is the} (\textit{binary code of}) \textit{a proof of} $\varphi $ \textit{from axioms
in the theory} $T$.\textit{ Is there a polynomial-time
computable function} $f$ \textit{such that for all} $\mathcal{L}_{\mathrm{Set%
}}$\textit{-sentences }$\varphi $, \textit{the following holds}:

\begin{center}
\textrm{ZFI}$_{<}\vdash _{\pi }\varphi \Rightarrow $ $\mathrm{ZFC}+\Lambda
\vdash _{f(\pi )}\varphi ?$
\end{center}

\begin{itemize}
\item We suspect that Question 7.3 has a positive  answer.
\end{itemize}

\noindent \textsc{Department of Philosophy, Linguistics, and the Theory of
Science \newline
\noindent University of Gothenburg, Gothenburg, Sweden}\newline
\noindent \texttt{email: ali.enayat@gu.se}

\end{document}